\documentclass[11pt, reqno, oneside]{amsart}

 \usepackage{curves} 
\usepackage{graphicx} 
\usepackage{amsmath, amsfonts, amssymb, amscd,graphics, graphpap, tikz, color, xcolor, cancel,enumerate,enumitem, hyperref, mathabx, mathrsfs, mathtools}
\usepackage[scr=rsfs,cal=boondox]{mathalfa}

\usepackage{geometry}
\usetikzlibrary{arrows.meta,arrows}
\usetikzlibrary{decorations.pathreplacing,decorations.markings}
\tikzset{arrow data/.style 2 args={%
      decoration={%
         markings,
         mark=at position #1 with \arrow{#2}},
         postaction=decorate}
      }%
      
\usepackage[mathscr]{euscript}
  \usepackage[normalem]{ulem}
  \usepackage{soul}

\numberwithin{equation}{section}

\theoremstyle{plain}
\newtheorem{theo}{Theorem}[section]
\newtheorem{lem}[theo]{Lemma}
\newtheorem{prop}[theo]{Proposition}

\newtheorem{cor}[theo]{Corollary}

\theoremstyle{definition}
\newtheorem{rem}[theo]{Remark}

\newtheorem{example}[theo]{Example}
\newtheorem{definition}[theo]{Definition}

\newenvironment{pf}{\noindent{\it Proof. }}{$\hfill\square$\par\medskip}

\theoremstyle{plain}

\theoremstyle{definition}

\newcommand{\beq}{\begin{equation}}
\newcommand{\eeq}{\end{equation}}
\newcommand{\beqn}{\begin{equation*}}
\newcommand{\eeqn}{\end{equation*}}
\renewcommand{\a}{\alpha}
\renewcommand{\b}{\beta}

\renewcommand{\d}{\delta}

\newcommand{\f}{\varphi}

\newcommand{\h}{\eta}

\renewcommand{\k}{\kappa}

\renewcommand{\o}{\omega}
\newcommand{\q}{\vartheta}
\renewcommand{\r}{\rho}
\newcommand{\s}{\sigma}

\newcommand{\x}{\xi}
\newcommand{\z}{\zeta}
\newcommand{\D}{\Delta}

\newcommand{\G}{\Gamma}

\renewcommand{\L}{\Lambda}

\newcommand{\bC}{\mathbb{C}}
\newcommand{\bD}{\mathbb{D}}
\newcommand{\bR}{\mathbb{R}}

\newcommand{\bT}{\mathbb{T}}

\newcommand{\gd}{\mathfrak{d}}

\newcommand{\gr}{\mathfrak{r}}

\newcommand\SO{\mathrm{SO}}

\newcommand{\cA}{\mathscr{A}}
\newcommand{\cB}{\mathscr{B}}
\newcommand{\cC}{\mathcal{C}}

\newcommand{\cF}{\mathscr{F}}

\newcommand{\cH}{\mathscr{H}}

\newcommand{\cJ}{\mathscr{J}}
\newcommand{\cK}{\mathscr{K}}
\newcommand{\cL}{\mathscr{L}}
\newcommand{\cM}{\mathscr{M}}

\newcommand{\cQ}{\mathscr{Q}}

\newcommand{\cS}{\mathscr{S}}

\newcommand{\cU}{\mathscr{U}}
\newcommand{\cV}{\mathscr{V}}
\newcommand{\cW}{\mathscr{W}}

\newcommand{\p}{\partial}

\renewcommand{\square}{\kern1pt\vbox
{\hrule height 0.6pt\hbox{\vrule width 0.6pt\hskip 3pt
\vbox{\vskip 6pt}\hskip 3pt\vrule width 0.6pt}\hrule height0.6pt}\kern1pt}

\DeclareMathOperator\sign{sign}

\renewcommand\Re{\operatorname{Re}}

\renewcommand\={:=}

\newcommand{\wt}{\widetilde}
\newcommand{\wh}{\widehat}

\newcommand{\bt}{\begin{theo}}
\newcommand{\et}{\end{theo}}
\newcommand{\bp}{\begin{prop}}
\newcommand{\ep}{\end{prop}}
\newcommand{\bc}{\begin{cor}\ \ }
\newcommand{\ec}{\end{cor}}
\newcommand{\bl}{\begin{lem}\ \ }
\newcommand{\el}{\end{lem}}
\newcommand{\bd}{\begin{definition}}
\newcommand{\ed}{\end{definition}}
\newcommand{\n}{\nabla}

\newcommand{\be}{\begin{equation}}
\newcommand{\ee}{\end{equation}}

\def\<#1,#2>{\langle\,#1,\,#2\,\rangle}
\newcommand{\arr}{\begin{array}{rlll}}
\newcommand{\ea}{\end{array}}
\newcommand{\bea}{\begin{eqnarray}}
\newcommand{\eea}{\end{eqnarray}}
\newcommand{\bean}{\begin{eqnarray*}}
\newcommand{\eean}{\end{eqnarray*}}

\newcommand\wc{\widecheck}

\newcommand{\Ga}[3]{\G_{#1 #2}^{\phantom{#1} #3}}
\newcommand{\GGa}[3]{{\mathbf \G}_{#1 #2}^{\phantom{#1} #3}}
\newcommand{\dGa}[3]{\d{\G}_{#1 #2}^{\phantom{#1} #3}}

\newcommand{\HH}[3]{{\mathbf H}_{#1 #2}^{\phantom{#1} #3}}

\newcommand{\RR}[4]{\operatorname{R}_{#1 #2 #3}^{\phantom{#1 #2 #3} #4}}

\font\smallsmc = cmcsc9
\font\smalltt = cmtt8
\font\smallit = cmti8

\def\sideremark#1{\ifvmode\leavevmode\fi\vadjust{
\vbox to0pt{\hbox to 0pt{\hskip\hsize\hskip1em
\vbox{\hsize2cm\tiny\raggedright\pretolerance10000
\noindent #1\hfill}\hss}\vbox to8pt{\vfil}\vss}}}

\renewcommand{\sf}{shearfree }

\renewcommand{\Re}{\operatorname{Re}}

\newcommand{\cappa}{{\bf k}}

\newcommand{\ps}{{\operatorname{p}}}
\newcommand{\qs}{{\operatorname{q}}}

\newcommand{\Ric}{{\operatorname{Ric}}}
\newcommand\bex{\begin{example}}
\newcommand\eex{\end{example}}
\newcommand\br{\begin{rem}}
\newcommand\er{\end{rem}}

\makeatletter
\@namedef{subjclassname@1991}{2020 Mathematics Subject Classification}
\makeatother

\title[Einstein manifolds  with  optical geometries of Kerr type]{ Einstein  manifolds\\ with 
   optical geometries of Kerr type}
\author{Masoud Ganji}
\author{Cristina Giannotti}
\author{Gerd Schmalz}
\author{Andrea Spiro}

\begin{document}
\setcounter{section}0

\begin{abstract}We classify the Ricci flat Lorentzian $n$-manifolds satisfying    three particular conditions,    encoding and combining some    crucial features of  the  Kerr metrics and the Robinson-Trautman optical structures.  We prove that: (a)   If  $n  > 4$, there is no Lorentzian manifold satisfying the considered Kerr type conditions,  in unexpected contrast with what occurs for the metrics satisfying   (very similar) Taub-NUT type conditions;  (b)   If $n = 4$   there are  two   large classes of  such Kerr type manifolds. Each  class consists of   manifolds  fibering over  open Riemann surfaces, equipped  with  a metric  of  constant Gaussian curvature  $\k = 1$ or $\k = -1$.  The    first class    includes a three parameter family of  metrics  admitting  real analytic  extensions to   $(\bR^3 \setminus\{0\}) \times \bR = (S^2 \times \bR_+) \times \bR$  and a large class of  other metrics  not admitting  this kind of  extensions. The metrics of this first class  admitting such extensions  are all isometric to  the  
well known    Kerr metrics, with  the three parameters corresponding to the  three space-like components of the angular momentum of the gravitational field.
The second class contains a subclass of metrics defined on $\big(\bD\times \bR_+\big)\times \bR$, where $\bD$ is the Lobachevsky Poincar\'e disc. This subclass is in bijection with the holomorphic functions on  $\bD$ satisfying an appropriate open condition.
 These and other results  are   consequences of a very simple  way to  construct totally   explicit examples of   Ricci flat Lorentzian manifolds.  
\end{abstract}
\maketitle
\section{Introduction}
The main result of this paper is the classification of the  even dimensional Lorentzian manifolds, admitting families of Ricci flat metrics that are compatible with a  twisting optical structure of Robinson-Trautman type  and satisfying certain  constraints,  which encode some of the most important features of the  $4$-dimensional  Kerr metrics. The results and their  proofs  have two noticeable  consequences: (1) a no-go theorem for higher dimensional Einstein manifolds with twisting optical structures  and (2) a simple method of construction of  a very  large class of  explicit  Ricci flat $4$-dimensional Lorentzian manifolds, parameterised by the solutions of a   linear elliptic p.d.e. in two variables.\par
\smallskip
Before going into the details of our results,  we  need to recall a few  known facts on Kerr metrics. These  metrics are  Lorentzian  metrics,  depending on real  parameters  $m$ and $a$, which are defined on   the $4$-dimensional manifold  $M^{\text{Kerr}} = (\bR^3 \setminus\{0\}) \times \bR =  (S^2 \times \bR_+) \times \bR $ and having the following coordinate expressions: 
\beq \label{standardKerr}
\begin{split} g^{\text{Kerr}} &= - \frac{1}{2} (\r^2 + a^2 \cos^2 \x) \big(d \x^2 + \sin^2 \x d \psi^2\big) - \\
&-  \big(d v + a \sin^2 \x d \psi\big) \vee \big(d\r +  a \sin^2 \x d \psi\big) 
+ \frac{1}{2} \left(1 - \frac{2 m \r}{\r^2 + a^2 \cos^2 \x}\right) (d v + a \sin^2 \x d \psi)^2 
\ .
 \end{split}
\eeq
Here    $\psi$,   $\xi$  are the  usual spherical coordinates of $S^2$,   $\r$, $v$ are  the  standard Cartesian coordinates of   $(0, + \infty) \times \bR$ and $m, a$ are  the parameters which are  physically interpreted as  the    {\it Schwarzschild mass}   of the gravitational field, and the ratio $a = \cJ/m$ between the  {\it angular momentum} $\cJ \neq 0$  of the  gravitational field  and the mass $m$  (see e.g. \cite{Ke2} (\footnote{ Note that   \eqref{standardKerr} is  actually  a metric which is  a  re-scaling  of  the  one in   \cite{Ke}  by the factor $- \frac{1}{2}$. It is  a normalisation that  better fits with  our purposes.  Notice  also that  we  are considering the  signature   $(-,-,-,+)$  instead  of  the signature  $(+, +, +, -)$  used  in \cite{Ke}.})). \par
\smallskip
There is however an alternative and very useful    presentation of     $g^{\text{Kerr}}$  or, more precisely,  of  its  restriction  to the open subset  
$$M^{\text{Kerr}}_{(+)}  = (S^2_{(+)} \times \bR_+)\times \bR  \subset M^{\text{Kerr}}\!\!\!,\  S^2_{(+)}  \= \text{upper hemisphere} =    \left\{ (\xi, \psi)\ :  \xi \in \left[0, \frac{\pi}{2}\right) \right\} \ ,$$
 that makes  manifest that $g^{\text{Kerr}}$  is  compatible with    a particular      {\it twisting optical structure}. \par
\smallskip
The   notion of optical structure was  introduced by Robinson and Trautman  to study    Lorentzian manifolds  admitting   electromagnetic plane waves (or   appropriate higher dimensional generalisations)  propagating  along  a prescribed foliation by   null  geodesics(see \cite{RT,RT1, Tr} -- see also  \cite{FLT, AGS, AGSS, AGSS1} and  references therein for  coordinate-free  presentations  and discussions of  the equivalences among the diverse  presentations of the same  objects).  In detail,  an {\it  optical structure}
on a  manifold  $M$ is a  pair $\cQ = (\cK, \{g\})$ given by  
\begin{itemize}[leftmargin=15pt]
\item[--] a $1$-dimensional distribution $\cK \subset TM$,  
\item[--]  the  collection $\{g\}$ of  all Lorentzian metrics $g$,   such that: (1) the  curves that are  tangent to  $\cK$ determine   a  {\it  geodesic \sf null congruence} for $M$;  (2)  the  $g$-orthogonal distribution $\cK^{\perp_g}$   is a  fixed   distribution $\cW = \cK^{\perp_g} \subset TM$,  independent of $g$  (\footnote{A {\it geodesic \sf congruence} is   a foliation   by null geodesics that  are tangent to a  vector field $\ps_o$, whose flow  leaves invariant   the $g$-orthogonal   distribution $\cW = \cK^{\perp_g}$ and the  conformal class of    $g|_{\cW \times \cW}$.}). 
\end{itemize} 
The metrics in  the class  $\{g\}$  are said to be    {\it compatible}   with the optical structure $\cQ$. \par
\smallskip
An important  class  of   optical structures is given by those on a special class of bundles, which we now very briefly review (they are extensively discussed in  \cite{AGSS},  where they are called {\it Robinson-Trautman \sf structures}).  Consider   a $2k$-dimensional  {\it quantisable}  K\"ahler manifold  $(N, J, g_o)$.  This means  that there exists:  (a) a principal  bundle   $\pi: \cS \to N$, with $1$-dimensional structure group $A = \bR_+ \simeq \bR$  or $S^1$ and (b)  a connection $1$-form $\theta: T \cS \to  \operatorname{Lie}(A) = \bR$   such that    $d \theta = \pi^*\o_o$, where  $\o_o \= g_o(J, \cdot)$  is the K\"ahler form of $(N, J, g_o)$.   Then take the trivial  principal $\bR$-bundle $\pi^{\cS}: M = \cS \times \bR \to \cS$ over $\cS$ and denote by  $\q$ and $\ps_o$    the   pulled-back $1$-form   $\q = \pi^{\cS*} \theta$  on $M$ and the vertical vector field  $\ps_o = \frac{\p}{\p u}$  generating the   right $\bR$-action on $M = \cS \times \bR$, respectively.   
It is known  (see \cite{AGSS, AGSS1} and \S \ref{sect2} of this paper) that the pair $\cQ = (\cK, \{g\})$ given by the  $1$-dimensional distribution $\cK = \langle \ps_o \rangle \subset T M$  and the class $\{g\}$ of all Lorentzian metrics on $M$  having the form
\beq\label{comp1}  g = \s (\pi \circ \pi^{\cS})^* g_o + \q \vee \bigg(\ps_o^* + \gd +  \frac{\wt \b}{2} \q\bigg)\ ,\qquad \text{where}\eeq
\begin{itemize}[leftmargin = 15 pt]
\item[(i)]  $\ps_o^*$   is  $\bR$-invariant, satisfies  $ \ps_o^*(\ps_o) = 1$  and is such that   $\cH \= \ker \ps_o^* \cap \ker \q$  is the horizontal distribution in $T \cS \subset TM \simeq T \cS +_M \bR$, 
\item[(ii)]  $\s, \wt \b$  are freely specifiable smooth  functions with  $\s$ nowhere vanishing, and $\gd$ is a freely specifiable $1$-form with $\ker \gd$ containing the vertical distribution of $\pi \circ \pi^\cS: M \to N$, 
\end{itemize}
is an optical structure on $M = \cS \times \bR$. For such an optical structure, the common  $g$-orthogonal distribution $\cW  = \cK^{\perp_g}$  is the kernel distribution  $\cW = \ker \q$ and  it is such that  $\dim \ker d \q|_{\cW \times \cW} = 1$. Note that $1$ is  the minimal possible value  that might occur for  the dimensions of the kernels  $\ker d \a_x|_{\cW \times \cW} $, $x \in M$,  for an arbitrary  $1$-form $\a$ that vanishes identically  on $\cW$.  The optical structures  for which   this property occurs are called  {\it  twisting} (\cite{RT, RT1, AGSS, FLT}).\par   \par
\medskip
Coming back to the Kerr metrics,   we can immediately observe   that:
\begin{itemize}[leftmargin = 20pt]
\item[(a)] $M^{\text{Kerr}}_{(+)}$ is a  trivial $\bR$-bundle over  the   $\bR_+$-bundle $\pi: \cS = S^2_{(+)} \times \bR_+ \to S^{(+)}$;  
\item[(b)]  the restricted Kerr metric  $g^{\text{Kerr}(+)} \= g^{\text{Kerr}}|_{M_{(+)}}$ has  the form  \eqref{comp1}  with  $\gd \equiv 0$, the function $\s$,  $\wt \b$ given by     
\beq \label{113}   \s =  - \frac{\r^2}{4 a \cos \xi} - \frac{a \cos \xi}{4}  \ ,\qquad    \wt \b =  - 1  - \frac{2 m \r}{\r^2 + a^2 \cos^2 \x}\ ,\eeq
and the  K\"ahler metric $g_o$ and the  $1$-forms  $\q$, $\ps_o^*$  given by  
\beq  \label{thego}g_o \=   2 a \cos \xi  (d \xi^2 + \sin^2 \xi d \psi^2)\ ,\qquad  \q \=  d v +    a \sin^2 \xi d \psi\ , \qquad  \ps_o^* \=   dv - d\r\ .  
\eeq
\end{itemize}
It follows  that {\it for any  choice of  $a, m$, the corresponding (restricted)  Kerr metric $g^{\text{\rm Kerr}(+)}$  on    $M^{\text{Kerr}}_{(+)} = (S^2_{(+)} \times \bR_+ )\times \bR$  is compatible with the twisting optical structure $\cQ$ constructed  as above on an $\bR$-bundle over   the quantisable    K\"ahler manifold }
$$\left(\ S^2_{(+)}\ , \ g_o \=    2 a \cos \xi  (d  \xi^2 + \sin^2 \xi d \psi^2)\right)\ .$$
The Kerr metrics \eqref{standardKerr}  are   the  unique real analytic extensions of such $\cQ$-compatible metrics on  $M^{\text{Kerr}}_{(+)} $  to the whole manifold  $M^{\text{Kerr}} =  (\bR^3 \setminus \{0\}) \times \bR$ and are Ricci flat, as is well known. \par
\medskip
At this point it is important to recall that another well-known class of $\cQ$-compatible metrics on    $M^{\text{Kerr}}_{(+)} $  is given by the classical  $4$-dimensional  Taub-NUT metrics. Actually, as it has been shown in  \cite{HLN},  the  classical  Kerr  and Taub-NUT metrics can be locally combined into a  single smoothly parameterised $\cQ$-compatible family of Ricci flat metrics on $M^{\text{Kerr}}_{(+)} $.\par
\smallskip
 Looking for  new  Einstein manifolds  with  optical geometries   and inspired by  the method of  construction   in   \cite{AGSS}  of the higher dimensional versions of the Taub-NUT metrics  \cite{BB, PP, CG,  AC},  we decided to focus on a class  of manifolds with optical structures which naturally includes $M^{\text{Kerr}}_{(+)}$, namely   the connected  Lorentzian manifolds   $M = \cS \times I $,  $I \subset \bR$,   where   $\cS$ is the total space of an $ \bR$-bundle  $\pi: \cS  = N \times \bR \to N$ (\footnote{ We recall that, by a classical result on principal bundles (see \cite[Cor. 12.3]{St1}), each $\bR$-bundle  is  trivial.}) over a  quantisable K\"ahler manifold $(N, J, g_o)$,  equipped  with the optical structure $\cQ$ determined by  the above recipe.  We call  them {\it Kerr manifolds}  (see below, Definition  \ref{Kerrmanifoldsdef}). \par
 \smallskip
  In this paper we tackle the  classification problem  of  the Kerr manifolds  admitting  $1$-parameter families $g^{(m)}$   of  compatible Einstein metrics, parameterised by a real number  $m$ and  for which  the data $\gd$, $\s$, $\wt \b$ appearing in  \eqref{comp1} have the following form
  \begin{itemize}[leftmargin = 18pt]
  \item[(1)]  $\gd \equiv 0$, 
  \item[(2)]   $\wt \b = B + m \wt  \b_o$  for some constant   $B \neq 0$ and a  smooth  $m$-independent function $\wt  \b_o$,  
  \item[(3)]  $\s$  is independent of  $m$ 
  \end{itemize}
  and are constrained by the condition $\bT(\s) = \bT(\wt \b_o) = 0$, with   $\bT = \qs_o - B \ps_o$ where $\qs_o$, $\ps_o$ are null and  vertical  with respect to the projection $\pi\circ\pi^\cS: M \to N$ and  such that $g(\qs_o , \ps_o) = 1 = \q(\qs_o) $. This constraint is equivalent to impose that  each metric $g^{(m)}$  is invariant  with respect to $\bT$, a condition which generalises the time-invariance condition satisfied by the classical Kerr metrics. We call the families of metrics  $g^{(m)}$  of this kind  {\it Kerr families of gravitational fields}.  For each of them, the  Einstein metric $g^{(m = 0)}$  is called {\it background (of Kerr type)}. 
Such a name  is motivated by the  fact that, for any  fixed  $a \neq 0$,   the Kerr metric  corresponding to   $m = 0$ is
\beq \label{15}
\begin{split} g^{\text{Kerr}} &= - \frac{1}{2} (\r^2 + a^2 \cos^2 \x) \big(d \x^2 + \sin^2 \x d \psi^2\big) - \\
&  - \big(d v + a \sin^2 \x d \psi\big) \vee \big(d\r + a \sin^2 \x d \psi\big) 
+ \frac{1}{2} (d v + a \sin^2 \x d \psi)^2 
\ ,
 \end{split}
\eeq
a metric that   can be transformed into  the standard expression of the flat  Minkowski metric under an  appropriate change of coordinates. This is  not   a  coincidence: in Proposition \ref{flatness}  we prove that  {\it any background metric on any $4$-dimensional Kerr manifold  is  flat}, so that  it can be  reasonably  considered   as a metric for  a  background geometry in   a  physical sense.\par
We  stress the fact that the  conditions   that  characterise the Kerr families are very close to those that define the Taub-NUT metrics and their higher dimensional versions:  To pass from one set of conditions to the other   it suffices to set  $B = 0$  and  remove the requirement of linear dependence of $\wt \b$ from a   ``mass'' parameter $m$.  In fact, by  the classification  in \cite[Thm. 5.2]{AGSS}, the  class of the  Einstein  metrics  satisfying these  modified conditions are precisely the above mentioned  generalisations of the  $4$-dimensional Taub-NUT metrics  for arbitrary dimensions and  with arbitrary Einstein constants.\par
 \medskip
In this paper  we succeed   in  classifying  the  Kerr manifolds admitting Kerr families  of {\it Ricci flat} gravitational fields and    explicitly  determine coordinate expressions for such Ricci flat   Kerr families.  It is reasonable to  expect that   a similar classification can be reached by the same tools  also for the case of  non-Ricci flat  Einstein metrics.   The following theorem   summarises our main results.

 \begin{theo} \label{teoremone}   Let $M =  \cS \times I$, $I \subset \bR$,  be a  Kerr  manifold  of dimension $n = 2k + 2\geq 4$. We denote by  $\pi^\cS: M = \cS \times I \to \cS$  the standard projection onto the total space    of  the  $\bR$-bundle    $\pi: \cS = N \times \bR  \to N$  over a K\"ahler manifold $(N, J, g_o)$,  equipped with a connection $1$-form $\theta$   such that $\pi^* \o_o = d \theta$. We also denote by $\q$ the pull-back $\q \= (\pi^{\cS})^* \theta$ on $M$ and by $\cQ_o$ the canonical optical structure on $ \cS \times \bR$ described above. Then: 
 \begin{itemize}[leftmargin = 20pt]
 \item[(i)] If $ n > 4$,  there is  no background and no Kerr family  of gravitational fields on $M$;  
 \item[(ii)] If $n = 4$, there exists a Kerr family $g^{(m)}$ of Ricci flat gravitational fields on $M$ if and only if 
 \begin{itemize}[leftmargin = 10pt]
 \item[(a)] $(N, J)$ admits a K\"ahler  metric $\wt g_o$ of constant Gaussian curvature $\k = \pm 1$;  
 \item[(b)]  there is a  globally defined K\"ahler potential $\f: N \to \bR$  for $g_o$ (i.e. a real function such that $g_o = d d^c \f(\cdot, J \cdot)$)  satisfying the constraint $ \sign(\f) \equiv - \k$ and  the linear elliptic   equation  
 \beq  \label{equazione} d d^c \f(\cdot, J \cdot) +  2  \k \,\f \,\wt g_o = 0\ ; \eeq
 \item[(c)]  modulo a bundle equivalence of $\pi^\cS: M =  \cS \times I\to \cS$,     there is   a choice for the coordinate $r$ for the  factor $ I \subset \bR$  of  $ M = \cS \times I$ and a (locally defined) coordinate  set $(x, y, v, r)$  with $(x, y)$ coordinates for $N$ and $v$ for the fiber of $\pi: \cS  = N \times \bR \to N$,    such that the  $1$-form $\q$ and the Lorentzian metrics $g^{(m)}$  take the form  
 \beq\label{comp2} \begin{split}   \q &= d v +  (\pi \circ \pi^{\cS})^* (d^c \f)\\
 g^{(m)}&=   -   \frac{\k}{2}\left(r^2 + \f^2  \right)  (\pi \circ \pi^{\cS})^* \wt g_o + \q \vee \bigg( d(v + r) + \frac{1}{2}\left( - 1 + \frac{\cappa  m \, r }{r^2 + \f^2} \right)\q \bigg)\ ,
 \end{split}
  \eeq
  where $\cappa$ is just a rescaling constant  for the parameter $m$.
\end{itemize}
 \end{itemize}
 \end{theo}  
 Claim (i)  is  a  consequence  of a no-go theorem   (see Corollary \ref{nogo1}) according to which  {\it no Kerr manifold of dimension greater than $4$  admits  an Einstein metric that  is compatible with its twisting optical structure $\cQ$}. We have to admit that such a result   was  a surprising  outcome  to us, being  somehow in contrast with the above mentioned  abundance  of  higher dimensional  analogs of the Taub-NUT metrics (\footnote{However,  one can notice   that  our no-go theorem  is consistent with   Taghavi-Chabert's  classification of higher dimensional Ricci flat metrics with  mWAND structures \cite{TC} -- see also \cite{Or, OS}.}).  For what concerns claim (ii), it is a   combination of various results  that we prove for the $4$-dimensional Kerr manifolds -- see Corollary \ref{333},  Theorem \ref{aaa}, Lemma \ref{thelemma51} and  formula \eqref{6666}. \par
 \smallskip
 The  geometrical  constraints provided by  the claims (i),  (ii-a)  and (ii-b) of Theorem \ref{teoremone} are crucially rooted on the ``twisting'' property of the considered optical structures  and on the invariance under the vector field $\bT = \qs_o - B \ps_o$ with $B \neq 0$.  The existence of such constraints is in  apparent contrast with the following well-known  facts:  (1) The  classical Kerr metrics are well defined over  the whole bundle   $M^{\text{Kerr}} = (\bR^3 \setminus\{0\}) \times \bR \to S^2$ -- such a bundle  fibers over the {\it compact} K\"ahler surface $S^2$, which   therefore  cannot admit any globally defined K\"ahler potential;  (2) There do exist higher dimensional analogs of   Kerr black holes (\cite{MP, Hor}).  However, there is no actual  contradiction  --  everything is consistent.  In fact   the metrics  in  (1) and (2) are simply  not {\it globally} compatible with any twisting  optical structures. More precisely: (a)  
  The classical Kerr metrics are compatible with a prescribed  optical structure on $M^{\text{Kerr}} = S^2 \times \bR_+ \times \bR $, but such an optical  structure is twisting only if it is restricted  to  $M^{\text{Kerr}}_{(+)} = S^2_{(+)} \times \bR_+ \times \bR $ (or to  $M^{\text{Kerr}}_{(-)}  = S^2_{(-)} \times \bR_+ \times \bR$);  (b) In   \cite{PPO} (see also   \cite{MP, FS})  it is proven that   the   geodesic null congruences   of  Mayer  and Perry's  higher dimensional black holes  are not shearfree and are therefore   not compatible with  any  (neither  twisting nor  with no twist) optical structure as in the  above definition. \par
 We thus think  that  the two conditions, given by the   compatibility  with a prescribed   (twisting)  optical structure   and the  invariance  under  a vector field  $\bT = \qs_o - B \ps_o$ with $B \neq 0$,  are somehow  conflicting  in dimensions greater than $4$.  In fact,  removing   the    condition $B \neq 0$    allows  the construction of  the   higher dimensional    Taub-NUT metrics, while  removing  the compatibility  requirement    with a  prescribed optical structure  opens  the way to  the construction of   Mayer and Perry's  higher  dimensional  black holes.\par
  \medskip 
 Despite of the    just  discussed  limitations, claim (ii) of Theorem  \ref{teoremone}   shows the existence of    a very large class of  families of $4$-dimensional Ricci flat gravitational fields of Kerr type, being such families
  easily constructible for  any  choice of a  non-compact Riemann surface $(N, J)$. In fact, they are  in   bijections  with the  solutions $\f$  to \eqref{equazione}  with  the prescribed  sign property.  Using  isothermal coordinates  for  $\wt g_o$,  in which the metric  is given by  the classical formula  $\tilde{g}_o=\frac{4 (dx^2 + dy^2)}{(1+\k (x^2 + y^2))^2}$, such  conditions on $\f$  take the form
\beq \label{firsteq**}   \Delta \f  + \k \frac{8 \f}{\left(1 + \k (x^2 + y^2)\right )^2} =0\qquad \&\qquad  \k   \f < 0\ .\eeq 
Classical results  on second order elliptic equations yield that the solutions to   \eqref{firsteq**} in a disk $\bD(0, \gr_o) \subset \bR^2$ of radius $\gr_o < 1$, are real analytic and in bijection with the  Fourier series of  the non-negative   real analytic functions on the boundary of such a   disk.  Using this, in \S \ref{expli-sol} we succeed in writing explicitly   all solutions to \eqref{firsteq**}   over a  closed disk as  sums of  series. 
 Moreover, considering  the complex coordinate $z= x + i y $, we  get that the   solutions    to  \eqref{firsteq**}  over  a closed disk   $\overline{\bD(0, \gr_o)} $ coincide with the  functions  of the form  
\beq\label{19}
\begin{split}
\f(z) &= \Re \left(F(z)- \frac{2|z|^2}{1+|z|^2}\frac{1}{z} \int_0^z  F(\z)d\z\right) \ \ \text{in case}\  \k=1\ \ \text{or}\\
 \f(z) &= \Re \left( \frac{1+|z|^2}{1-|z|^2} F(z)+ z\frac{dF(z)}{dz}\right)\ \ \text{in case}\   \k=-1\ ,  \ \ 
\end{split}
\eeq
 for  some holomorphic function   $F$ 
  constrained just by the sign condition $\k \f < 0$. Notice that $F$ may be defined as a holomorphic function on any subset of $N$.  In  this regard we would like  to mention that  the contents of Theorem \ref{teoremone} (ii)  are very likely   related with  Debney, Kerr and Schild's classification of Ricci flat  algebraic degenerate metrics with non-zero complex expansion rate \cite{DKS}.  Indeed, on the one hand, as  a consequence of our Proposition \ref{flatness},  any of our metrics \eqref{comp2} admits the   presentation \eqref{216} and  is therefore   a {\it Kerr-Schild metric} (see e.g. \cite[Ch. 32]{SKMHH} for the definition and references to the vast literature on this  class of Lorentzian metrics). On the other hand,   the  Ricci flat metrics considered in \cite{DKS}  are proved to be  not only Kerr-Schild  and   compatible with  an appropriate geodesic shearfree null congruence,  but    also  including a distinguished subclass of metrics with a one-dimensional isometry group and  parameterised by  the real parts of  the  solutions of  a certain ordinary complex differential equation,  a result  which is strongly reminiscent   of ours. Investigations on the  relation between these two noticeable classes of metrics are left to future activities. \par
\medskip
 As it is natural to expect,  among the   solutions to \eqref{firsteq**} for  $\k = + 1$,  there are   potentials $\f$ which,  under an  appropriate change of coordinates,  take the form $\f(\xi, \psi) \: = a \cos \xi$,  $a \in \bR$, and are thus corresponding to  the  classical Kerr metrics \eqref{standardKerr}.  Moreover,  the metrics    which satisfy the conditions:  
 \begin{itemize}
 \item[(i)]  the corresponding   potential  $\f$ is as in  \eqref{19}  with $F$ polynomial,
 \item[(ii)] they  admit real analytic extensions to the whole  $M^{\text{Kerr}} = (\bR^3 \setminus\{0\}) \times \bR$,   
 \end{itemize}
 are precisely   the same   metrics that are  obtained  from  \eqref{standardKerr} through a space-like rotation of the flat background.  We expect that   (i)    is not a really necessary condition and that only the Kerr metrics (and their rotated versions) are  metrics  of the form \eqref{comp2}  that admit real analytic extensions to the whole $M^{\text{Kerr}} $,  all the others being forced to develop  curvature singularities  of various kinds.  It would be consistent with   the  known  uniqueness theorems for the Kerr metrics  (see e.g. \cite{Ro1}).  Investigations  in this regard are left to  future works.\par 
\medskip
For what concerns the case  $\k = -1$,    a very large number of  local functions satisfying  \eqref{firsteq**}   exist and, in contrast with the case $\k = +1$,    a huge part   of them (namely, those  corresponding to   complex  functions $F$ that are holomorphic over the whole unit disk $\bD  \= \bD(0, 1)$ and for which  the function $\f$ defined in  \eqref{19} is positive   on at least  a subset of $\bD$)  admit unique  real analytic extensions  over    the whole  2-dimensional  Lobachevsky  space $(\bD, \wt g_o)$.  The associated  Ricci flat metrics are   well defined on the Kerr manifold  $M^{\text{Lob}} \= (\bD^2 \times \bR) \times \bR$ and, at  least  on the regions in which  $\sup |\f|$  is sufficiently small, they can be considered as  deformations of a Lorentzian metric  that has  $\bR \times \SO^o_{1,2}$ as connected component of the  isometry  group (this mirrors the fact that the Kerr metrics can be considered as deformations of the Schwarzschild metric, i.e. of a metric for  which the connected component of the isometric group is  $\bR \times \SO_3$). It follows that  most of the metrics of this second class are 
not even locally isometric to the classical Kerr metrics. \par
\smallskip
We  expect that,  by appropriately choosing the  holomorphic function $F$  that  determines the potential  $\f$  by \eqref{19},  the corresponding metric has an isometry group of dimension not greater than  $1$  (i.e., with connected component given just by the flow of  the  vector field $\bT$) and that it  is not isometric  to  any  of the previously known exact solutions to the Einstein equations.  Investigations on this issue and  on diverse aspects  of the geometries  of the Ricci flat Lorentzian manifolds given in this paper are left to future work.  However,  the  coordinate-free presentations of these manifolds and of their compatible metrics  seem to us  to be quite novel   and  suggest diverse  intriguing questions. For instance, it should be noticed that the  Ricci flat metrics  on    $M^{\text{Lob}} = (\bD \times \bR) \times \bR$ described above are not static in time and not asymptotically spherically symmetric. Therefore  they are not  suited to  model stationary gravitational fields generated by localised masses. Nonetheless, they might be considered as  non-static cosmological models or, more precisely,  as  limit forms  (for    stress-energy tensors tending to $0$)   of   non-static cosmological models,    quite different from  the most commonly considered cosmological models. 
\par
\smallskip
As a concluding remark, we would like to  emphasise that  some points  of our  proofs are crucially  based  on   long  sequences of straightforward computations.   We are  aware this makes the verification difficult. In order to avoid  any possible errors,  we  scrupulously  re-checked  each   sequence of   computations with the help of   the symbolic manipulation program {\it Maple}. Note that, with the help  of this or similar  software,  everybody  can    quickly verify {\it a posteriori}  the Ricci flatness of any of the  $4$-dimensional  metrics given  in  our classification.   Shorter  proofs   than ours are surely desirable and we hope that they will come up  in the near future.
\par
\medskip
The structure of the paper is as follows. The section  \S 2 is devoted to   a brief review of  some  facts  on  optical structures,  to  giving a convenient parameterisation of the compatible  metrics of an optical structure  and to introducing  the definitions and  first properties of  Kerr manifolds,  backgrounds of Kerr type and   Kerr families of gravitational fields. In \S 3, we prove our no-go theorem for Einstein  Kerr manifolds in higher dimensions and a corollary that implies (i) of Theorem \ref{teoremone}. In \S 4 and \S 5, we classify the $4$-dimensional Kerr manifolds admitting backgrounds of Kerr type and  Kerr families of Ricci flat gravitational fields, proving in this way  (ii) of Theorem \ref{teoremone}. In  section \S 6,  explicit expressions for the potentials of the metrics of our classification are determined and the above  mentioned discussion  of  the sets of solutions corresponding to $\k = \pm 1$  is given in  greater detail. In two appendices, we provide the details of  some computations for the proofs in  \S 3 and \S 4.\\[5pt]
 \noindent{\it Acknowledgments.}
We are very grateful to Marcello Ortaggio for his invaluable assistance in providing us with several pieces of information and constructive feedback.\\[5pt]
{\it Notation.}
For any pair of $1$-forms  $\a, \b$ on a manifold $M$,   the symbol  $\a \vee \b$ stands for the symmetric tensor product $\a \vee \b  \= \frac{1}{2} \left(\a \otimes \b + \b \otimes \a\right)$.  Given a vector field $X$ and a pseudo-Riemannian metric $g$ on $M$, we denote $X^\flat \= g(X, \cdot)$. 
\par
\medskip
\section{Kerr  families of gravitational fields   on Kerr type backgrounds}
\subsection{Optical geometries}
\label{section1}
Let $(M, g)$ be  a  Lorentzian $n$-manifold.  A {\it  null congruence on $M$}  is a   foliation  of  the manifold  by    curves, which are   tangent to a nowhere vanishing null  vector  field at all points.  
The  following objects are naturally associated with  a null  congruence:
 \begin{itemize}[leftmargin = 15pt]
 \item[--] The class  $\{\ps\}$ of all  null vector fields, which are  tangent to the curves of the foliation; 
 \item[--] The codimension  one  distribution $ \cW:=\ps_o^{\perp}$ which  is   $g$-orthogonal to  some (hence,  to all) nowhere vanishing vector field  $\ps_o \in \{\ps\}$;
 \item[--] The degenerate metric $ h= g|_{\cW\times \cW}$   induced  on $\cW$.
 \end{itemize}
 The null congruence is   called {\it \sf} if  there is a vector field $\ps$ in  the class $\{\ps\}$, whose  (local)  flow  preserves  $\cW$ and the conformal class  $[h]$ of  $h$,  i.e.
\beq \label{twocond}      \cL_{\ps} X  \in   \cW \quad \text{for any} \ X \in\cW\ \quad \text{and}\quad   \cL_{\ps} h = f h\quad\text{for some function}\ \   f \in \cF(M). \eeq
This is in turn equivalent to the condition  (see e.g. \cite[Lemma 2.1]{AGSS})  
\beq \label{defsf2}   \cL_{\ps} g = f g + {\ps}^\flat \vee \h\qquad \text{for some function $f$ and a $1$-form}\ \ \h\ \eeq
and, if this is the case, then  $\ps$ is geodesic and all curves of the congruence are geodesics (\cite[Prop. 2.4]{AGSS}). Due to this, often the  \sf null congruences are alternatively   called {\it geodesic \sf null congruences}. \par
\smallskip
The property of a foliation by curves to  be a   \sf congruence  with respect to some metric  $g$ is encoded in the following definition,  which have appeared  (in possibly  different but all equivalent ways)  in \cite{RT, AGSS, AGSS1}. Actually,  the following  is a very mild generalisation of the previous, because it  deals with   Lorentzian metrics  that might  have any of the two possible  types of signatures,  the mostly  plus and the  mostly minus.  
\begin{definition}[] An {\it optical geometry} on $M$ is a  quadruple 
$\cQ = (\cW, [h]_{\pm}, \cK, \{g\})$, 
 given by 
\begin{itemize}[leftmargin = 20pt] 
\item[(i)] a codimension one distribution $\cW \subset TM$,  
\item[(ii)] the union  $[h]_{\pm} = [h] \cup [-h]$  of  two conformal classes   $[h]$, $[-h]$, one  consisting of  semipositive  degenerate   metrics  $h$ on  $\cW$ with  one-dimensional kernels,  the other  made of their opposite metrics  $- h$, 
\item[(iii)] the $1$-dimensional distribution $\cK$, given by the kernels  $\cK|_x = \ker h_x \subset \cW_x$, $x \in M$,  of  the degenerate metrics $h$, 
\item[(iv)]  the class $\{g\}$ of all  Lorentzian metrics  (with either  mostly plus  or  mostly minus signs)   inducing on $\cW$ one of the   degenerates metrics in  $ [h]_{\pm}$,  
\end{itemize}
satisfying the following conditions:  {\it the class $\{g\}$ is not empty and  there exists a vector field  $\ps$ in   $\cK$ such that \eqref{twocond} holds}. The  metrics in  $\{g\}$ are said  to be {\it compatible with $\cQ$}.
\end{definition}
\begin{rem} If  $\cQ = (\cW, [h], \cK, \{g\})$ is an optical geometry,   both  pairs   $(\cW, [h])$  and $( \cK_h = \ker h, \{g\})$  provide sufficient information to  recover the other   terms of the quadruple.  Thus, an optical geometry can be equivalently   defined  as  a  pair of   the form  $(\cW, [h])$ or   of the form  $(\cK, \{g\})$, provided that its  two terms   satisfy appropriate conditions.  The pairs  $(\cW, [h])$ of such a kind  have been named   {\it \sf structures} in  \cite{AGSS},   while  the corresponding  pairs $(\cK, \{g\})$ are   {\it Robinson and Trautman's optical geometries} considered  in   \cite{RT}.    
\end{rem} 
\par
\medskip
 \subsection{Kerr manifolds} \label{sect2}
 \begin{definition} \label{Kerrmanifoldsdef}A {\it Kerr  structure of  dimension $n = 2 k$} is a quadruple 
 $$\cM = (M = \cS\times \bR , \pi: \cS \to N, (J, g_o),   \cH)\ ,$$ 
 given by 
 \begin{itemize}[leftmargin = 22pt]
 \item[(a)]  An   $n$-dimensional  manifold $M$,  identifiable with  a Cartesian product   $M = \cS \times \bR$ between $\bR$ and   a total space $\cS$ of an $\bR$-bundle $\pi: \cS \to N$  as  in (b);  
 \item[(b)]  A principal $\bR$-bundle $ \pi: \cS \to N$ over a  K\"ahler manifold $N$; 
 \item[(c)] A K\"ahler structure $(J, g_o)$ on $N$; 
 \item[(d)] A connection $ \cH \subset T\cS$ on $ \pi: \cS \to N$, which is the kernel distribution of   an $\bR$-invariant  connection $1$-form $\theta^{\cS}$  satisfying the condition
 $d \theta^{\cS} =  \pi^* \o_o $, in which  
    $ \o_o \= g_o(J\cdot, \cdot)$  is the  K\"ahler form of $(N, J, g_o)$. 
 \end{itemize}
 A manifold $M$ as in (a), equipped with a  Kerr structure $\cM = (M = \cS\times \bR , \pi: \cS \to N,$ $(J, g_o),   \cH)$,  is   called {\it Kerr manifold}.   By a small abuse of language,  in what follows we  are going to   call ``Kerr manifolds''  also any    connected open subset   $M' = \cS \times I$,   $I \subset \bR$,   of a Kerr manifold $M = \cS \times \bR$ as defined above.
  \end{definition}
   \begin{rem} The identification $M = \cS \times \bR$, given  by  a  Kerr  structure   on a manifold $M$,  defines   a natural   $\bR$-action on such a manifold and, consequently,  a local $\bR$-action on any  connected open subset $M' = \cS \times I$, $I \subsetneq \bR$. On the other hand, any  such    $M'  \subset M$  is  clearly diffeomorphic to  $\cS \times \bR$ and hence it can be equipped with  an  $\bR$-action  of its own.   However  such $\bR$-action
 is  essentially different from  the   local $\bR$-action defined above.  So,  in what follows, we    write  $M' = \cS \times I$, $I \subset \bR$,  whenever we need to stress that the (local) $\bR$-action to be considered on $M'$ is the  one inherited from the    $\bR$-action of the larger manifold $M = \cS \times \bR$.
  \end{rem}
A Kerr manifold  $M = \cS \times \bR$  is naturally equipped with the following optical geometry. 
 Let  $\frac{\p}{\p u}$ be the standard coordinate vector field of $\bR$ and denote by  $\ps_o ( \equiv  \frac{\p}{\p u})$ the   naturally corresponding  vector  field on  $M = \cS \times \bR$.  Let also $\qs^{\cS}_o$ be the  fundamental vertical vector field of the principal bundle $\pi^{\cS}: {\cS} \to N$ corresponding to the standard basis  ``$1$'' of  $Lie(\bR) = \bR$ (i.e. such that  the flow $\Phi^{\qs_o}_s$, $s \in \bR$,  is the   family of diffeomorphisms,  determined by the right actions on the $\bR$-bundle $\cS$  of the real numbers $s  \in \bR$) and denote by $\qs_o (\equiv  \qs_o^\cS)$ the  corresponding vector field on the Cartesian product  $M =   \cS \times \bR$. Finally,  let $\cW_o$, $\cK_o$,  $h_o$ and $\wc g_o$ be the  two distributions,  the semi-positive scalar metric  on   $\cW_o$  and the Lorentzian metric that are   defined at  the points $x = (y, u) \in M = \cS \times \bR$ by 
 \beq 
\begin{split} & \cW_o|_{(y,u)} \= \cH_y + \langle \ps_o|_u\rangle\ ,\qquad \cK_o|_{(y,u)} \=  \langle \ps_o|_u\rangle\ ,\\[5pt]
&h_o (X + \lambda \ps_o, Y + \mu \ps_o) \= g_o(\pi^{\cS}_*(X), \pi^{\cS}_*(Y))\ ,\qquad X, Y \in \cH\ ,\ \lambda, \mu \in \bR\ ,\\[5pt]
& \wc g_o(Z + \lambda \qs_o + \lambda' \ps_o, Z' +  \mu \qs_o +  \mu' \ps_o) \= h_o(X,Y) +\frac{1}{2}( \lambda \mu' + \lambda' \mu)\ ,\qquad Z, Z' \in \cH\ , \\
& \hskip 12 cm \lambda, \lambda', \mu, \mu'  \in \bR\ .
\end{split}
\eeq	
Since $\wc g_o$    induces $h_o$ on $\cW_o$,   the class $\{g\}$ of the  metrics   inducing an element  of  $[h_o]_{\pm}  = \{\s h_o, \s > 0\ \text{or} \ \s < 0\}$ on $\cW_o$,  is not empty. Moreover,  by construction,  $\cL_{\ps_o} X  \in   \cW_o$ for any $X \in \cW$,  $ \cL_{\ps_o} h_o = 0$ and $\cK_o|_x = \ker h_o|_x$,  $x \in M$. We therefore have that  
 $\cQ_o = (\cW_o, [h_o]_\pm, \cK_o, \{g\})$ is an optical geometry. \par
 \begin{definition}   Let $M$ be a Kerr manifold, equipped with  a Kerr structure  $\cM = (M = \cS\times \bR , \pi: \cS \to N,$ $(J, g_o),   \cH)$. The above defined quadruple $\cQ_o $ is  called  {\it $\cM$-canonical optical geometry of $M$}.  The  optical geometries of this kind are called {\it of Kerr type}. 
\end{definition}
\noindent{\bf Remark.} The  \sf structures  corresponding to    optical geometries of Kerr type  are typical  examples of  {\it \sf structures of K\"ahler-Sasaki type}  in the sense of  \cite{AGSS}.  
\par
\medskip 
      \subsection{Parameterisations of   compatible metrics  and adapted frame fields} \label{parameterisation}
      Consider a  Kerr manifold $M$, equipped with  the  Kerr structure $\cM = (M = \cS\times \bR , \pi: \cS \to N, (J, g_o),   \cH)$ and the corresponding    $\cM$-canonical  optical geometry 
      $\cQ_o = (\cW_o, [h_o], \cK_o, \{g\})$. Let also $\ps_o$ and $\qs_o$  be the vector fields  defined above and $\wc \pi = \pi \circ \pi^{\cS}: M \to N$ the    projection  onto  $N$,  determined  by composing the  natural  projection  $\pi^{\cS}: M = \cS \times \bR \to \cS$ with   the $\bR$-bundle projection $\pi: \cS\to N$.  Each tangent space $T_x M$, $x = (y, u)$,  admits the following natural direct sum decomposition 
      \beq \label{decomp} T_{(y, u)}  M = \cH_y +  \bR \ps_o|_u + \bR  \qs_o|_y = \cW_o 
       + \bR \qs_o|_y\ .\eeq
  We remark that a  Lorentzian metric  $g$  on $M$ is  compatible with $\cQ_o$ (for simplicity, later, we  call it just  {\it compatible}) if and only if $g|_{\cW \times \cW} = \s \wc \pi^* g_o$ for a smooth real function $\s> 0$ or $\s < 0$  at all points.  This means that a  compatible metric $g$ is  uniquely determined by the following four objects:  the    function $\s \neq 0$,  the two functions
      $$x \longmapsto g_x(\ps_o, \qs_o) \ ,\qquad x \longmapsto g_x(\qs_o, \qs_o)$$ 
      and the tensor field  in $\cH^*$ that is  defined at each point $x$  by  
      $v \in \cH_x \longmapsto g_x(\qs_o, v)$.
      Note that  the function  $x \longmapsto g_x(\ps_o, \qs_o)$ is a nowhere vanishing function (otherwise  the metric would  be   no longer  non-degenerate). These observations motivate the following 
      \begin{definition} The {\it canonical datum} for a compatible metric $g$  is a quadruple $(\s, \a, \b, \gd)$ given by three smooth  real functions $\s$, $\a$, $\b$  and a tensor field $\gd \in \cH^*$,  with $\s$,  $\a$  nowhere vanishing.  Given a canonical datum, there is a uniquely associated  compatible metric $g$, namely   the  Lorentzian metric defined by 
      \beq g|_{\cW \times \cW} = \s h_o\ ,\qquad g(\qs_o, \ps_o) = \frac{\s \a}{2}\ ,\qquad g(\qs_o, \qs_o) =\frac{ \s \b}{2}\ ,\qquad g(\qs_o, \cdot)|_{\cH} = \frac{ \s \gd}{2} \ .\eeq
      \end{definition}
Identifying   $\gd\in \cH^*$  with the   $1$-form of $M$ that coincides with $\gd$ on $\cH^*$ and  vanishes identically on the complementary  spaces $\langle \ps_o|_x , \qs_o|_x \rangle$, and denoting by $\q$ and $\ps_o^*$ the $1$-forms defined by 
      \beq \q_x(\qs_{o}|_x) = 1\ ,\quad \q_x|_{\cW_x} = 0\ ,\quad \ps_{ox}^*(\ps_o|_x) = 1\ ,\quad \ps_{o x}^*|_{\cH_x +  \bR \qs_o|_x} = 0\quad \text{for any}\ x\in  M\ ,\eeq 
    the compatible Lorentzian metric   corresponding to  $(\s, \a, \b, \gd)$ is
      \beq g =\s \wc \pi^* g_o + \q \vee \bigg(\s \a \ps^*_o + \s \gd +   \frac{\s \b}{2} \q\bigg) =  \s\bigg(\wc \pi^* g_o + \q \vee \bigg(\a \ps^*_o + \gd +   \frac{\b}{2} \q\bigg)\bigg)\ .\eeq
\indent      {\it Note that any compatible metric has this form}.      
 \par
 \medskip
Now, for any (local) vector field $X$ on the K\"ahler manifold  $(N, J)$  there is a unique vector field  in  the distribution $\cH$  on $M = \cS \times \bR$, determined by the spaces $\cH_y = \cH_{(y, u)}$ of  the decompositions \eqref{decomp},    projecting onto  $X$. We denote such a vector field by $\wh X$ and  call it {\it the lifted vector field corresponding to $X$}.\par
\smallskip
The following lemma shows  how to locally express each lift  $\wh X$ in terms of the corresponding vector field $X$ on $N$ and a    potential    for  the K\"ahler metric $g_o$  (we recall that, if $\o_o = g_o(J \cdot, \cdot)$ is the K\"ahler form of $g$ and  $d^c$ is  the differential operator    on  functions defined by   $d^c f  = - d f \circ J$,  a {\it potential}  is a  (local) function $\f$ such that $\o_o = d d^c \f$).\par
\begin{lem} \label{thelemma51} There exists 
\begin{itemize}
\item[(i)] an open cover $\{\cV_\cA\}$ of  the  K\"ahler manifold $(N, J)$ of  the Kerr structure $\cM$, 
\item[(ii)] an associated    family of trivialisations $\wc \pi{}^{-1}(\cV_\cA) \simeq \cV_\cA \times \bR \times \bR$ of $\wc \pi: M \to N$, 
\item[(iii)]  a corresponding  family of  K\"ahler potentials $\{\f_\cA: \cV_A \to \bR\}$ on $(N, J, g_o)$,
\end{itemize}
 such that for any vector field $X$ on  a set  $\cV_A$, the corresponding  lifted vector field $\wh X$ on $\wc \pi{}^{-1}(\cV)  \simeq \cV \times \bR \times \bR$ is  given by   
\beq \label{theEhat} \wh X =  X -  d^c \f(X) \qs_o =  X + d \f (J  X) \qs_o  \ .  \eeq
\end{lem} 
\begin{pf} We recall that  the local sections $f: \cV \to   \wc \pi{}^{-1}(\cV)$ of the $\bR^2$-bundle $\wc \pi:   M \to N$  are in one-to-one correspondence with 
the trivialisations $\xi: \wc \pi{}^{-1}(\cV)  \to \cV \times \bR \times \bR$  of $M$.   It is the correspondence that associates   $f$ with   the unique $(\bR \times \bR)$-equivariant map 
$\xi^{(f)}: \wc \pi{}^{-1}(\cV)  \to \cV \times \bR \times \bR$   satisfying 
$( \xi^{(f)} \circ f)(y)  = (y, 0, 0)$, $y \in \cV$.
Since $M$ is the Cartesian product $M  = \cS \times \bR$, any  local section $f$ has the form 
$$f(y) = (f^\cS(y), f^\bR(y)) \ \text{for a local section}\  f^\cS: \cV \to \pi^{-1}(\cV) \subset \cS \  \text{and a  map} \  f^\bR: \cV \to \bR\ .$$
We consider only the  local sections with $f^\bR = 0$.  For any such section, the corresponding trivialisation identifies  any point $x \in M$ with a triple $(y, \r, u = 0)$ in $(\cV \times \bR) \times \bR$ and any  tangent space $T_{x} M$  with the direct sum 
$T_{x} M \simeq T_y \cV + T_\r \bR + T_0 \bR =  T_y \cV + \bR + \bR $.
Under this identification,  for  any  $x \in \wc \pi{}^{-1}(\cV)$, the vectors $\qs_o|_x$ and ${\ps_o}|_x$ are represented by the unit elements ${\bf 1}^{(I)}$ and ${\bf 1}^{(II)}$ of the  Lie algebras $\bR^{(I)}$, $\bR^{(II)}$  (both isomorphic to $\bR$) of the two copies of the Lie group $\bR$.   By  the definition of the  distribution $\cH$  of $M = \cS \times \bR$, the   lift $\wh X \in \cH$ of a vector field $X$ of $\cV$ is identifiable with   a vector field of $\cV \times \bR \times \bR$ of the form 
$ \wh X = X  +   \k^f(X) {\bf 1}^{(I)}$  where $\k^f$ is  the $1$-form on $\cV$ (determined by  the section $f = (f^{\cS}, 0)$),  which   is the potential $1$-form $\k^f \=  -  f^{\cS*} \theta$ on $\cV$ for the connection $1$-form $\theta$. In particular,   
\beq \label{condition5.4**} d \k^f = -  \o_o \ .\eeq
 For any   $y_o$ in $N$,  we may 
consider a simply connected neighbourhood   $\cV$,   in  which there are  a trivialisation for $\wc \pi{}^{-1}(\cV) \subset M$ {\it and}   a K\"ahler potential $\f$. In this neighbourhood, considering   the   local section $f(y)= (f^{\cS}(y), 0)$ that determines the trivialisation and  due to  \eqref{condition5.4**},  we have that the corresponding $\k^f$ satisfies  $d(\k^f + d^c \f) =- \o_o + \o_o =   0$. Hence  there exists   a  function $ \psi: \cV \to \bR$ such that 
$\k^f + d ^c \f =  - d \psi  $ and  the modified section  $f' = (f^{\cS}  + \psi, 0)$ is such that, in the corresponding  trivialisation, each  lifted vector field has the form $ \wh X = X +  (\k^{f}(X) + d\psi(X))\qs_o = X - d^c \f(X)\qs_o  $, as desired.
\end{pf}
For any trivialisable open subset $\wc \pi^{-1}(\cV) \simeq \cV \times \bR \times \bR$, associated with an open subset $\cV$ of the cover  $\{\cV_A\}_{A \in \cJ}$  of the Lemma \ref{thelemma51},  we
may consider a local frame field  $(E_1, \ldots, E_{n-2})$  on $\cV \subset  N$ and the uniquely associated frame field   $(X_A)_{A = 1, \ldots n}$ on $\wc \pi^{-1}(\cV) $ given by 
\begin{multline*} X_i = \wh E_i =  E_i - d^c \f(E_i) \qs_o = E_i + J_i^j E_j(\f) \qs_o \  , \  1 \leq i \leq n-2\ ,\\ X_{n-1} = \ps_o\ ,\qquad X_n = \qs_o\ ,\end{multline*}
where  $J_i^j$ denote  the components of the complex structure $J$ in the frame field $(E_i)$. 
The frame fields of this kind are called {\it adapted to the $\cM$-canonical optical structure} (or  just {\it adapted}, for short).  
The dual coframe fields of the frame field $(E_i)$ on $N$ and of  the frame field $(X_A) = (X_i = \wh E_i, X_{n-1} = \ps_o, X_n = \qs_o)$ on $M$, are denoted by $(E^i)$ and 
$(X^i = \wh E^i,  X^{n-1} = \ps^*_o, X_n = \qs^*_o )$,  respectively. \par
\medskip
\subsection{Kerr families of gravitational fields on Kerr type backgrounds}
 \label{notation}
\begin{definition} \label{ansatz} 
A {\it background of Kerr type} is a Kerr manifold $M$, equipped with a Kerr structure $\cM = (M = \cS\times \bR$ , $\pi: \cS \to N, (J, g_o),   \cH)$ and  an Einstein   metrics $\h$, which is compatible with  the $\cM$-canonical optical geometry and  associated with a canonical datum of the form   
$\big(\s, \a = \frac{1}{\s}, \b = \frac{B}{\s}, \gd = 0\big)$
 for a  constant   $B \neq 0$  and a function $\s \neq 0$  satisfying the constraint  $\bT(\s) = 0$ for  $\bT \=   \qs_o -  B \ps_o$.
In other words, $\h$ is an Einstein metric of the form 
\beq \label{back} \h =\s \, \wc \pi^* g_o + \q \vee \bigg(\ps^*_o  +  \frac{B}{2} \q\bigg)\ ,\qquad\text{with}\qquad  \wc \pi = \pi \circ \pi^{\cS}: M \to N\ , \eeq
with $B \neq 0$ constant and $\s\neq 0$ nowhere vanishing function  satisfying the above condition.
 \end{definition}
 \par
 \medskip
 The reason  for considering  Lorentzian manifolds of this kind comes from the following proposition, according to which  the $4$-dimensional backgrounds of Kerr type  are just flat Lorentzian  manifolds admitting   an optical structure of Kerr type. We recall that, according to the narration of R. P. Kerr  of the origin of his celebrated metrics \cite{Ke3} {\it the classical  Kerr metrics were discovered    studying  deformations of flat metrics compatible  with optical structures  of such a  kind}.\par
 \begin{prop} \label{flatness} Let $M$ be  a $4$-dimensional Kerr manifold, with  Kerr manifold structure $\cM$. Any  background metric   $\h =\s \wc \pi^* g_o + \q \vee (\ps^*_o  +  \frac{B}{2} \q)$  of  Kerr type on $M$  has identically vanishing  curvature. 
\end{prop}
\begin{pf}  Let $(N, J, g_o)$ be the $2$-dimensional K\"ahler manifold, which occurs in the definition of  the Kerr manifold structure  of $M$ and 
$\cV \subset N$ one of the  open sets of the  cover $\{\cV_A\}_{A \in \cJ}$, with  associated trivialisations $\wc \pi{}^{-1}(\cV_A)  \simeq \cV_A \times \bR \times \bR$  and  K\"ahler potentials $\f_A: \cV \to \bR$,  as   in   Lemma \ref{thelemma51}.  If $(v, u)$  are   coordinates on   the standard fiber  $\bR  \times \bR $ of $\wc \pi{}^{-1}(\cV) $,  the $1$-forms   $\ps^*_o$  and $\q$ have the form
$$ \ps_o^* = d u\ ,\qquad  \q  = d v + d^c \f\ , $$
Consider   the modified fiber coordinates 
\begin{equation} R =  \frac{1}{\sqrt{2|B|}} u+ \operatorname{sign}(B) \sqrt{\frac{|B|}{2}} v \ ,\ \  T = - \frac{1} {\sqrt{2 |B|}} u\ .
  \end{equation}
In these coordinates, we have that   $\bT =  \operatorname{sign}(B) \sqrt{\frac{|B|}{2}} \frac{\p}{\p T}$ (so that, since $\bT(\s) = 0$,  the function  $\s$ is  independent of $T$)  and 
\beq
\begin{split}
& \ps_o^* = du  = - \sqrt{2 |B|} d T \ ,\qquad \q =   d v + d^c \f  =  \operatorname{sign}(B) \sqrt{\frac{2}{|B|}} d R +\operatorname{sign}(B) \sqrt{\frac{2}{|B|}} d T + d^c \f \ ,\\ 
&  \q \vee (\ps_o^* + \frac{B}{2} \q) = \\
& \hskip 1 cm = \frac{B}{2} \left(\sqrt{\frac{2}{|B|}}d R +  \operatorname{sign}(B) d^c \f\right)\vee  \left( \sqrt{\frac{2}{|B|}}d R +  \operatorname{sign}(B) d^c \f\right) -    \operatorname{sign}(B) dT \vee d T \ .
\end{split}  \eeq
We therefore see that  
\beq \label{back1}
\begin{split}
& \h  = \s\big( (\pi^{\cS} \circ \pi)^* g_o \big)   + \q \vee \left( \ps_o^* + \frac{B}{2} \q\right) =  \\
&\ \  =  \operatorname{sign}(B)\bigg\{\bigg( \operatorname{sign}(B)\s\big( (\pi^{\cS} \circ \pi)^* g_o \big)   + \\
&\hskip 1 cm  +  \frac{|B|}{2} (\sqrt{\frac{2}{|B|}} d R + \operatorname{sign}(B)d^c \f) \vee (\sqrt{\frac{2}{|B|}}d  R + \operatorname{sign}(B)d^c \f)\bigg)  -   d T\vee d T \bigg\}\ .
\end{split}
\eeq
This expression for $\h$ and the fact that $\s$ is $T$-independent show that the $4$-dimensional Lorentzian manifold 
$(\wc \pi{}^{-1}(\cV) , \operatorname{sign} (B) \h)$ is locally isometric to the  Cartesian product of  the  $1$-dimensional Minkowski space   $\bR^{0,1} = (\bR, - dT \vee dT)$ and  a  totally geodesic $3$-dimensional submanifold (which is  either Riemannian or Lorentzian,  depending on whether   $B\s >  0$ or $B\s < 0$). Let us denote such a totally geodesic submanifold by  $(\wc \cS,  g^{\wc \cS} )$. 
 Being   
 \beq \label{trivia} (\wc \pi{}^{-1}(\cV) , \operatorname{sign} (B)  \h)= (\wc \cS \times \bR, g^{\wc \cS} + (- d T \vee d T))\eeq
   a  Cartesian product  of two geodesic submanifolds,     one of which is $1$-dimensional and thus  flat,
    the Ricci tensor  of this $4$-dimensional metric 
     satisfies an Einstein equation  $\Ric^{( \operatorname{sign} (B)  \h)} = \L  \operatorname{sign} (B)  \h$ for some constant $\L$   if and only if  
 \begin{itemize}[leftmargin = 20pt]
 \item[(a)] the Einstein constant $\L$  is $0$ (because  $\Ric^{( \operatorname{sign} (B) \h)}\bigg(\frac{\p}{\p T}, \frac{\p}{\p T} \bigg) = 0$  since $(\bR, - d T \vee d T)$  is  totally geodesic and flat) and   
 \item[(b)]    the Ricci tensor  $\Ric^{(g^{\wc \cS})}$  of 
 $\big(\wc \cS,  g^{\wc \cS} \big)$ is  identically vanishing  (this is  because $\wc \cS$ is a totally geodesic submanifold and the Einstein constant is $\L = 0$, by (a)). 
 \end{itemize}
 By a classical fact on  $3$-dimensional pseudo-Riemannian manifolds (see  e.g. \cite[Prop. 1.120]{Be}), the Riemann curvature tensor $R^{(g^{\wc \cS})}$ of $\big(\wc \cS,  g^{\wc \cS} \big)$ is uniquely determined by its Ricci tensor   $\Ric^{(g^{\wc \cS})}$ and  the latter  vanishes  if and only if $R^{(g^{\wc \cS})} = 0$. This means that $\big(\wc \cS,  g^{\wc \cS} \big)$ is flat and  the curvature of \eqref{trivia} vanishes identically. This proves the claim.
\end{pf}
\par 
\medskip
We are now ready to introduce the class of Lorentzian metrics, on  which  we focus  and that  can be considered as  natural generalisations of the classical Kerr metrics. \par
 \begin{definition} \label{Kerrfamiliesdef}
 Let $(M, \cM, \h)$  be a background of Kerr type, hence with $\h$ of the form \eqref{back}.  A {\it  Kerr  family of gravitational fields}   on  $(M, \cM, \h)$  (or  in  any  of its connected open subsets $M' = \cS \times I$, $I \subset \bR$)
 is a  one-parameter   family 
$g^{(m)}$  of compatible Einstein Lorentzian metrics  on  $M$, which is  determined by canonical data 
$\big( \s^{(m)} \= \s,  \a^{(m)} \= \frac{1}{\s}, \b^{(m)} \= \frac{B + m \wt \b_o}{\s},  \gd^{(m)} \= 0\big)$, i.e. of metrics of the form  
\beq \label{Kerr-Schild} g^{(m)} =\s \wc \pi^* g_o + \q \vee \bigg(\ \ps^*_o  +  \frac{B +  m \wt \b_o}{2} \q\bigg) = \h +  \frac{m \wt \b_o}{2}  \q \vee  \q \ ,\eeq 
where  $\wt \b_o$ is an   $m$-independent function   satisfying the condition  $\bT(\wt \b_o) = 0$.   
\end{definition}
From   \eqref{Kerr-Schild},  we see that a Kerr family   $g^{(m)}$ is  a one-parameter deformation, linearly depending on  $m$, of the background  metric $\h$: 
\beq \label{216} g^{(m)} = \h + \d g^{(m)}\ ,\qquad \text{with}\qquad  \d g^{(m)} \= m \frac{\wt \b_o}{2}  \q \vee  \q\ . \eeq
Motivated by a well-known terminology  for    the classical Kerr metrics,  we call 
  the decomposition  $g^{(m)} =  \h +  \frac{m \wt \b_o}{2}  \q \vee  \q $     the   {\it Kerr-Schild presentation} of the   family $g^{(m)}$. \par
\medskip
\begin{rem} \label{importantrem}  According to    \eqref{back1},   a  $4$-dimensional background metric $\h$ locally decomposes into $\h = -  (\sign(B) d T \vee d T) + \sign(B) g^{\wc \cS}$, where $  \sign(B) g^{\wc \cS}$ is the induced metric of  any totally geodesic   hypersurface $\wc \cS = \{ T = \text{const.}\}$.  This decomposition shows that the vector field $\bT = \sign(B) \sqrt{\frac{|B|}{2}} \frac{\p}{\p T}$ has a constant  squared $\h$-norm  $-\frac{B}{2}$ and that  its nature changes from  time-like to space-like according to  whether $\sign(B) g^{\wc \cS}$ is (positive/negative) defined  or  is  Lorentzian. In the cases in which $B{\cdot} \s > 0$, the metric $\sign (B) g^{\wc \cS}$ is either positive or negative  defined and in both cases $\bT$ is time-like. In the other  cases, $g^{\wc \cS}$ has a Lorentzian type signature and $\bT$ is space-like. Thus, in   coordinates $(x^i)_{i = 0, \ldots, 3}$  in which  $\bT = \frac{\p}{\p x^0}$,   the  coordinate $x^0$ can be considered as   ``time'' whenever $B {\cdot} \s > 0$  and as a  ``spacelike coordinate'' otherwise. Accordingly, the constraints  $\bT(\s) = 0$, $\bT(\b_o) = 0$  for a Kerr family of metrics correspond to imposing  that   the metrics $g^{(m)}$ are   {\it  time invariant} or   {\it   invariant under (certain) space-like  translations}.  For  what concerns  the classical Kerr  metrics,  \eqref{113} shows that they  correspond to a background $\h = g^{(m = 0)}$  in which  $B {\cdot} \sigma = (-1) \sigma > 0$  and  $\sigma < 0$, meaning that    the classical Kerr metrics  are ``time invariant'', as  is well known.
 \end{rem}
\par
\medskip
\section{A no-go theorem for  the  Einstein metrics  that are compatible  with optical geometries of Kerr type}
\subsection{A dimensional constraint for  Kerr manifolds with Kerr type backgrounds}
The purpose of this section is to  prove (i) of Theorem \ref{teoremone}. As we mentioned in the Introduction, that claim is  actually a corollary  of the following more general result.\par
\begin{theo}[Dimensional constraint] \label{nogo}Let  $M$ be a Kerr manifold, equipped with the Kerr structure $\cM = (M = \cS\times \bR$ , $\pi: \cS \to N, (J, g_o),   \cH)$ and  the corresponding $\cM$-canonical optical structure $\cQ_o$. If  $n = \dim M$ is greater than $4$, then  there is no $\cQ_o$-compatible  Lorentzian   metric $g$ satisfying  the following two conditions:
\begin{itemize}[leftmargin = 18pt]
\item[(a)] It has the form 
\beq \label{back*} g =\s \, \wc \pi^* g_o + \q \vee \bigg(\ps^*_o  +  \frac{\wt \b}{2} \q\bigg)\ ,\qquad\text{with}\qquad  \wc \pi = \pi \circ \pi^{\cS}: M \to N\eeq
for some   $\s, \wt \b: M \to \bR$,    where $\s$ is nowhere vanishing   and such that  $\bT(\s) = \bT(\wt \b) = 0$ for  the vector field   $\bT \=   \qs_o -  B \ps_o$ with  a constant $B \neq 0$; 
\item[(b)] Its Ricci curvature   satisfies 
\beq\label{semi-Einstein} \Ric(\ps_o, X) = 0\qquad \text{for any vector field} \ X \in \cW\ .\eeq
\end{itemize}
In particular, if $\dim M > 4$, then  there is no Einstein  metric with the properties    in  (a). 
\end{theo}
The proof   is given in  the next subsection.   As a   corollary,  here is the    mentioned  result:
\begin{cor}[No-go Theorem] \label{nogo1}  No Kerr manifold of dimension  greater than $4$ admits a   background metric of Kerr type. \end{cor}
\begin{pf} A background metric $\h$  of Kerr type  is Einstein and it satisfies the  property  (a) of  Theorem \eqref{nogo} with  $\wt \b = B$. Thus  no such metric exists if   $\dim M > 4$.
\end{pf}
\par
\medskip
\subsection{Proof of   the dimensional constraint}
Consider the restriction of a metric \eqref{back*} to  an open  subset  $\cU \subset M$ of the form  $\cU = \wc \pi^{-1}(\cV) \simeq \cV \times \bR \times \bR$,  where $\cV$ is an open subset of  $N$ of the  open cover $\{\cV_A\}_{A \in \cJ}$  described in  Lemma \ref{thelemma51}. 
On  such  open set,  let us 
consider an adapted frame field 
$(X_A)_{A = 1, \ldots n} = (\wh E_1, \ldots , \wh E_{n-2}, \ps_o, \qs_o)$, with  $\wh E_i$   lifts of  vector fields   $E_i$  on $\cV \subset N$. The corresponding  dual coframe field  is denoted by $(X^A) = (X^i = \wh E^i,  X^{n-1} = \ps^*_o, X_n = \qs^*_o )$. We organise the  proof Theorem \ref{nogo} for $g|_{\cU}$ into two steps.  We   first
 translate   the conditions  $\Ric(\ps_o, \ps_o) = 0$ and $\Ric(\ps_o, \wh E_i) = 0$  (which are equivalent to \eqref{semi-Einstein})  into equations  on   the potential $\f$ for $g_o$,  on the constant $B$ and on  the  function  $\s$    in     \eqref{back*}.  Then we   prove  that  such  equations  have no   solution if $\dim M > 4$.\par
 \smallskip
 \subsubsection{Notational issues}\label{section321}  Before going into the  details of the proof, we need to fix some convenient notation.
\begin{itemize}[leftmargin = 15pt]
\item We denote  by  $c_{ij}^k: \cV \to \bR$  the  functions  given  by the relations $[E_i, E_j] = c_{ij}^k E_k$.  
\item  The components   of  $g_o$,  $\o_o = g_o(J\cdot, \cdot)$ and  $J$  in the frame field $(E_i)$ on $\cV \subset N$,   are denoted by 
$$
g_{ij} \= g(E_i, E_j)\ ,\quad  \o_{ij} \=  \o_o(E_i, E_j) = g_o(J E_i, E_j)\ ,\quad  J_i^j = E^j(J(E_i)) = g^{jk} \o_{ik}\ .
$$
\item  Whenever  a tensor  component refers to   the vector fields   $X_{n-1} = \ps_o$ or $X_n = \qs_o$ or to the corresponding $1$-forms $X^{n-1} = \ps_o^*$ and $X^n = \qs_o^*$,  we mark it with  indices  ``$\ps_o$'' or ``$\qs_o$'' instead of  $n-1$ or $n$. This is  to make formulas  more readable (for instance,    the component $\Ric(X_{n-1}, X_n)$ of the Ricci tensor   is denoted by     $ \Ric_{\ps_o \qs_o} $ instead of  $\Ric_{n-1\,n}$).
\item  We denote by  $\Ga i j m$, $1 \leq i, j, m \leq n-2$,    the Christoffel symbols of the Levi-Civita connection $\n^{g_o}$  of  $(\cV, g_o)$  with respect to  the frame field $(E_i)$ 
and we denote by    $\GGa A B C$, $ 1 \leq A, B, C \leq n$,     the Christoffel symbols of the Levi-Civita connection $\n^{g}$  of  $\big(\cU =  \wc \pi^{-1}(\cV), g\big)$   with respect to  the adapted  frame field $(X_A)$. 
\item  $\RR A B C D$  are the curvature  components of  $\big(\cU =  \wc \pi^{-1}(\cV), g\big)$  in the   frame field $(X_A) $  (\footnote{Following   \cite{KN},  we  assume that   $R$  is defined by   $R_{X Y} Z :=  \n_{X} \n_{Y} Z -  \n_{Y} \n_{X} Z - \n_{[X, Y]} Z $.}). 
 \item $v$, $u$   are coordinates for  the  fiber  $\bR \times \bR$ of $\cU = \wc \pi^{-1}(\cV)$,  in which  $\ps_o$ and $\qs_o$  have the form 
 \beq \ps_o = \frac{\p}{\p u}\ ,\qquad \qs_o = \frac{\p}{\p v}\ .\eeq
 \item  Given a real function    $f: \cV \times \bR \times \bR \to \bR$, we denote   $f_\ell \= E_\ell(f)$, $1 \leq \ell \leq n -2$. In this way, by \eqref{theEhat},    the vector fields $\wh E_i$  take the form  
\beq  \label{theEhat*} \wh E_i = E_i - d^c \f(E_i) \qs_o = E_i + J_i^j E_j(\f) \qs_o = E_i + J_i^j \f_j \qs_o\ .\eeq 
 \end{itemize}
 \par
 We  now recall  that  if   $X_{A_o}$,  $X_{B_o}$ are two  commuting vector fields of the adapted frame field $(X_A)$,     any corresponding   curvature component $\RR {A_o} {B_o} C D$   is equal to  
  \beq \label{Riemannform}  \RR {A_o} {B_o} C D = X_{A_o} (\GGa {B_o} C D) -  X_{B_o} (\GGa {A_o} C D) - \GGa {A_o} C  F \GGa {B_o} F D +    \GGa {B_o} C F \GGa {A_o} F  D . \eeq
  Since
$$ [\wh E_i, \wh E_j] =   \wh{[E_i, E_j]}  - \o_o(E_i, E_j) \qs_o = c_{ij}^k \wh E_k - \o_{ij} \qs_o\ ,\qquad [\wh E_i, \ps_o] = [\wh E_i, \qs_o] =  [\ps_o, \qs_o] = 0\ ,  $$
the formula \eqref{Riemannform} holds  unless $1\le A_o=i,B_o=j \le n-2$. For  such pairs,  instead of \eqref {Riemannform},  one has to use 
  \beq \label{Riemannform1} \RR i j C D =  \wh E_i (\GGa j C D) -  \wh E_j (\GGa i C D) -  \GGa i C  F \GGa  j  F D +   \GGa j C F \GGa  i F D- c^k_{ij} \GGa k C D +  \o_{ij} \GGa {\qs_o} C D\ . \eeq
\par
\smallskip
\subsubsection{The Christoffel symbols of  the metric  \eqref{back*}}
  In \cite{AGSS1} (see also  \cite[Appendix]{AGSS}),  the complete  list   of the   Christoffel symbols of   the Levi-Civita connection of  a compatible metric  of  a manifold of K\"ahler-Sasaki type   is   given.  Since any   metric  \eqref{back*}  is of such kind, we may use the results  of that paper  to determine the Christoffel symbols $\GGa A B C$. Actually, since  \eqref{back*}  is a compatible metric of  a particularly special  form (in fact, its  associated canonical datum  has the form  $\big(  \s,  \a= \frac{1}{\s}, \b= \frac{ \wt \b}{\s},  \gd = 0\big)$ with  $\s$ and $\wt \b$ such that $\bT(\s) = \bT(\wt \b) = 0$)   the list  in  \cite[Prop. 3.1]{ AGSS1}  radically  simplifies and  we get:
     \begin{align}
\nonumber & \GGa i j m =
  g^{mk} g_o(\n^o_{E_i}   E_j, E_k)
+ \frac{1}{2 \s }  \wh E_i(\s) \d_j^m  + \frac{1}{2 \s }  \wh E_j(\s) \d_i^m
 -   \frac{1}{2 \s} g_{ij}     g^{mk} \wh E_k(\s) = \\
\nonumber &  = \Ga i j m
+ \frac{1}{2 \s } \s_i \d_j^m  + \frac{1}{2 \s }  \s_j \d_i^m
 -   \frac{1}{2 \s} g_{ij}     g^{mk} \s_k +\\
   \label{519} & \hskip 4 cm + B \frac{\ps_o(\s) \f_\ell }{ 2 \s } \left(   J^\ell_i \d_j^m +  J_j^\ell \d_i^m
 -  g_{ij}     g^{mk}  J^\ell_k \right) \ ,\\
 \label{6.12} & \GGa i j {\ps_o}  =    -  g_{ij}  \bigg( B  \ps_o(\s) - \wt \b   \ps_o(\s) \bigg)  \ ,\quad
 \GGa i j {\qs_o}  =  -\frac{   \o_{ij}}{2}  -  g_{ij}   \ps_o(\s)\ ,\\[12 pt]
\label{cli} & \GGa i {\ps_o} m = \GGa {\ps_o} i m =   \underset{=  \frac{J^m_i}{4 \s}}{ \frac{g^{mk}  \o_{ik}}{4 \s} } +   \frac{\d^m_i}{2 \s}\ps_o(\s) \ ,\qquad  \GGa i {\ps_o} {\ps_o} =   \GGa  {\ps_o} i {\ps_o} =
 \GGa i {\ps_o} {\qs_o} =  \GGa {\ps_o} i  {\qs_o} = 0 \ ,\\[12 pt]
\nonumber  &  \GGa i {\qs_o} m =  \GGa {\qs_o} i m =     \underset{=  \frac{J^m_i}{4 \s}(\wt \b -B)}{   \frac{g^{mk}  \o_{ik}  (\wt \b -B) }{ 4 \s}}  +\bigg(  \underset{=  \frac{ J^m_i}{4 \s}}{ B  \frac{g^{mk}  \o_{ik}  }{ 4 \s}} 
+ B\frac{ \d_{i}^m}{2 \s} \ps_o(\s)\bigg) =  \frac{J^m_i}{4 \s}(\wt \b -B) + B  \GGa i {\ps_o} m\ ,\\
 \label{via} & \hskip 3 cm  \
 \GGa i {\qs_o} {\ps_o} = \GGa  {\qs_o} i {\ps_o} = \frac{1}{2} \wt \b_i - \frac{1}{2}J^\ell_i \f_\ell \ps_o(\wt \b)   \ ,\ \  \GGa i {\qs_o} {\qs_o} =  \GGa  {\qs_o} i {\qs_o} = 0\ ,\end{align}
\begin{align}
 &  \GGa {\ps_o}  {\ps_o}  m  =   \GGa {\ps_o}   {\ps_o} {\ps_o}  =   \GGa {\ps_o}   {\ps_o} {\qs_o}  =  0 \ ,\qquad    \\[12 pt]
   &  \GGa {\ps_o}  {\qs_o}  m  = \GGa {\qs_o}  {\ps_o}  m  =  0 \ ,\qquad
 \GGa {\ps_o}   {\qs_o} {\ps_o}  =  \GGa {\qs_o}   {\ps_o} {\ps_o}  =  \frac{1}{2 }\ps_o(\wt \b) \ ,\qquad \GGa {\ps_o}   {\qs_o} {\qs_o} =  \GGa {\qs_o}   {\ps_o} {\qs_o}   = 0 \ ,\\
\nonumber & \GGa {\qs_o}  {\qs_o}  m =   -  \frac{g^{mk}}{4\s}\wh E_k(\wt \b) = -  \frac{g^{mk}}{4\s} \wt \b_k  - B \frac{g^{mk}}{4\s}J_k^\ell  \f_\ell \ps_o(\wt \b) \ ,\\
  \label{6.17}  &  \hskip 6 cm  \GGa {\qs_o}   {\qs_o} {\ps_o}  =\frac{\qs_o(\wt \b)}{2}+ \frac{\ps_o(\wt \b) \wt \b }{2} \ ,\ 
 \GGa {\qs_o}   {\qs_o} {\qs_o}  =    - \frac{1}{2 }\ps_o(\wt \b) \ .
  \end{align}
  Combing these formulas with   \eqref{Riemannform} and \eqref{Riemannform1},  we  are ready to  translate the conditions $\Ric(\ps_o, \ps_o) =  \Ric(\ps_o, \wh E_i) = 0$ into  equations for $\s$, $B $ and $\f$. 
\par
\smallskip
  \subsubsection{The condition $\Ric(\ps_o, \ps_o) = 0$ as an equation  on $\s$} \label{sect423}
    We recall that  
   $\Ric(\ps_o, \ps_o) =  \Ric_{{\ps_o}\ps_o}   {=} \ \RR m {\ps_o} {\ps_o} m +  \RR {\qs_o} {\ps_o} {\ps_o} {\qs_o}$.
Thus, writing  the curvature  components  in terms of the Christoffel symbols  \eqref{519} -- \eqref{6.17} and neglecting all   trivially vanishing terms, we  get 
   \begin{align}
 \nonumber \Ric(\ps_o, \ps_o)  &{=}-  \ps_o (\GGa m {\ps_o} m)-  \GGa m {\ps_o}  \ell \GGa  {\ps_o} \ell m  = \\
\nonumber  &=  - \frac{n-2}{2}\ps_o \left (  \frac{ \ps_o(\s)}{ \s}  \right ) - \left (   \frac{J_{m}^{\ell} }{4 \s}+ \frac{ \ps_o(\s)\d_{m}^{\ell}}{2 \s}  \right)\left (    \frac{J_{\ell}^m }{4 \s}+ \frac{ \ps_o(\s)\d_{\ell}^m}{2 \s}\right ) = \\
\label{6.23ter} & =\frac{n-2}{4\s^2}\left (-2\s\ps_o(\ps_o(\s)) +\ps_o(\s)^2+ \frac{1 }{4}\right ) \ .
  \end{align}
 Since $n \geq 4$,   the condition $\Ric(\ps_o, \ps_o) = 0$ is equivalent to the equation 
\beq \label{main}  - 2 \s \ps_o(\ps_o(\s)) +  \left(\ps_o(\s)\right)^2+ \frac{1}{4}=0\ .\eeq
 Differentiating   along $\ps_o$  yields that  any solution of \eqref{main} satisfies also the equation   
 $$\ps_o\left(\ps_o \left(\ps_o(\s)\right)\right) = \frac{\p^3\s }{\p u^3} = 0\ .$$  Together  with   $0 = \bT(\s)= - B\frac{\p \s}{\p u} + \frac{\p \s}{\p v}$, we get that  $\s$ has necessarily  the form
\begin{equation} \label{2} \s = C_0 + C_1 (u + B v) + C_2(u + B v)^2\ \end{equation}
for  some smooth functions $C_\a: \cV \times \bR \times \bR \to \bR$, $\a = 0,1,2$,  which are independent of the fiber coordinates  $u$ and $v$ (hence, identifiable with functions $C_\a: \cV \subset N \to  \bR$). 
Plugging \eqref{2} into \eqref{main},  we  see   that the  $C_\a$  satisfy  the algebraic equation
\begin{multline}  \label{thesigma} \left( C_1 + 2 C_2(u + B v) \right)^2 - 4 C_2  \left( C_0 + C_1 (u + B v) + C_2(u + B v)^2\right) + \frac{1}{4} = \\
= C_1^2  - 4 C_0 C_2 + \frac{1}{4} = 0\ .
\end{multline}
 In particular, $C_2 , C_0: \cV \to \bR$ are nowhere vanishing.  If for each given solution  $\s$  to   \eqref{main}, determined by the functions $C_\a: \cV \to \bR$, we  introduce the   associated  function  $r:   \cV \times \bR \times \bR  \to \bR$   defined by     $ r(y, v, u) \=u + B v + D(y)$ with $D(y) = \frac{C_1(y)}{2 C_2(y)}$, we conclude that
  $\s$ takes   the form 
 \beq  \label{thesigma1}  \begin{split}
 \s &= C_0 + C_1 r -   \frac{C_1^2}{2 C_2}  + C_2 r^2 -\frac{C_1 C_2}{ C_2} r  + \frac{C_1^2 C_2}{4 C_2^2} 
   \overset{\eqref{thesigma}}=  C_2 r^2 + \frac{1}{16 C_2}\ .  
 \end{split}
 \eeq
\par
\smallskip
\subsubsection{The conditions  $\Ric(\wh E_i, \ps_o) = 0$ as  equations on $\f$, $B$ and $\s$}
Utilising the Ricci tensor definition and expressing curvature components with Christoffel symbols  (and  removing any   term that is manifestly zero according to   \eqref{519}--\eqref{6.17}) we get
\begin{multline} \Ric(\wh E_i, \ps_o) = \Ric_{i\ps_o}   {=} \ \RR m i {\ps_o} m + \RR {\ps_o} i {\ps_o} {\ps_o} + \RR {\qs_o} i {\ps_o} {\qs_o}  = \\
  = \wh E_m (\GGa i {\ps_o} m) -  \wh E_i (\GGa m {\ps_o} m) -  \GGa m {\ps_o}  \ell \GGa i  \ell  m  +  \GGa i {\ps_o} \ell  \GGa m  \ell  m  - c^r_{mi} \GGa r {\ps_o} m\ .   
\label{6.21}
\end{multline}
We  may now  plug in    \eqref{519} -- \eqref{6.17} into \eqref{6.21}   and  simplify the  expression using  the fact  that $\n^{g_o} J = 0$ (because $(N, J,  g_o)$ is K\"ahler) and that $\n$  is torsion free.  After   a few straightforward (but quite  tedious)  computations,  we get   
  \begin{multline} \label{5.32ter} 
  \Ric(\wh E_i, \ps_o) =  \frac{n-3}{8 \s^2} \left(  \frac{n-6}{n-3} J^\ell_{i} \s_\ell       -  4 \s E_i \left( \ps_o(\s) \right)    +  4 \ps_o(\s) \s_i     -\right.\\
  \left. - \frac{n-6}{n-3}  B \f_i \ps_o \left(\s\right)      -   4 \s   B  J_i^\ell \f_\ell \ps_o \left(\ps_o(\s) \right) +    4  B  J_i^\ell \f_\ell \left( \ps_o \left(\s\right) \right)^2  \right). \ 
\end{multline}
The reader is referred to   Appendix  \ref{appendix-1}  for  the  details  of  those computations. 
\smallskip
Having $n \geq 4$,  it follows  that  each equation  $\Ric(\wh E_i, \ps_o) = 0$  is  equivalent to  
\begin{multline} \label{5.32ter*} 
 \frac{n-6}{n-3} J^\ell_{i} \s_\ell       -  4 \s E_i \left( \ps_o(\s) \right)    +  4 \ps_o(\s) \s_i     -\\
  -  \frac{n-6}{n-3}   B  \f_i \ps_o \left(\s\right)      -   4 \s   B  J_i^\ell \f_\ell \ps_o \left(\ps_o(\s) \right) +   4   B   J_i^\ell \f_\ell \left( \ps_o \left(\s\right) \right)^2   
 = 0\ .
\end{multline}
\par
\smallskip
 \subsubsection{The   proof that   $\Ric(\ps_o, \ps_o) = 0$ and $\Ric(\ps_o, \wh E_i) = 0$   are not compatible if   $n > 4$} \label{sect425}
 Assume that   $g$  is such that  $\Ric(\ps_o, \ps_o) = 0$, i.e. 
 that $\s$ has the form \eqref{thesigma1}.   Let us also   adopt the  the short-hand notation 
$A \= C_2$,  $D 	=\frac{C_1}{2 C_2}$,  $A_{i} \= E_i(A)$,   $D_i \= E_i(D)$ and $\s_i = E_i(\s)$. Since    $ r = u + B v  + D$, 
 the following relations hold: 
\beq \label{thesigma2}
\begin{split}
&\s_\ell = E_\ell\bigg(A r^2 + \frac{1}{16 A}\bigg) = A_\ell r^2 + 2 A  D_\ell r - \frac{1}{16 A^2} A_\ell \ ,\\
&\ps_o(\s) = \frac{\p \s}{\p u} = 2 A r  \frac{\p r}{\p u}  = 2 A r\ ,\quad\ps_o(\ps_o(\s)) = 2 A\ ,\quad
E_i(\ps_o(\s))  = 2 A_i r + 2 A D_i\ .
\end{split}
\eeq
If we now assume  that $g$ is also a solution to   $\Ric(\wh E_i, \ps_o) = 0$, $1 \leq i \leq n-2$, plugging   \eqref{thesigma1} and \eqref{thesigma2} into  \eqref{5.32ter*}  we get  
$$0 = \frac{n-6}{n-3} J^\ell_{i} \left(A_\ell r^2 + 2 A r D_\ell - \frac{1}{16 A^2} A_\ell\right)    - $$
\begin{multline*}
    -  4\left( A r^2 + \frac{1}{16 A }\right)  \left(2 A_i r + 2 A D_i\right)   
    + 8 A r \left( A_i r^2 + 2 A r D_i - \frac{1}{16 A^2} A_i\right)     - \\
  - \frac{n-6}{n-3}  B \f_i  2 A r     
+  8 \left( A r^2  - \frac{1}{16  A}\right)   B  J_i^\ell \f_\ell  A
+  16   B   J_i^\ell \f_\ell A^2 r^2 = 
\end{multline*}
\begin{multline} \label{5.32terter}
=  \frac{n-6}{n-3} J^\ell_{i}  A_\ell r^2 + 2\frac{n-6}{n-3} J^\ell_{i}   A  D_\ell r -   \frac{n-6}{n-3} J^\ell_{i}  \frac{1}{16 A^2} A_\ell  - \\
   - \xcancel{ 8 A  A_i r^3 } -  8 A^2 D_i r^2    -      \frac{A_i}{2  A} r -    \frac{1}{2}D_i   
   + \xcancel{ 8 A A_i r^3 } + 16 A^2 r^2 D_i -    \frac{A_i}{2 A} r     - \\
  - 2  \frac{n-6}{n-3}   B  \f_i   A r     
 + 8 A^2   B  J_i^\ell \f_\ell  r^2  -   \frac{1}{2}   B  J_i^\ell \f_\ell 
+  16   B   J_i^\ell \f_\ell A^2 r^2\ .
\end{multline}
The right hand side of the last equality  is a second order polynomial  in the function $r = r(y, v, u) $ (here,  $y = (y^i)$ stands for  coordinates  on $\cV \subset N$). By definition of  such a function, for any fixed value of $y$,   the values of  $r$  run unconstrained in an unbounded interval of $\bR$. This implies that   \eqref{5.32terter}  is satisfied if and only if  the coefficients of the  monomials   $r^2$, $r$ and $1$  vanish identically, i.e. if and only if the following equations are satisfied for any $1 \leq i \leq n-2$:   
\beq \label{4.26}
\begin{split}
 & \frac{n-6}{n-3} J^\ell_{i}  A_\ell    + 8 A^2 \left( D_i   + 
   B   J_i^\ell \f_\ell \right)  = 0\ ,\\
&
     2 A \left( \frac{n-6}{n-3} J^\ell_{i}  D_\ell -  \frac{n-6}{n-3} B \f_i  -   \frac{A_i}{ 2A^2}    
  \right)       = 0\ ,\\
&
  - \frac{1}{16 A^2}  \frac{n-6}{n-3} J^\ell_{i} A_\ell  -   \frac{1}{2}  \left(D_i      
 + B  J_i^\ell \f_\ell \right) = 0\ .
\end{split}
\eeq
Here, the third line is equivalent to the first and we may  neglect it.  Moreover,  by multiplying the first line by $J^i_m$ and contracting for $i$,   the   system   \eqref{4.26}   becomes  equivalent to   
\beq \label{6.32} \left( \begin{array}{cc} - 1 & 8  \frac{n-3}{n-6}\\[5pt]
- 1 &   2\frac{n-6}{n-3} \end{array}\right) \left( \begin{array}{c}  \frac{A_i}{ A^2}\\[5 pt]  J^\ell_{i}   D_\ell 
\end{array} \right) =   \left(\begin{array}{c}  8    \frac{n-3}{n-6}  \\[5 pt]  2 \frac{n-6}{n-3}  \end{array}\right)B \f_i  \ .
\eeq
 Note that the matrix  $\left(  \smallmatrix - 1 & 8  \frac{n-3}{n-6}\\[5pt]
- 1 &   2\frac{n-6}{n-3} \endsmallmatrix \right)$ is invertible if and only if 
$$8   \frac{n-3}{n-6}  \neq 2 \frac{n-6}{n-3}  \qquad \Longleftrightarrow \qquad 4(n-3)^2 \neq (n-6)^2 \qquad 
\Longleftrightarrow\qquad  n \neq 4\ .$$
Hence, if $n > 4$,  for a given  $1 \leq i \leq n-2$   the  equations   \eqref{5.32terter} are  equivalent to   
\beq  \left( \begin{array}{c}  \frac{A_i}{ A^2}\\[5 pt]  J^\ell_{i}   D_\ell 
\end{array} \right) =  \left( \begin{array}{cc} - 1 & 8  \frac{n-3}{n-6}\\[5pt]
- 1 &  2 \frac{n-6}{n-3} \end{array}\right) ^{-1} \left(\begin{array}{c}  8    \frac{n-3}{n-6}  \\[5 pt]   2   \frac{n-6}{n-3}  \end{array}\right)B \f_i  = 
 \left( \begin{array}{c}  0\\[7pt] 1\end{array}\right) B \f_i \ . 
\eeq
 This means that   $B \f_i J^i_\ell = - D_\ell$ or, equivalently,  that  $  B d^c \f = - d D$.  
 From this it follows that   $0 = -  d^2 D = B d d^c \f  = B \o_o$. This  is  impossible, because $\o_o$ is non-degenerate and $B \neq 0$. From this contradiction,  Theorem \ref{nogo} follows.
   \par
   \medskip
    \section{Classification of   Kerr $4$-manifolds admitting Kerr type  backgrounds}
    \subsection{The fundamental constraints on   the underlying K\"ahler surface $(N, J, g_o)$}
      From now on we  focus on the $4$-dimensional case. Let  $M$ be a   Kerr $4$-manifold, equipped with a Kerr structure $\cM = (M = \cS\times \bR$ , $\pi: \cS \to N, (J, g_o),   \cH)$ and  the associated $\cM$-canonical optical structure  $\cQ_o = (\cW_o, [h_o]_\pm, \cK_o, \{g\})$.   The following theorem   provides    (topological and metrical)  conditions on the  $2$-dimensional K\"ahler manifold  $(N, J, g_o)$ that are necessary for  $M$  to admit a Kerr type background metric.  In the subsequent  section we prove  that  such necessary  conditions are also sufficient  and a detailed   description of all background metrics of a Kerr manifold is given.  According to Proposition \ref{flatness}, all such background metrics are locally isometric  to each other (they are all flat).  However,  as we will shortly see, they all have  different  presentations as compatible metrics  and correspond to different  canonical data.\par
 \medskip
      \begin{theo} \label{111}    Let $M$ be  a Kerr $4$-manifold with  Kerr structure $\cM$ and  $\cM$-canonical optical structure $\cQ_o$. The Ricci tensor of a metric $\h$  as in (a) of Theorem \ref{nogo}
is such that 
          \beq \label{the51} \Ric(X, Y) = 0\qquad \text{for any pair of vector fields}\ X, Y \in  \cW_o\  \eeq
    if and only if   the following  three conditions hold:
    \begin{itemize}[leftmargin = 22 pt]
    \item[(1)] The  Riemann surface $(N, J)$ admits a Riemannian metric $\wt g_o$ of constant Gaussian curvature $\k = \pm  1$, which is K\"ahler with respect to   the complex structure $J$ (i.e., in real terms, there exists 
 an atlas of real coordinate charts, each of them simultaneously isothermal   for $g_o$  and $\wt g_o$);
    \item[(2)] There exists a K\"ahler potential   $\f: N \to \bR$  for $g_o$ (i.e. such that $\o_o = d d^c \f$) satisfying
    \beq\label{42}  \operatorname{sign} \f = -   \k\ ,\qquad g_o + 2 \k \f\, \wt g_o = 0\ ;\eeq
    \item[(3)]   The Riemann surface $N$  is not compact and,  up to a  bundle automorphism of the trivial bundle  $\pi^\cS: M = \cS \times \bR \to \cS$, there is  a globally defined coordinate $r$ for the fiber of $\pi^\cS: M \to \cS$ 
    such that  
     \beq\label{form}  \s =  -  \frac{r^2}{4 B \f}- \frac{ B \f}{4}\ .\eeq
     \end{itemize}   
    \end{theo}
    
    The proof   is divided into two steps and it is given in the following subsections. We first show that a sub-system  of the equations corresponding to the   \eqref{the51} implies that  $\s$ has  locally  the form \eqref{form} for some  potential $\f$  for   the metric $g_o$. Second, we prove that the remaining equations  corresponding to  \eqref{the51} are satisfied if and only if  there is an open cover $\{\cV_A\}$ for $N$  and a  corresponding family of  local potentials  on the  sets $\cV_A$, each of them  satisfying \eqref{42} for a metric $\wt g_o|_{\cV_A}$ with constant scalar curvature $\k = \pm 1$.  From  these facts  it directly follows  that: (i) The metrics $\wt g_o|_{\cV_A}$  combine into a globally  defined metric of constant Gaussian curvature over $N$ with  the same isothermal coordinates of $g_o$;  (ii) There exists a function $\f: N\to \bR$, which  is well defined at each $p \in N$  by 
     $$\f(p) \= - \k \frac{g_o(X_p, X_p)}{2 \wt g_o(X_p, X_p)}$$
     for a  freely specifiable  non-trivial local vector field $X$   on a neighbourhood of   $p$; (iii)  for any open set $\cV_A$ of the above open cover  of $N$,  the restriction $\f|_{\cV_A}$   is  a potential  for $g_o|_{\cV_A}$ and, consequently, the  function $\f$  is a  potential for $g_o$, which is globally defined over $N$.  From (i) --(iii), the claims (1) -- (2)  follow.  Moreover: 
  \begin{itemize}[leftmargin = 15 pt]
 \item The K\"ahler manifold $(N, J, g_o)$ is not  compact (because it admits a global potential); 
 \item Since $\f$ and $\s$ are  globally defined on $M$,  the  relation \eqref{form} (which was originally shown to be valid just in appropriate coordinates sets) defines
 a global real function $r: M \to \bR$, that  provides a global  coordinate for the fiber of  $\pi^\cS: M = \cS \times \bR \to \cS$.
 \end{itemize}
 From these two observations, (3) follows. \par
 \smallskip
 We   also  stress the following two  facts: 
 \begin{itemize}[leftmargin = 22 pt]
 \item[(a)]   If $\cV \subset N$ is an open set with coordinates   $(x, y, v, r)$    on $\wc \pi^{-1}(\cV) = \cV \times \bR \times \bR$  with  $(x,y)$  selected in such a way  that  the coordinate components of  $\wt g_o$ have the classical  expressions    $\wt g_{o\,ij} =    \frac{4 \d_{ij}}{\left(1 + \k( x^2 + y^2)\right )^2}$, then  (2) holds if and only if  $\f|_{\cV}$ satisfies one of the following  pairs of conditions: 
     \beq \label{eqpot} 
     \begin{split} 
     & \Delta \f  + \frac{8 \f}{\left (1 + x^2 + y^2\right )^2} =0 \qquad \& \qquad \ \f < 0\ \qquad  \text{in case } \ \k = 1\qquad \text{or}\\
    &  \Delta \f  - \frac{8 \f}{\left (1 - x^2 - y^2\right )^2} =0\qquad  \& \qquad  \ \f > 0\  \qquad \text{in case } \ \k = -1\ .
    \end{split}
     \eeq
     \item[(b)]  For both possibilities for  $\k$, we have that     $\operatorname{sign}(\s) =  \operatorname{sign}(B \k)$ and the signature of $g$ is  mostly plus or mostly minus, according to  whether $\operatorname{sign}(B \k)$ is equal to $-1$  or $+1$.
        \end{itemize}           
These  facts will greatly assist in the discussions of the following sections.
          \par
    \smallskip
    \subsubsection{Step 1 of the Proof of Theorem \ref{111}} \label{section511}
  Consider an open cover $\{\cV_A\}_{A \in \cJ}$ of $N$  as in   Lemma \ref{thelemma51},   and pick an open set $\cV$ in this cover and   a (local) potential    $\f: \cV \to \bR$   for   $g_o$.  We adopt  the notational conventions of \S \ref{section321}. With no loss of generality,  we    also assume that  $\cV \subset N$  admits isothermal coordinates $(x, y)$ for  $g_o$ and  that $\big(E_1 \= \frac{\p}{\p x} , E_2 \= \frac{\p}{\p y}\big)$. In this way   the components 
   $J_i^j$,   $\o_{ij}$ and   $g_{ij}$  of the  complex structure, the K\"ahler form and the K\"ahler metric have the form
    \beq 
\begin{split}   (J_i^j) &= \left( \begin{matrix} 0 & -1\\1 & 0\end{matrix}\right)\ ,\\
 \big(\o_{ij} \big) & = \big(d d^c\f (E_i, E_j) \big) =  \big(- E_i( J_j^\ell \f_\ell) + E_j( J_i^\ell \f_\ell)  \big) = \left(\begin{array}{cc}  0 & \lambda \\
- \lambda & 0 
   \end{array}\right) \ ,
\\
 (g_{ij}) &= \left(g(E_i, E_j) \right)=  \left(\o_{i m} J^m_j\right) =   \lambda \,\d_{ij}, \qquad \text{where we set}\ \lambda \= \D \f\ .
\end{split}
\eeq
Finally,  as usual, we denote  by $(X_A) = (\wh E_1, \wh E_2, \ps_o , \qs_o)$ the adapted frame field on $\wc \pi^{-1}(\cV) \simeq \cV \times \bR^2$,   corresponding  to the  frame field $(E_1, E_2)$. \par
\smallskip  
By the results of \S \ref{sect423}, the condition $\Ric(\ps_o, \ps_o) = 0$ (which is a consequence of  \eqref{the51}) is satisfied on the open set $\wc \pi^{-1}(\cV)$ if and only if $\sigma = \s(x,y, v, u)$ has the form $\sigma =  C r^2 + \frac{1}{16 C}$ 
 for a nowhere vanishing function $C: \cV \to \bR$ and   $r \= 
 u + B v +  D(x, y)$, where $B$ is the  constant appearing in the Definition \ref{ansatz} and $D=D(x,y)$ is an arbitrary function. 
 Replacing the fibre coordinate $u$ by $r = 
 u + B v +  D(x, y)$ yields new coordinates $(x,y,v,r)$ on  $\wc \pi^{-1}(\cV)$, in which the vector fields $\ps_o$,  $\qs_o$ take the form
 $\ps_o = \frac{\p}{\p r}$,  $\qs_o = \frac{\p}{\p v} 
 + B \frac{\p}{\p r}$ and  $\s  =    C r^2 + \frac{1}{16 C}$ becomes a 
 function of the  new (independent) variables  $(x,y, v, r)$.
 \par
 \smallskip
The  computations  in \S \ref{sect425} (which  hold in {\it any} dimension) show that $\Ric(\ps_o, \wh E_i) = 0$, $i = 1,2$, hold true if and only if  the system 
\eqref{4.26} is satisfied.  Since we are now in the case  $n = 4$ and  we are in a  coordinates system, in which  $\s$ has  no summand that is linear with respect to the coordinate   $r$ (meaning that we may  assume  $D = \frac{C_1}{2 C_2}$ to be   zero), the first two equations of  \eqref{4.26} are equivalent each other and both equivalent to 
\beq \label{4.26*}
 J^\ell_{i}  \left(\frac{C_\ell}{C^2}   - 
  4 B  \f_\ell \right)  = 0\quad \Bigg( \Longleftrightarrow \quad d\left( \frac{1}{C} +  4  B \f\right) = 0 \ \Longleftrightarrow \   \frac{1}{C} =  - 4  B \f + \text{const.}\ \Bigg)
\eeq
Since the potentials are determined up to constants, we may assume   $ \frac{1}{C} =  - 4  B \f$ and  conclude that 
$\Ric(\ps_o, X) = 0$ for any $X \in \cW_o$ if and only if $\s$ has the form \eqref{form} for an appropriate   K\"ahler  potential $\f$ on $\cV$.\par
\smallskip
\subsubsection{Step 2 of the Proof of Theorem \ref{111}} 
Let us now assume that $\s$ has the form \eqref{form} and hence that $\Ric(\ps_o, X) = 0$ for any $X \in \cW_o$. Then    \eqref{the51} holds if and only if   $\Ric_{ij} = \Ric(\wh E_i, \wh E_j) = 0$ for any  $1 \leq i, j \leq 2$. 
Writing   the  curvature     in terms of   the Christoffel symbols  \eqref{519} -- \eqref{6.17}) (see \eqref{Riemannform},  \eqref{Riemannform1}) and    recalling that the $E_i = \frac{\p}{\p x^i}$ commute,  we get the following  (here, as in the previous section, $\lambda \= \D \f$): \\
\phantom{aaa} \scalebox{0.8}{\vbox{
 \begin{align} 
\nonumber   \Ric_{11}\  & {=}   \left(E_m +J_m^\ell \f_\ell \qs_o\right)   (\GGa 1 1 m) -  \left(E_1 + \f_2 \qs_o\right)(\GGa m 1 m) -  \GGa m 1  \ell \GGa 1  \ell  m    -  \GGa m 1  {\ps_o} \GGa 1 {\ps_o}  m -  \GGa m 1  {\qs_o} \GGa 1 {\qs_o}  m +\\
\nonumber   &  +   \GGa 1 1 \ell  \GGa m \ell  m +   \GGa 1 1 {\ps_o}  \GGa m {\ps_o}  m +   \GGa 1 1 {\qs_o}  \GGa m {\qs_o}  m  -  \lambda   \GGa {\qs_o} 1 2    +  \ps_o (\GGa 1 1 {\ps_o})  +  \qs_o (\GGa 1 1 {\qs_o}) -  \GGa {\ps_o} 1  \ell \GGa 1  \ell {\ps_o} - \\
  & - \GGa {\qs_o} 1  \ell \GGa 1  \ell {\qs_o} \ , \\
\nonumber   \Ric_{12}\  &{=}   \left(E_m +   J_m^\ell\f_\ell \qs_o\right)   (\GGa 1 2 m) -  \left(E_1 +  \f_2 \qs_o\right)(\GGa m 2 m) -  \GGa m 2  \ell \GGa 1  \ell  m    -  \GGa m 2  {\ps_o} \GGa 1 {\ps_o}  m  -  \GGa m 2  {\qs_o} \GGa 1 {\qs_o}  m +\\
 \nonumber  & +   \GGa 1 2 \ell  \GGa m \ell  m +   \GGa 1 2 {\ps_o}  \GGa m {\ps_o}  m  +   \GGa 1 2 {\qs_o}  \GGa m {\qs_o}  m   -\lambda  \GGa {\qs_o} 2 2   +   \ps_o (\GGa 1 2 {\ps_o})  +   \qs_o (\GGa 1 2 {\qs_o})  -  \GGa {\ps_o} 2  \ell \GGa 1  \ell {\ps_o}- \\
  &  -  \GGa {\qs_o} 2  \ell \GGa 1  \ell {\qs_o} \ ,
\\
\nonumber   \Ric_{22}\ & {=}   \left(E_m +  J_m^\ell\f_\ell \qs_o\right)   (\GGa 2 2 m) -   \left(E_2 -   \f_1 \qs_o\right)(\GGa m 2 m) -  \GGa m 2  \ell \GGa 2  \ell  m    -  \GGa m 2  {\ps_o} \GGa 2 {\ps_o}  m  -  \GGa m 2  {\qs_o} \GGa 2 {\qs_o}  m +\\
\nonumber   &+   \GGa 2 2 \ell  \GGa m \ell  m +   \GGa 2 2 {\ps_o}  \GGa m {\ps_o}  m+   \GGa 2 2 {\qs_o}  \GGa m {\qs_o}  m   + \lambda   \GGa {\qs_o} 2 1   +  \ps_o (\GGa 2 2 {\ps_o})  +  \qs_o (\GGa 2 2 {\qs_o})  -  \GGa {\ps_o} 2  \ell \GGa 2  \ell {\ps_o} - \\
& -  \GGa {\qs_o} 2  \ell \GGa 2  \ell {\qs_o} \ .
 \end{align}
 }
 }\\
Since $\ps_o = \frac{\p}{\p r}$,  $\qs_o =\frac{\p}{\p v} + B \frac{\p}{\p r}$,  and no Christoffel symbol depends on  $v$,  the differential operator $\qs_o$ acts on the $\GGa A B C$  just as  $   B\frac{\p}{\p r} =   B \ps_o$. Thus   $\Ric_{ij} = \Ric(\wh E_i, \wh E_j) = 0$ 
reduce to\\
\phantom{aaa} \scalebox{0.8}{\vbox{\begin{align} 
\nonumber  \Ric_{11}= &  \left(E_m + J_m^\ell \f_\ell   B\frac{\p}{\p r}\right)   (\GGa 1 1 m) -  \left(\frac{\p}{\p x} + B \f_y \frac{\p}{\p r}\right)(\GGa m 1 m) -  \GGa m 1  \ell \GGa 1  \ell  m    -  \GGa m 1  {\qs_o} \GGa 1 {\qs_o}  m +\\
 \label{6.74} &\hskip 1.5 cm  +   \GGa 1 1 \ell  \GGa m \ell  m+   \GGa 1 1 {\qs_o}  \GGa m {\qs_o}  m   - \lambda   \GGa {\qs_o} 1 2  + B  \frac{\p}{\p r} (\GGa 1 1 {\qs_o}) -  \GGa {\qs_o} 1  \ell \GGa 1  \ell {\qs_o} = 0 \ ,\\
\nonumber   \Ric_{12}= & \left(E_m + J_m^\ell E_\ell(\f) B	\frac{\p}{\p r}\right)   (\GGa 1 2 m) -  \left(\frac{\p}{\p x} + B \f_y \frac{\p}{\p r}\right)(\GGa m 2 m) -  \GGa m 2  \ell \GGa 1  \ell  m    -  \GGa m 2  {\qs_o} \GGa 1 {\qs_o}  m +\end{align}}}\\
 \phantom{aaa} \scalebox{0.8}{\vbox{\begin{align} 
 \label{6.75}  &\hskip 1.5 cm  +   \GGa 1 2 \ell  \GGa m \ell  m+   \GGa 1 2 {\qs_o}  \GGa m {\qs_o}  m   -  \lambda  \GGa {\qs_o} 2 2   + B\frac{\p}{\p r} (\GGa 1 2 {\qs_o}) -  \GGa {\qs_o} 2  \ell \GGa 1  \ell {\qs_o} = 0\ ,
\\
\nonumber   \Ric_{22}= &  \left(E_m + J_m^\ell E_\ell(\f)  B	 \frac{\p}{\p r}\right)   (\GGa 2 2 m) -  \left(\frac{\p}{\p y} - B \f_x \frac{\p}{\p r}\right)(\GGa m 2 m) -  \GGa m 2  \ell \GGa 2  \ell  m    -  \GGa m 2  {\qs_o} \GGa 2 {\qs_o}  m +\\
\label{6.76}  &\hskip 1.5 cm  +   \GGa 2 2 \ell  \GGa m \ell  m+   \GGa 2 2 {\qs_o}  \GGa m {\qs_o}  m   + \lambda   \GGa {\qs_o} 2 1   + B \frac{\p}{\p r}  (\GGa 2 2 {\qs_o}) -  \GGa {\qs_o} 2  \ell \GGa 2  \ell {\qs_o} = 0\ .
 \end{align}
 }
 }\\
Now,  plugging the condition  $\wt \b =  B$ into  the expressions  \eqref{519} -- \eqref{6.17},  we get that the only   non-zero Christoffel symbols  $\GGa A B C$ of  $\h $ are\\
\scalebox{0.8}{\vbox{
     \begin{align}
 & \GGa 1 1 1 =\Ga 1 1 1
+ \frac{1}{2 \s }   \s_x    + B \frac{\s_t \f_y }{2  \s } &
 & \GGa 1 1 2 =\Ga 1 1 2
 -   \frac{1}{2 \s}  \s_y  + B  \frac{\s_t \f_x }{2 \s }\\
  & \GGa 1 2 1  =  \GGa 2 1 1  =\Ga 1 2 1
  + \frac{1}{2 \s }  \s_y  - B  \frac{\s_t\f_x }{2 \s }  &
 & \GGa 1 2 2 = \GGa 2 1 2 = \Ga 1 2 2
+ \frac{1}{2 \s }  \s_x  + B  \frac{\s_t \f_y }{2 \s }
 \end{align}}}\ \\
 \scalebox{0.8}{\vbox{
     \begin{align}
 & \GGa 2 2 1 =\Ga 2 2 1
 -   \frac{1}{2 \s}  \s_x  - B   \frac{\s_t \f_y }{2 \s } &
 & \GGa 2 2 2 =\Ga 2 2  2
+ \frac{1}{2 \s }  \s_y  - B \frac{\s_t \f_x }{2 \s } \\
 &  \GGa 1 1 {\qs_o}  = - \lambda   \s_t \ ,\qquad
 \GGa 1 2 {\qs_o}  =   -  \frac{  \lambda}{2} \ , &&
  \GGa 2 1 {\qs_o}  =   \frac{\lambda}{2} \ ,\qquad   \GGa 2 2 {\qs_o}  =  -\lambda  \s_t\ , \\
 & \GGa 1 {\ps_o} 1 = \GGa {\ps_o} 1 1 =   \frac{\s_t }{2 \s} \ ,&&
  \GGa 1 {\ps_o} 2 = \GGa {\ps_o} 1 2 =     \frac{1}{4 \s} \ , \\
 & \GGa 2 {\ps_o} 1 = \GGa {\ps_o} 2 1 =  -   \frac{1}{4 \s}  \ ,&&
  \GGa 2 {\ps_o} 2 = \GGa {\ps_o} 2 2 =  \frac{\s_t}{2 \s} \ , \\
   & \GGa 1 {\qs_o} 1 = \GGa {\qs_o} 1 1 =   B \frac{\s_t }{2 \s} \ ,&&
  \GGa 1 {\qs_o} 2 = \GGa {\qs_o} 1 2 =   B  \frac{1}{4 \s} \ ,\\
 &  \GGa 2 {\qs_o} 1 = \GGa {\qs_o} 2 1 =   - B  \frac{1}{4 \s}   \ ,&&
  \GGa 2 {\qs_o} 2 = \GGa {\qs_o} 2 2 = B \frac{\s_t}{2 \s} \ .
  \end{align}
  }}\\
  Using this,   a sequence of  tedious but straightforward computations and simplifications lead to  the following
\begin{lem}\label{lemmappend} The equation \eqref{6.75} is identically satisfied, while \eqref{6.74} and \eqref{6.76} are both equivalent  to the equation 
    \begin{multline} \label{commonequation}
   \s^2 \Ric^N_{11}        -\frac{1}{2} \s_{xx} \s + \frac{1}{2} \s_x^2       -\frac{1}{2} \s _{yy} \s +\frac{1}{2} \s_y^2   
  -   \frac{ B}{4} (\f_{xx} + \f_{yy})  \s    -\\
  - B (\f_{xx} + \f_{yy}) \s^2\s_{rr} - \frac{ B^2}{2} \s \s_{rr} \left(  (\f_y)^2   + (\f_x)^2\right) +  \\
+   \frac{B}{2} \s_y\f_y    +   \frac{B}{2}\s_x   \f_x   +    \frac{ B^2}{8} \left( (\f_x)^2  + (\f_y)^2\right)  = 0
 \end{multline}
where   $\Ric^N$ denotes the Ricci tensor of the K\"ahler metric $g_o$ of $N$ and  $\Ric^N_{ij} \= \Ric(E_i, E_j)$. 
\end{lem}
We give the  details of the  proof of this lemma  in Appendix \ref{appendixA}.   Plugging the expressions
\beq  \noindent 
\begin{split}
& \s =  -  \frac{1}{4 B \f }r^2 - \frac{ B \f}{4}\ ,\qquad \s_r =   - \frac{1}{{2 B\f}} r\ ,\qquad \s_{rr} =  -  \frac{1}{{ 2 B\f}}\ ,\\
& \s_x =    \frac{\f_x}{{ 4 B\f}^2}r^2  -  \frac{B\f_x}{4}	 \ ,\qquad \s_y =   \frac{\f_y}{{ 4 B\f}^2}r^2  - \frac{B\f_y}{4}\ ,\\
& \s_{xx} =    \frac{\f_{xx}}{ 4 B\f^2}r^2  - \frac{\f^2_x}{2 B\f^3} r^2  -  \frac{ B\f_{xx}}{4}\ ,\qquad  \s_{yy} =   \frac{\f_{yy}}{ 4 B\f^2}r^2  -  \frac{\f^2_y}{2 B\f^3} r^2  -  \frac{ B\f_{yy}}{4}\ .
\end{split}
\eeq
into \eqref{commonequation},  we  
see that the left hand side   is an even  polynomial of degree $4$ in the variable $r$. The equation is therefore satisfied if and only if the coefficients of the monomials $r^4$, $r^2$ and $ r^0$ are identically vanishing. These three conditions  are
\beq
\begin{split}
&  \frac{1}{16B^2\f^2}\left( \Ric^N_{11}         -  \frac{\f^2_x + \f_y^2}{2\f^2}  
+   \frac{\f_{xx} + \f_{yy} }{ \f}  \right)   = 0\ ,\ 
    \frac{1}{8}\left(\Ric^N_{11}        -   \frac{\f^2_x + \f_y^2}{2 \f^2}     +  \frac{\f_{xx} + \f_{yy}}{ \f} \right) = 0\ ,  \\
&  \frac{\f^2}{16^2}\bigg( \Ric^N_{11}    -  \frac{\f^2_x + \f_y^2}{ 2 \f^2}   
+    \frac{ \f_{xx}  + \f_{yy}}{\f}  \bigg) = 0\ ,
\end{split}\eeq
which are  manifestly   equivalent   each other and all   to  the single equation 
\beq \label{THEEQ} \Ric^N_{11}        -   \frac{\f^2_x + \f_y^2}{2 \f^2}     +  \frac{\Delta \f}{ \f} = 0\  \bigg( \Longleftrightarrow \  - \frac{1}{2} \D(\log\D \f)   + \frac{\Delta \f}{2 \operatorname{sign}(\f) |\f|} + \frac{1}{2}  \D \log|\f|  = 0\ \bigg)\eeq
(here we used  the classical  formula  for curvatures in isothermal coordinates).  Introducing the auxiliary function  $\psi \= \frac{\D \f}{|\f|} > 0$,   the equation  \eqref{THEEQ} becomes equivalent to the system
\beq \label{THESYS}
\D \f - \operatorname{sign}(\f) \f \psi\ =0,\qquad
 - \frac{1}{2} \D(\log \psi) +  \operatorname{sign}(\f)  \frac{1}{2} \psi   = 0\ .
\eeq
If we  consider the metric $\wt g_o|_{\cV}$ on $\cV$ with components $(\wt g_{o\,ij})  = \left( \begin{array}{cc} \frac{\psi}{2} & 0 \\0 & \frac{\psi}{2} \end{array}\right)$,  the second equation in \eqref{THESYS} simply means that, in the regions where $\operatorname{sign}(\f)$ is constant,   the Gaussian curvature $\k$ of $\wt g_o|_{\cV}$ is constant and equal to $\k = - \sign(\f) = \mp 1$.  At the same time, since  $(x, y)$ are isothermal coordinates for both metrics, the first equation in \eqref{THESYS}   can be re-written as  $g_o|_{\cV} +  2 \k \f \wt g_o|_{\cV}$.  This yields the above described results for the  second step. 
\par
\medskip
\subsection{Classification of  backgrounds of Kerr type} 
As we previously  planned,  we can now prove the main  classification result of this section. 
   \begin{cor}\label{333} Let $M$ be a $4$-dimensional Kerr manifold, equipped  with a Kerr structure $\cM$ and the corresponding $\cM$-canonical optical structure $\cQ_o$,  as in Theorem \ref{111}. Such a Kerr manifold admits a background metric $\h$ of Kerr type if and only if the following  hold:
   \begin{itemize}[leftmargin = 15pt]
   \item[(1)] The  K\"ahler surface $(N, J, g_o)$ satisfies: (a) it  is  non-compact, (b) it admits  a Riemannian metric $\wt g_o$  of constant Gaussian curvature $\k =   \pm 1$, which is K\"ahler for the complex manifold $(N, J)$ and (c)   there exists a  global K\"ahler potential   $\f: N \to \bR$  for $g_o$ satisfying
    \beq\label{42*}  \operatorname{sign} \f = - \k\ ,\qquad g_o + 2 \k \f\, \wt g_o = 0.\eeq
  \item[(2)]  The $\bR$-bundle $\pi: \cS \to N$  is trivial and, up to a  bundle automorphism of $M$, there is  a standard coordinate $r$ for the fiber of $\pi^\cS: M = \cS \times \bR\to \cS$ 
    such that   $  \s =  -  \frac{r^2}{4 B \f}- \frac{ B \f}{4}$.
     \end{itemize}  
  If $M$ satisfies the above conditions, for any  $\f$  satisfying $\o_o = d d^c \f$ and  \eqref{42}, there is  a  family of   background  metrics $\h$, parameterised by  $B$,   with     $\s$  as in (2).
    \end{cor}
\begin{pf}  By Proposition \ref{flatness}, any background metric  $\h$ on $M$ is flat and  satisfies \eqref{the51}. Thus, the necessity of  the above conditions is a  consequence of Theorem \ref{111}. Conversely, assume  that there are $\wt g_o$ and $ \f$  as in the hypotheses and consider an open cover $\{\cV_A\}_{A \in \cJ}$  as in Lemma \ref{thelemma51} and on each  corresponding trivialisable set $\cU = \wc \pi^{-1}(\cV) = \cV \times \bR \times \bR$,  let $\s: \cV_A \times \bR \times \bR \to \bR$ be as \eqref{form}.  By condition (2) these locally defined functions combine into  a global function $\s: M \to \bR$. Let $\cV$ be  one of the    sets of the  open cover  $\{\cV_A\}$  and denote by $(\wh E_1, \wh E_2, \ps_o, \qs_o)$ the  adapted frame field on $\wc \pi^{-1}(\cV) $ in which  $\wh E_1$, $\wh E_2$ are  lifts of the two coordinate vector fields,  given  by $g_o$-isothermal coordinates $(x, y)$ on $\cV$. By  Theorem \ref{111}, the Ricci tensor of a  metric $\h$ as in (a) of Theorem \ref{nogo} with $\wt \b = B \neq 0$  being constant and $\s$  determined as above,  the equation \eqref{the51} vanishes on the section $\cW_o$.  Therefore, in order to prove that such $\h$ is  Einstein on each open set $\wc \pi^{-1}(\cV)$, it  suffices  to  show  that on any such open set it satisfies the Einstein equations that are not included in  \eqref{the51}, that is   the equations
\beq \label{the525}  \Ric(\qs_o, \qs_o) = \Ric(\qs_o, \ps_o) = \Ric(\qs_o, \wh E_1) = \Ric(\qs_o, \wh E_2) = 0\ . \eeq
 Let us  start with  the first two equations,  $  \Ric(\qs_o, \qs_o) =0$ and  $  \Ric(\qs_o, \ps_o) =0$. 
  By definition, 
   $\Ric(\qs_o, \qs_o) =  \RR m {\qs_o} {\qs_o} m +  \RR {\ps_o} {\qs_o} {\qs_o} {\ps_o}$ and $ \Ric(\ps_o, \qs_o) =  \ \RR m {\ps_o} {\qs_o} m +  \RR {\qs_o} {\ps_o} {\qs_o} {\qs_o} $.   Expressing  the  components of the Riemann curvature in terms of the Christoffel symbols  \eqref{519} -- \eqref{6.17} (with $\wt \b \= B$!) and   omitting all terms that are manifestly zero,  the two equations can be written as \\ \phantom{aaaaaa} \scalebox{0.8}{\vbox{
\begin{align}
\nonumber \Ric_{{\qs_o}{\qs_o}}\  & \ {=} -  \qs_o (\GGa m {\qs_o} m) -  \GGa m {\qs_o}  \ell \GGa  {\qs_o} \ell m  = - B^2  \ps_o (\GGa m {\ps_o} m) -  B^2 \GGa m {\ps_o}  \ell \GGa  {\ps_o} \ell m =  \\
\nonumber&=  B^2\left(-\ps_o \left (     \frac{J_{m}^m }{4 \s}+ \frac{ \ps_o(\s)\d_{m}^m}{2 \s}\right )-\left (    \frac{J_{m}^{\ell} }{4 \s}
+ \frac{ \ps_o(\s)\d_{m}^{\ell}}{2 \s}\right ) \left(    \frac{J_{\ell}^m }{4 \s}+ \frac{ \ps_o(\s)\d_{\ell}^m}{2 \s} \right)\right) = \\
\nonumber  &= B^2\left(- 2 \frac{ \ps_o(\ps_o(\s))}{2 \s}+2  \frac{ \ps_o(\s)^2}{4 \s^2}+ \frac{2}{16\s^2}\right) =
\\
\label{6.25ter} & =   \frac{B^2 }{   \s^2}\left (-2\s\ps_o(\ps_o(\s)) +\ps_o(\s)^2+ \frac{1 }{4}\right ) =0\ ,\\
 \nonumber \Ric_{{\ps_o}{\qs_o}}\  &\ {=}  -  \ps_o (\GGa m {\qs_o} m)-  \GGa m {\qs_o}  \ell \GGa  {\ps_o} \ell m \\
  \nonumber  &=  -  B \ps_o \left ( \frac{\ps_o(\s)}{ \s} \right ) - B \left ( \frac{J_{m}^{\ell} }{4 \s}+ \frac{ \ps_o(\s)\d_{m}^{\ell}}{2 \s}  \right)\left (   \frac{J_{\ell}^m }{4 \s}+ \frac{ \ps_o(\s)\d_{\ell}^m}{2 \s}\right ) = \\
 \label{6.23ter*}& =  -\frac{ B  }{ 2\s^2}\left (- 2\s\ps_o(\ps_o(\s)) +\ps_o(\s)^2+ \frac{1 }{4}\right ) =0\ .
  \end{align}
  }}\\
   Since they are both  equivalent to \eqref{main}, they are  identically satisfied because  $\h$ satisfies $\Ric(\ps_o, \ps_o)= 0$ by \eqref{the51}. \par
\smallskip
For what concerns the equations
$ \Ric(\wh E_i, \qs_o) = 0$, $i = 1,2$,  using once again the expressions \eqref{Riemannform} and \eqref{Riemannform1} and neglecting the terms  that  are trivial, according to  \eqref{519} -- \eqref{6.17}, the two equations reduce to
\begin{multline}  \Ric(\wh E_i, \qs_o) \  {=} \ \RR m i {\qs_o} m + \RR {\ps_o} i {\qs_o} {\ps_o} + \RR {\qs_o} i {\qs_o} {\qs_o}  = \\
=\wh E_m (\GGa i {\qs_o} m) -  \wh E_i (\GGa m {\qs_o} m) -  \GGa m {\qs_o}  \ell \GGa i  \ell  m   +   \GGa i {\qs_o} \ell  \GGa m \ell  m   -  \GGa {\ps_o} {\qs_o}  {\qs_o} \GGa i  {\qs_o}  {\ps_o} = 0\ .
 \label{6.22}
\end{multline}
 Using once again  the explicit expressions for the Christoffel symbols \eqref{519} -- \eqref{6.17} with $\wt \b \= B$, 
a straightforward check shows that  these two equations are equivalent to the equations   $\Ric(\wh E_i, \ps_o) = 0$, which we discussed in \S \ref{section511}. Since \eqref{the51} holds,   we conclude  that also the equations \eqref{6.22} are  identically satisfied.
\end{pf} 
\par
\medskip
 \section{Kerr families of Ricci flat gravitational fields}
 As in the previous section,  $M$ is a $4$-dimensional Kerr manifold,  with a Kerr structure $\cM = (M = \cS\times \bR$ , $\pi: \cS \to N, (J, g_o),   \cH)$ and   associated $\cM$-canonical optical structure  $\cQ_o = (\cW_o, [h_o]_\pm, \cK_o, \{g\})$.  The  next theorem together with   Corollaries \ref{nogo1} and  \ref{333}  represent our main result.  In \S \ref{expli-sol} we will show that it can be re-stated as in Theorem \ref{teoremone}.\par
      \begin{theo} \label{aaa}    Given   $M$, $\cM$  and $\cQ$, assume  that  $M$ admits a background metric $\h$ of Kerr type, determined   by a constant  $B\neq 0$ and functions  $\f: N \to \bR$,   $\s: M \to \bR$ as in  Corollary \ref{333}. Then  there exists  a constant $m_o \neq 0$ and a  Ricci flat compatible metric $g^{(m = m_o)} \neq \h$ of the form \eqref{Kerr-Schild}    in a connected open subset $M' = \cS \times I$, $I \subset \bR$,  of $M$ if and only if the following two conditions hold:  
     \begin{itemize}[leftmargin = 20pt]
     \item[(1)]  $B = -1$; 
     \item[(2)]  Using  the  coordinate $r$ for the fiber of $\pi^\cS: M = \cS \times \bR \to \cS$,  so that    $p_o = \frac{\p}{\p r}$  and $\s$ takes the form \eqref{form},   the function   $\wt \b_o$ takes  the form 
                \beq \label{lella}   \wt \b_o = \frac{\cappa r}{r^2 + \f^2}\ ,\eeq
      for some  constant $\cappa \neq 0$.
      \end{itemize}
      In this case,  all metrics  $g^{(m)}$ of the form  \eqref{Kerr-Schild},  with  $B = -1$ and  $\wt \b_o$  as  in \eqref{lella},  are   Ricci flat and constitute a Kerr family of gravitational fields.
           \end{theo}
 The proof is essentially just a straightforward  consequence   of  imposing the  condition   $\Ric^{(m)} = 0$ on the Ricci tensor of a  metric $g^{(m)}$
 having the form \eqref{Kerr-Schild} with the  same constant $B$ and the  same function $\s$ of   the  background  $\h$. Considering coordinates $(x,y, v, r)$ and  an  adapted frame field  $\cB = \bigg(\wh E_1, \wh E_2,  \ps_o = \frac{\p}{\p r}, \qs_o = \frac{\p}{\p v} + B \frac{\p}{\p r}\bigg)$  in an open set $\wc \pi^{-1}(\cV_A) \subset M$ of Lemma \ref{thelemma51}, with  $E_1 = \frac{\p}{\p x}$, $E_2 = \frac{\p}{\p y}$,  the   Ricci flatness  condition is   equivalent to the three  following systems of 
 equations:  
 \begin{align}
\label{62bis}  & \Ric^{(m)}(\ps_o, \ps_o) = 0 \ ,\qquad \Ric^{(m)}(\ps_o, \wh E_i) = 0\ ,\\
\label{63bis}  & \Ric^{(m)}(\ps_o, \qs_o) = 0\ ,\qquad  \Ric^{(m)}(\wh E_i, \wh E_j) =  0\ ,\\
 \label{64bis} & \Ric^{(m)}(\wh E_i, \qs_o) = 0\ ,\qquad \Ric^{(m)}(\qs_o, \qs_o) = 0\ .
 \end{align}
The   proof  is divided  into three steps, which we are going to  perform in  \S \ref{firststep1}, \ref{secondstep2} and \ref{thirdstep3}, respectively.    First, we  show that, since $\h = g^{(m = 0)}$ is  a background metric,  the equations \eqref{62bis} are identically satisfied by any  metric of the form \eqref{Kerr-Schild} (regardless of what are the constant $B$ and the function  $\wt \b_o$). Second, we  prove that  there exists an $m_o \neq 0$ such that $g^{(m_o)}$  satisfies    \eqref{63bis}  if and only if   $\wt \b_o$ has the form  \eqref{lella} for a freely specifiable function  $\cappa = \cappa(x,y)$ on the open set  $\cV_A \subset N$.  In the  last step, we show that, assuming that  \eqref{63bis} is satisfied (and hence  that $\wt \b_o$ has the form \eqref{lella} for some function $\cappa$), the  remaining equations \eqref{64bis}  are  satisfied if and only if  $\cappa = \cappa(x,y)$ is   constant  and   $B = -1$.  This will prove the first claim and imply the last claim  as well, because it will show   that  all metrics $g^{(m)}$ of the form \eqref{Kerr-Schild} satisfy (1) and (2) if and only if those conditions are satisfied by just one value $m = m_o \neq 0$. \par
\medskip
 \subsection{First step}\label{firststep1}
 From now on, we  use the notation  $\GGa A B C$ and $\wt {\GGa A B C}$ for  the Christoffel symbols with respect to the adapted frame field  $\cB$ of the background  metric $\h = g^{(m = 0)}$ and of the metric $g^{(m_o)}$, respectively. We  also
  denote  $\HH A B C \=\wt{ \GGa A B C} - \GGa A B C$,  so that we may write that $\wt {\GGa A B C} = \GGa A B C + \HH A B C$. 
  From  \eqref{519} -- \eqref{6.17} we have  (recall  that $\qs_o(\wt \b_o) = B \ps_o(\wt \b_o)$!)
   \begin{align}
  \label{519*}  & \HH i j m = 0\\
   \label{6.12*} &\HH i j {\ps_o}  =    m_o g_{ij}   \wt \b_o   \ps_o(\s)  \ ,\quad
 \HH i j {\qs_o}  = 0\ ,\\
\label{cli*} & \HH i {\ps_o} m = \HH {\ps_o} i m =  0 \ ,\qquad  \HH i {\ps_o} {\ps_o} =   \HH  {\ps_o} i {\ps_o} =
 \HH i {\ps_o} {\qs_o} =  \HH {\ps_o} i  {\qs_o} = 0 \ ,\\
\nonumber  & \HH  i {\qs_o} r =  \HH {\qs_o} i r =   m_o  \frac{J^r_i}{4 \s} \wt \b_o\ ,\\
 \label{via*} & \hskip 2 cm  \
 \HH i {\qs_o} {\ps_o} = \HH  {\qs_o} i {\ps_o} =   \frac{m_o}{2} \wt \b_{oi} - \frac{m_o}{2}J^\ell_i \f_\ell \ps_o(\wt \b_o)   \ ,\ \  \HH i {\qs_o} {\qs_o} =  \HH  {\qs_o} i {\qs_o} = 0\ ,
\\
 &  \HH {\ps_o}  {\ps_o}  m  =   \HH {\ps_o}   {\ps_o} {\ps_o}  =   \HH {\ps_o}   {\ps_o} {\qs_o}  =  0 \ ,\qquad    \\
   &  \HH {\ps_o}  {\qs_o}  m  = \HH {\qs_o}  {\ps_o}  m  =  0 \ ,\qquad
 \HH {\ps_o}   {\qs_o} {\ps_o}  =  \HH {\qs_o}   {\ps_o} {\ps_o}  =  \frac{m_o}{2 }\ps_o(\wt \b_o) \ ,\qquad \HH {\ps_o}   {\qs_o} {\qs_o} =  \HH {\qs_o}   {\ps_o} {\qs_o}   = 0 \ ,\\
\nonumber & \HH {\qs_o}  {\qs_o}  r =   -  m_o \frac{g^{rk}}{4\s} \wt \b_{ok} - m_oB  \frac{g^{rk}}{4\s}J_k^\ell  \f_\ell \ps_o(\wt \b_o)  \ ,\\
  \label{6.17*}  &  \hskip 2 cm \HH {\qs_o}   {\qs_o} {\ps_o}  = \frac{m_o B}{2}	 \ps_o(\wt \b_o) +  m_o \frac{\ps_o(\wt \b_o)  (B + m_o \wt \b_o)}{2} \ ,\ 
 \HH {\qs_o}   {\qs_o} {\qs_o}  =    -  \frac{m_o}{2 }\ps_o(\wt \b_o) \ .
  \end{align}
  \par
  \medskip
 As announced above,  in this step we just  observe that  the  equations  $\Ric^{(m_o)}(\ps_o, \ps_o) = \Ric^{(m_o)}(\ps_o, \wh E_i) = 0$ are  identically satisfied for any choice of $\wt \b_o$ and $B$.  Indeed,  after  expanding $\Ric^{(m_o)}(\ps_o, \ps_o)$ in terms of the Christoffel symbols   $\wt {\GGa A B C} $, one  can immediately see  that  $\Ric^{(m_o)}(\ps_o, \ps_o)$ is independent of the parameter $m_o$ and is therefore identically equal to $0$,  because  $\h = g^{(m = 0)}$ is Ricci flat.  For instance, using the fact that $c^r_{mj} \equiv 0$, we immediately get 
  \begin{align}
\nonumber  \Ric^{(m_o)}(\ps_o, \wh E_i)   &{=}    \wh E_m (\GGa i {\ps_o} m)+\xcancel{\wh E_m (\HH i {\ps_o} m)}  -  \wh E_i (\GGa m {\ps_o} m) -  \xcancel{ \wh E_i (\HH m {\ps_o} m)} \\ 
  \nonumber &-\left (\GGa m {\ps_o}  \ell +\xcancel{\HH m {\ps_o}  \ell }\right ) \left(\GGa i  \ell  m+\xcancel{\HH i  \ell  m} \right ) +
 \left ( \GGa i {\ps_o} \ell+\xcancel{\HH i {\ps_o} \ell }  \right ) \left (\GGa m  \ell  m+\xcancel{\HH m  \ell  m} \right ) = \\
   \nonumber &= \Ric^{(m = 0)}(\ps_o, \wh E_i) =0
\end{align}
Similar computations hold for $\Ric^{(m_o)}(\ps_o, \ps_o) = 0$. 
\par
\medskip
    \subsection{Second step} \label{secondstep2}
  Let us expand  the left hand side of the first equation  in   \eqref{63bis} (i.e.   the Ricci component $\Ric^{(m_o)}(\ps_o, \qs_o) $)   in terms of the Christoffel symbols $\wt {\GGa A B C} = \GGa A B C + \HH A B C$. Using the Ricci flatness of $\h$,  the   equation  reduces to 
  \begin{align}
\nonumber   \Ric^{(m_o)}({\ps_o}, {\qs_o})\ & {=}
  -  \ps_o (\GGa r {\qs_o} r+\xcancel{\HH r {\qs_o} r}) -  \left (\GGa r {\qs_o}  \ell +\HH r {\qs_o}  \ell \right )\left (\GGa  {\ps_o} \ell r +\xcancel{\HH  {\ps_o} \ell r} \right )+ \\
\nonumber  &+ \left (\GGa {\ps_o} {\qs_o} {\ps_o}+\HH {\ps_o} {\qs_o} {\ps_o}\right )\left (  \GGa r {\ps_o}  r+\xcancel{\HH r {\ps_o}  r}\right )    -    \ps_o (\GGa {\qs_o} {\qs_o} {\qs_o}+\HH {\qs_o} {\qs_o} {\qs_o}) =\\
\nonumber  &=\xcancel{\Ric_{{\ps_o}{\qs_o}}}-  \HH r {\qs_o}  \ell \GGa  {\ps_o} \ell r + \HH {\ps_o} {\qs_o} {\ps_o}  \GGa r {\ps_o}  r   -    \ps_o (\HH {\qs_o} {\qs_o} {\qs_o})=
\end{align}
\begin{align}
\nonumber  &=-m_o  \frac{J^{\ell}_r}{4 \s} \wt \b_o\left (\frac{J^r_{\ell}}{4\s}+\frac{\d^r_{\ell}}{2\s}\ps_o(\s) \right )+m_o\ps_o(\wt \b_o)\frac{\ps_o(\s)}{2\s}+\frac{m_o}{2}\ps_o(\ps_o(\wt \b_o))=\\
 \label{610**} & =\frac{m_o}{2} \left (\ps_o(\ps_o(\wt \b_o))+\ps_o(\wt \b_o)\frac{\ps_o(\s)}{\s} +\frac{\wt \b_o}{4\s^2}\right ) = 0\ .
  \end{align}
 Recalling that  $\ps_o = \frac{\p}{\p r}$  and that  $\s$  has the form $\s (x, y, r) = - \frac{r^2}{4 B \f(x,y)} - \frac{B \f(x,y)}{4}$,    the equation \eqref{610**} can be written as 
\beq \label{611**} \s \wt \b_{o\,rr}  +  \wt \b_{o\, r} \s_r +  \frac{\wt \b_o}{4 \s} = 0.\eeq
 The general solution of \eqref{611**} is 
\beq\label{formino} \wt \b_o(x, y, r) = \cappa_1(x,y) \frac{ \left( B^2 \f^2(x,y)  - r^2\right)}{r^2 + B^2 \f^2(x,y)}  + \cappa_2(x,y)   \frac{ r}{r^2 + B^2 \f^2(x,y)}\ ,\eeq
where $\cappa_1 = \cappa_1(x,y), \cappa_2 = \cappa_2(x,y)$   are two freely specifiable  functions,  independent of  $r$. Therefore,   if  $g^{(m_o)}$ is Ricci flat,   on each set  $\wc \pi^{-1}(\cV_A)$ such a metric  is completely determined by   the constant $B$,  
  the  potential $\f$  and a free choice of the  functions $\cappa_1, \cappa_2: \cV_A \to \bR$, provided that $\s$ is as in \eqref{form} and $\wt \b_o$ is as in \eqref{formino}. From now on, we assume all this, so that we  automatically have   $ \Ric^{(m_o)}(\ps_o, \qs_o) =  0$.  
Let us now consider the remaining  equations in   \eqref{63bis}.  For any $1 \leq i, j \leq 2$, we have   \par
 \begin{align} 
\nonumber  &  \Ric^{(m_o)}(\wh E_i, \wh E_j)  =   \xcancel{\Ric^{(m= 0)}(\wh E_i, \wh E_j)}   -   \HH m j  {\ps_o} \GGa i {\ps_o}  m   -  \GGa m j  {\qs_o}    \HH i {\qs_o}  m +\\
\nonumber  & +    \HH  i j {\ps_o} \GGa m {\ps_o}  m    +  \xcancel{ \GGa i j {\qs_o}  \HH m {\qs_o}  m }  +   \o_{mi} \HH {\qs_o} j m +\\
\nonumber &+  \ps_o (\HH i j {\ps_o}) -  \GGa {\ps_o} j  \ell   \HH i \ell  {\ps_o}     + \GGa i j {\qs_o}  \HH  {\ps_o} {\qs_o} {\ps_o}  -   \HH {\qs_o} j  \ell \GGa i  \ell {\qs_o}    + \GGa i j {\qs_o}  \HH {\qs_o} {\qs_o} {\qs_o}  = \\
\nonumber  & =   - \frac{ m_o  \wt \b_o \o_{ij}}{4 \s}    \ps_o(\s)  -    \frac{    m_o  \wt \b_o g_{ij}}{2 \s}(\ps_o(\s))^2   - \frac{m_o \wt \b_o g_{ij}}{8\s}   +   \frac{m_o \wt \b_o \o_{ij}}{4 \s}    \ps_o(\s)  + \\
\nonumber  & +  \frac{ m_o  \wt \b_o g_{ij}   }{ \s}(\ps_o(\s))^2    -   \frac{m_o \wt \b_o g_{j i} }{4 \s}  +  m_o g_{ij}   \ps_o (    \wt \b_o )  \ps_o(\s) +  m_o   \wt \b_o  g_{ij}     \ps_o (  \ps_o(\s))  + \\
\nonumber & +     \frac{ m_o    \wt \b_o  \o_{ij }}{4 \s}   \ps_o(\s)    -  \frac{  m_o  \wt \b_o  g_{ij}   }{2 \s}(\ps_o(\s))^2         -  \frac{   m_o  \o_{ij}}{4} \ps_o(\wt \b_o) -   \frac{m_o g_{ij} }{2 } \ps_o(\wt \b_o)  \ps_o(\s)  +  \\
\nonumber & +   \frac{  m_o   \wt \b_o   g_{j i }}{8 \s}  +     \frac{m_o \wt \b_o  \o_{j i } }{4 \s}    \ps_o(\s)    +     \frac{ m_o  \o_{ij}}{4} \ps_o(\wt \b_o)  +     \frac{m_o  g_{ij}  }{2 }\ps_o(\wt \b_o) \ps_o(\s)= \\
\nonumber  & =    m_o g_{ij} \bigg(    -   \frac{ \wt \b_o  }{4 \s}  +   \ps_o (    \wt \b_o )  \ps_o(\s) +     \wt \b_o     \ps_o (  \ps_o(\s))     \bigg)  
 \end{align}
If we now  plug  the  expressions for $\wt \b_o$ and $\sigma$ into this formula, we get that 
\beq   \Ric^{(m_o)}(\wh E_i, \wh E_j)  =\frac{ m_o \cappa_1 g_{ij}}{2 B \f}\ .\eeq
Since $m_o \neq 0$, this means that the  equations   $\Ric^{(m_o)}(\wh E_i, \wh E_j) = 0$, $i, j = 1, 2$,  are satisfied   if and only if  the function $\cappa_1 = \cappa_1(x, y)$ is identically zero   $\wt \b_o = \wt \b_o(x, y, r)$ has the form
\beq\label{formino-bis} \wt \b_o(x, y, r) = \cappa(x, y)   \frac{ r}{r^2 + B^2 \f^2(x, y)}\ ,\eeq
for a  freely specifiable  function  $  \cappa = \cappa(x,y)$. \par
\medskip
   \subsection{Third step}\label{thirdstep3} From now on,  we   assume not only  that  $\s$ is as in \eqref{form} but also   that $\wt \b_o$ is as in \eqref{formino-bis} so that  all equations  in   \eqref{62bis} and \eqref{63bis} are satisfied.
Let us  expand the  component  $\Ric^{(m_o)}(\qs_o, \wh E_i)$, $i = 1,2$,  in terms of the Christoffel symbols $ \wt {\GGa A B C} = \GGa A B C + \HH A B C$. We recall that, according to the hypotheses on the coordinates $(x, y)$ and on the frame field  $(E_1, E_2)$, all   functions $c_{ij}^k$ are zero and the components   $J_i^j$ are   equal to either $\pm 1$ or $0$ at all points. Therefore,  $\Ric^{(m_o)}(\wh E_i, \qs_o)$  is equal to  
   \par
              \scalebox{0.8}{\vbox{
   \begin{align}\nonumber  & \Ric^{(m_o)}(\wh E_i, \qs_o)\   =   \xcancel{ \Ric^{(m = 0)}(\wh E_i, \qs_o) }+\wh E_r (\HH i {\qs_o} r)    -  \HH r {\qs_o}  \ell \GGa i  \ell  r   -  \HH r {\qs_o}  {\ps_o} \GGa  i {\ps_o}  r   +\HH i {\qs_o} \ell  \GGa r \ell  r+ \\
\nonumber  &\ \  + \HH i {\qs_o} {\ps_o}  \GGa r {\ps_o}  r+\o_{ri} \HH {\qs_o} {\qs_o} r  +   \ps_o (\HH i {\qs_o} {\ps_o})-\wh E_i (\HH {\ps_o} {\qs_o} {\ps_o})   -\wh E_i (\HH {\qs_o} {\qs_o} {\qs_o}) -   \HH {\qs_o} {\qs_o}  \ell \GGa i \ell  {\qs_o}   = 
\end{align}
}
}
\ \\
 \phantom{aaaa}   \scalebox{0.8}{\vbox{
   \begin{align}
\nonumber   &\ \ = - m_o J^\ell_i\ \frac{\s_\ell}{4 \s^2}\wt \b_o  +   m_o J^\ell_i \frac{1}{4 \s}\wt \b_{o\ell }  + m_o \f_i B \frac{\ps_o(\s)}{4 \s^2} \wt \b_o - m_o \f_i \frac{B \ps_o (\wt \b_{o})}{4 \s} - \\
\nonumber & \ \ - m_o \frac{1}{4 \s}\wt \b_o  \Bigg ( J^{\ell}_r\Ga i \ell r
 + \frac{ \s_\ell  }{ \s }  J^{\ell}_i
  - B  \frac{\ps_o(\s) \f_i }{ \s }\Bigg )-  \left ( \frac{m_o}{2} \wt \b_{or} - \frac{m_o}{2}J^\ell_r \f_\ell  \ps_o(\wt \b_o)\right ) \left ( \frac{J^r_i}{4 \s} +   \frac{\d^r_i}{2 \s}\ps_o(\s)\right) +\\
\nonumber  &\ \ + m_o \frac{\wt \b_o J^{\ell}_i}{4\s} \left(\Ga r \ell  r+\frac{\s_\ell}{\s} +B \frac{ \f_t J^t_{\ell} \ps_o(\s)}{\s} 	\right )
 + \left ( \frac{m_o}{2} \wt \b_{oi} - \frac{m_o}{2}J^\ell_i \f_\ell \ps_o(\wt \b_o)\right )   \frac{\ps_o(\s)}{\s} +\\
\nonumber  &\ \  + m_o \frac{J^k_i }{4\s} \wt \b_{ok}  -  m_oB \frac{1 }{4\s} \f_i  \ps_o(\wt \b_o)+  \frac{m_o}{2} \ps_o \left (\wt \b_{oi}  -  J^\ell_i\ps_o(\wt \b_o)\f_\ell \right )- \xcancel{\frac{m_o}{2}	\wh E_i(\ps_o(\wt \b))}+\\
\nonumber  &\ \ + \xcancel{\frac{m_o}{2} \wh E_i(\ps_o(\wt \b_o))} -\left(  m_o\frac{g^{\ell k}}{4\s} \wt \b_{ok} + m_o B\frac{g^{\ell k}}{4\s}J_k^s \f_s \ps_o(\wt \b_o) \right ) \left (  \frac{   \o_{i\ell}}{2}  +  g_{i\ell}   \ps_o(\s)\right ) = \\
\nonumber   &\ \ = - m_o J^\ell_i\ \frac{\s_\ell}{4 \s^2}\wt \b_o  +   m_o J^\ell_i \frac{1}{4 \s}\wt \b_{o\ell } + m_o \f_i B \frac{\ps_o(\s)}{4 \s^2} \wt \b_o - m_oB \f_i \frac{ \ps_o (\wt \b_{o})}{4 \s} - \\
\nonumber &    \boxed{- m_o \frac{1}{4 \s}\wt \b_o  J^{\ell}_r\Ga i \ell r}
 - \xcancel{ m_o \frac{\s_\ell  }{4 \s^2}\wt \b_o  J^{\ell}_i}
  +  \xcancel{m_o \frac{1}{4 \s^2}\wt \b_o B  \ps_o(\s) \f_i } - \\
\nonumber & \ \  - \xcancel{m_o \wt \b_{o\ell } \frac{J^\ell_i}{8 \s} } - m_o \f_i  \ps_o(\wt \b_o)  \frac{1}{8 \s}   -  \xcancel{m_o \wt \b_{oi}  \frac{1}{4 \s}\ps_o(\s) } + m_o J^\ell_i \f_\ell  \ps_o(\wt \b_o) \frac{1}{4 \s}\ps_o(\s)       +
\\
\nonumber   &\ \ = \boxed{m_o \frac{1}{4 \s}\wt \b_o  \underset{ = 0} {\underbrace{\left( - J^{\ell}_r\Ga i \ell r +  J^{\ell}_i  \Ga r \ell  r\right)} } } - m_o J^\ell_i\ \frac{\s_\ell}{4 \s^2}\wt \b_o  +   m_o J^\ell_i \frac{1}{4 \s}\wt \b_{o\ell } +   m_o \f_i B \frac{\ps_o(\s)}{4 \s^2} \wt \b_o   - \\
\nonumber & \ \  - m_o \f_i  \ps_o(\wt \b_o)  \frac{(3B+1)}{8 \s}  - m_o(B+1) J^\ell_i \f_\ell  \ps_o(\wt \b_o) \frac{1}{4 \s}\ps_o(\s)     +  \frac{m_o}{2} \ps_o  (\wt \b_{oi})   -  J^\ell_i  \frac{m_o}{2}   \hskip - 10 pt \underset{ = (-  \ps_o(\wt \b_o)\frac{\ps_o(\s)}{\s} - \frac{\wt \b_o}{4\s^2})}{\underbrace{\ps_o (\ps_o(\wt \b_o))}} \hskip - 8 pt \f_\ell     = \\
\nonumber   &\ \ = \frac{m_o}{4 \s^2} \bigg( -  J^\ell_i\ \s_\ell \wt \b_o  +    J^\ell_i \s \wt \b_{o\ell } +  B \f_i  \ps_o(\s)  \wt \b_o - \\
\label{andrea's}  & \ \  -   \frac{(3B+1)}{2}  \f_i  \s \ps_o(\wt \b_o)  -   (B-1) J^\ell_i \f_\ell \s \ps_o(\s)  \ps_o(\wt \b_o)      +  2 \s^2  \ps_o  (\wt \b_{oi})   + J^\ell_i \f_\ell   \frac{\wt \b_o}{2 }    \bigg)\ .
\end{align}
}}\\
We now plug   the explicit formulas  for $\sigma$ and $\wt \b_o$ into the above expressions.  Such a replacement is straightforward, but very  arduous and it is convenient  to use  a  symbolic manipulation program,   such as  {\it Maple}, to perform it.   In this way  the equations $ \Ric^{(m_o)}(\wh E_i, \qs_o) = 0$,  $i = 1, 2$, become equivalent to the following equations on  $ \cappa = \cappa(x, y)$: 
\begin{multline}\label{AA}
\Bigg(B \cappa \f (B +1) (B \f -r )  (B \f  +r ) \frac{\partial \f }{\partial x} 
+ 5 r (B +1) \cappa  \left(B^2 \f ^2-\frac{r^2}{5}\right)  \frac{\partial \f }{\partial y}+\\
+\bigg((B^2 \f ^2-r^2) \frac{\partial \cappa }{\partial x}-2 B r \f   \frac{\partial \cappa }{\partial y}\bigg) (B^2  \f ^2+r^2) \Bigg) \frac{1}{2 (B^2 \f^2+r^{2})^3} = 0\ ,
\end{multline}
\begin{multline}\label{BB}
\Bigg(B \cappa \f  (B +1) (B \f -r ) (B \f +r ) \frac{\partial \f }{\partial y} 
-5 (B +1) \cappa  r \left(B^2  \f^2-\frac{r^2}{5}\right) \frac{\partial \f }{\partial x}
\\
+\left( \frac{\partial \cappa }{\partial y} (B^2\f^2-r^2) + 2 B r  \f  \frac{\partial \cappa}{\partial x}\right) (B^2 \f^2+r^2)\Bigg) \frac{1}{2 (B^2 \f^2 +r^{2})^{3}} = 0\ .
\end{multline}
 The numerators of the left hand sides of these  equations are polynomials of order $4$ in the variable $r$. The   monomials in $r$ of highest degree in these  equations have coefficients $- \frac{\p \cappa}{\p x} $ and $\frac{\p \cappa}{\p y}$, respectively.  Since  the    equations are  satisfied  for any value of $r$ only if such  coefficients   vanish identically,  it follows that they are  satisfied  {\it  only if   $\cappa(x,y)$ is  constant}.    If we now assume  that $\cappa$ is a non-zero constant ($\cappa \neq 0$  because otherwise $\wt \b_o$ would be  $0$ and  $g^{(m_o)} = \h$), the equations \eqref{AA} and \eqref{BB}  reduce to
\beq \label{AA*}
\begin{split}
& (B+1) \cappa \bigg( B  \f (B^2 \f^2 -r^2 )  \frac{\partial \f }{\partial x}
+  r   \left(5 B^2 \f ^2- r^2 \right)  \frac{\partial \f }{\partial y} \bigg)= 0\ ,\\
&(B +1)  \cappa \bigg(B  \f  (B^2 \f^2 -r^2) \frac{\partial \f }{\partial y}
-  r \left(5 B^2  \f^2- r^2 \right) \frac{\partial \f }{\partial x} \bigg)= 0\ .
\end{split}
\eeq
The left hand sides of these equations are now polynomials in $r$ of  order $3$ and  they vanish identically if and only if either $\f = \f(x,y)$ is constant (but this cannot be because $\f$ is a potential of a K\"ahler metric) or  $B = -1$. We  conclude  that all equations in \eqref{62bis}, \eqref{63bis} and \eqref{64bis} {\it with the  only exception of  $\Ric^{(m_o)}(\qs_o, \qs_o) = 0$},   are  satisfied  by  a  $g^{(m_o)} \neq \h$   if and only if   $ \wt \b_o (x ,y)$ is as in \eqref{lella} and $B = -1$.\par
\smallskip
It remains to check that  if  $\wt \beta_o$ and $B$ are  as  above,   the equation $\Ric^{(m_o)} (\qs_o, \qs_o) = 0$ is identically satisfied. As usual, let us   expand  $\Ric^{(m_o)}(\qs_o, \qs_o)$ in terms of the   $ \wt {\GGa A B C} = \GGa A B C + \HH A B C$. We get 
      \par
    \scalebox{0.8}{\vbox{
  \begin{align}\nonumber  & \Ric^{(m_o)}({\qs_o},{\qs_o})\  {=}  \xcancel{\Ric^{(m = 0)}({\qs_o},{\qs_o})} +  \wh E_r (\HH {\qs_o} {\qs_o} r) -  \xcancel{\qs_o (\HH r {\qs_o} r) }- 2 \HH r {\qs_o}  \ell \GGa  {\qs_o} \ell r    -  \HH r {\qs_o}  \ell \HH  {\qs_o} \ell r     +
\\
\nonumber  &\hskip 1.5 cm  +   \HH {\qs_o} {\qs_o} \ell  \GGa  r \ell  r  +   \xcancel{ \GGa {\qs_o} {\qs_o} \ell \HH  r \ell  r }   + \xcancel{ \HH {\qs_o} {\qs_o} \ell  \HH  r \ell  r }  +  \HH {\qs_o} {\qs_o} {\ps_o} \GGa r  {\ps_o} r +  \xcancel{  \GGa {\qs_o} {\qs_o} {\ps_o} \HH r  {\ps_o} r }+ \xcancel{\HH {\qs_o} {\qs_o} {\ps_o} \HH r  {\ps_o} r }+ 
\HH {\qs_o} {\qs_o} {\qs_o}  \GGa r {\qs_o}  r  +\\
\nonumber &\hskip 1.5 cm+   \ps_o (\HH {\qs_o} {\qs_o} {\ps_o}) -   \qs_o (\HH {\ps_o} {\qs_o} {\ps_o})    - \xcancel{ 2 \GGa {\ps_o} {\qs_o}  {\ps_o} \HH {\ps_o} {\qs_o}  {\ps_o} }
   - (\HH {\ps_o} {\qs_o}  {\ps_o}  )^2      + \HH {\qs_o} {\qs_o} {\qs_o} \HH {\ps_o} {\qs_o}  {\ps_o} =\\
 \nonumber  & =  
 -m_o \wh E_r \left( \frac{g^{rk}}{4\s} \wt \b_{ok} +B  \frac{g^{rk}}{4\s}J_k^\ell  \f_\ell \ps_o(\wt \b_o) \right ) 
  -   2 m_o  \frac{J^r_{\ell}}{4\s}\wt \b_o	  \left (B\frac{J^{\ell}_r}{4 \s} 
  + B \frac{\d^{\ell}_r}{2\s}\ps_o(\s)  \right ) + \frac{m_o^2}{8\s^2}\wt \b^2_o     - \\
\nonumber & 
 - m_o \left ( \frac{g^{\ell k}}{4\s} \wt \b_{ok} +B  \frac{g^{\ell k}}{4\s}J_k^t  \f_t \ps_o(\wt \b_o) 
 \right )  \GGa  r \ell  r +   m_o \left (B	 \ps_o(\wt \b_o) + 
  m_o \frac{\ps_o(\wt \b_o) \wt \b_o }{2}\right )\frac{\ps_o(\s)}{\s} 
  -\frac{m_o B}{2\s}\ps_o(\s)\ps_o(\wt \b_o)\\
\nonumber  & +m_o \ps_o \left ( 	B	 \ps_o(\wt \b_o)+m_o\frac{\ps_o(\wt \b_o) \wt \b_o }{2}\right ) 
-   \frac{B m_o}{2 }\ps_o \left (\ps_o(\wt \b_o) \right )    - m_o^2 \frac{\ps_o(\wt \b_o)^2}{4}      - m_o^2  \frac{\ps_o(\wt \b_o)^2}{4}   = \\
 \nonumber  & =  
 -m_o \wh E_r \left( \frac{g^{rk}}{4\s} \wt \b_{ok} +B  \frac{g^{rk}}{4\s}J_k^\ell  \f_\ell \ps_o(\wt \b_o) \right )  -m_o \left ( \frac{g^{\ell k}}{4\s} \wt \b_{ok} +B  \frac{g^{\ell k}}{4\s}J_k^t  \f_t \ps_o(\wt \b_o) 
 \right )  \GGa  r \ell  r +\\
 \nonumber & 
 +  m_o B\frac{ \wt \b_o }{4\s^2} +  m_o^2 \frac{\wt \b_o^2}{8\s^2} +  m_o B	\frac{  \ps_o(\wt \b_o) \ps_o(\s)}{\s} +   m_o^2\frac{ \ps_o(\s) \ps_o(\wt \b_o)   \wt \b_o }{2\s }
- m_o B\frac{ \ps_o(\s)\ps_o(\wt \b_o)  }{2\s}+\\
\label{klivia} & + 	m_o B \ps_o(\ps_o(\wt \b_o))
+m_o^2\frac{\ps_o(\ps_o(\wt \b_o)) \wt \b_o }{2}
+\xcancel{ m_o^2\frac{\ps_o(\wt \b_o)^2 }{2} } 
-   m_o B \frac{ \ps_o  (\ps_o(\wt \b_o)  )  }{2 }  - \xcancel{m_o^2 \frac{\ps_o(\wt \b_o)^2}{2}  }     
 \end{align}
}}\ \\
Since   \eqref{610**} implies  $\ps_o(\ps_o(\wt \b_o))+\ps_o(\wt \b_o)\frac{\ps_o(\s)}{\s} +\frac{\wt \b_o}{4\s^2} = 0$, the sum of the three terms in \eqref{klivia} that  are  proportional to $m_o^2$ vanishes identically, while  the sum of all other terms in the last two lines of  that expression reduces to $ m_o B\frac{ \wt \b_o }{8\s^2}$. We therefore get \par
 \scalebox{0.8}{\vbox{
  \begin{align}\nonumber  & \Ric^{(m_o)}({\qs_o},{\qs_o})\  = \\
  \nonumber & = 
 -m_o \wh E_r \left( \frac{g^{rk}}{4\s} \wt \b_{ok} +B  \frac{g^{rk}}{4\s}J_k^\ell  \f_\ell \ps_o(\wt \b_o) \right )  -m_o \left ( \frac{g^{\ell k}}{4\s} \wt \b_{ok} +B  \frac{g^{\ell k}}{4\s}J_k^t  \f_t \ps_o(\wt \b_o) 
 \right )  \GGa  r \ell  r + m_o B\frac{ \wt \b_o }{8\s^2} =\\
   \nonumber  & =  
 -m_o \frac{\wt \b_{ok}}{4\s}  E_r(g^{rk}) 
  +\xcancel{m_o g^{rk} \frac{\s_r}{4\s^2} \wt \b_{ok}}+\xcancel{m_o Bg^{rk}J^s_r\f_s \frac{\ps_o(\s)}{4\s^2} \wt \b_{ok}}  -m_o \frac{g^{rk}}{4\s}\wt \b_{okr}-m_o B \frac{g^{rk}}{4\s}J^s_r\f_s \ps_o(\wt \b_{ok})-
\\
\nonumber  &
-m_o B   E_r(g^{rk}) \frac{J_k^\ell}{4\s}  \f_\ell \ps_o(\wt \b_o)  +\xcancel{m_o B  \frac{g^{rk}}{4\s^2} \s_r J_k^\ell  \f_\ell \ps_o(\wt \b_o)  }
+\xcancel{m_o B^2  \frac{g^{rk}}{4\s^2}J^s_r\f_s \ps_o (\s)J_k^\ell  \f_\ell \ps_o(\wt \b_o) } - \\
\nonumber  &-m_o B  \frac{g^{rk}}{4\s}J_k^\ell  \f_{\ell r} \ps_o(\wt \b_o) 
-m_o B  \frac{g^{rk}}{4\s}J_k^\ell  \f_\ell  E_r\left (\ps_o(\wt \b_o) \right )
-m_o B^2  \frac{g^{rk}}{4\s}J_k^\ell  \f_\ell J^s_r\f_s\ps_o\left (\ps_o(\wt \b_o) \right )
- \\
\nonumber  & -m_o \frac{g^{\ell k}}{4\s} \wt \b_{ok} \Ga  r \ell  r
 -\xcancel{m_o \frac{g^{\ell k}}{4\s^2} \wt \b_{ok}\s_\ell}
 -\xcancel{m_o B \f_tJ^t_\ell\frac{g^{\ell k}}{4\s^2} \wt \b_{ok}\ps_o(\s)}-	  \\
\nonumber  &   -m_o B \frac{g^{\ell k}}{4\s}J_k^t  \f_t \ps_o(\wt \b_o)   \Ga  r \ell  r 
-\xcancel{m_o B  \frac{g^{\ell k}}{4\s^2}J_k^t  \f_t \ps_o(\wt \b_o)\s_\ell }
-\xcancel{m_o B^2 \f_tJ^t_\ell  \frac{g^{\ell k}}{4\s^2}J_k^s  \f_s \ps_o(\wt \b_o) \ps_o(\s)} +\frac{m_o B}{8 \s^2} \wt \b_o  = 
\end{align}
}
}
\par
 \scalebox{0.8}{\vbox{
\begin{align}
 \nonumber  & =  
  -m_o \frac{g^{rk}}{4\s}\wt \b_{okr}-m_o B \frac{g^{rk}}{2\s}J^s_r\f_s \ps_o(\wt \b_{ok}) 
 -m_o B  \frac{g^{rk}}{4\s}J_k^\ell  \f_{\ell r} \ps_o(\wt \b_o) 
-m_o B^2  \frac{g^{rk}}{4\s}J_k^\ell  J^s_r   \f_\ell\f_s\ps_o\left (\ps_o(\wt \b_o) \right ) +\frac{m_o B}{8 \s^2} \wt \b_o 
- \\
\label{Klivia1}  & -m_o \frac{\wt \b_{ok}}{4\s} \left( E_r(g^{rk}) +g^{\ell k} \Ga  r \ell  r\right) - m_o B  \frac{J_k^\ell}{4\s}  \f_\ell \ps_o(\wt \b_o) \left( E_r(g^{rk})  +g^{\ell k}  \Ga  r \ell  r \right)
\ .
 \end{align}
 }}\ \\
 We now  recall that, for each  $k ,r , \ell= 1,2$,  the expression $ E_r(g^{rk})  +g^{\ell k}  \Ga  r \ell r$ is identically $0$. One way (out of many others) to check this is   by using the isothermal  coordinates $(x, y)$ on  $\cV $, with   $E_1= \frac{\p}{\p x}$, $E_2 = \frac{\p}{\p y}$,  and   $g^{rs} = \frac{1}{\lambda} \d^{rs}$,   $\lambda = \D \f$. The Christoffel symbols $\Ga r \ell s$  of the metric $g_o$ in these coordinates are given by \eqref{6.78}. Using   this information,  a  straightforward check shows that $ E_r(g^{rk})  +g^{\ell k}  \Ga  r \ell r = 0$,  which we plug into 
   into \eqref{Klivia1} to get  \par
  \scalebox{0.8}{\vbox{
  \begin{align}\nonumber  & \Ric^{(m_o)}({\qs_o},{\qs_o})\  = \\
  \nonumber & \ \ = -m_o \frac{g^{rk}}{4\s}\wt \b_{okr}-m_o B \frac{g^{rk}}{2\s}J^s_r\f_s \ps_o(\wt \b_{ok}) 
 -m_o B  \frac{g^{rk}}{4\s}J_k^\ell  \f_{\ell r} \ps_o(\wt \b_o) 
-m_o B^2  \frac{g^{rk}}{4\s}J_k^\ell  J^s_r   \f_\ell\f_s\ps_o\left (\ps_o(\wt \b_o) \right ) +\frac{m_o B}{8 \s^2} \wt \b_o  = 
\end{align}
}
}
\ \\
\scalebox{0.8}{\vbox{
  \begin{align}
  \nonumber  &\ \  =  
    - \frac{m_o}{4\lambda \s}\wt \b_{oxx}- \frac{m_o}{4\lambda \s}\wt \b_{oyy}
\xcancel{  - \frac{m_o B}{2\lambda\s}\f_y \ps_o(\wt \b_{ox})}
\xcancel{  +\frac{m_o B}{2\lambda\s}\f_x \ps_o(\wt \b_{oy})}
 - \xcancel{ \frac{m_o B}{4\s} \f_{y x} \ps_o(\wt \b_o)}
+ \xcancel{\frac{m_o B }{4\lambda\s} \f_{x y} \ps_o(\wt \b_o) } - \\
\nonumber  &\ \ \ \ \ -  \frac{m_o B^2}{4\lambda\s}  \f_y^2 \ps_o\left (\ps_o(\wt \b_o) \right )
-  \frac{m_o B^2}{4\lambda\s}  \f_x^2 \ps_o\left (\ps_o(\wt \b_o) \right )
 +\frac{m_o B}{8 \s^2} \wt \b_o  = \\
  \label{Klivia2}  &\ \  =  
    - \frac{m_o }{4\lambda \s}\left ( \wt \b_{oxx}+\wt \b_{oyy}\right )
-  \frac{m_o B^2}{4\lambda\s}   \ps_o\left (\ps_o(\wt \b_o) \right ) \left ( \f_y^2+\f_x^2\right )
 +\frac{m_o B}{8 \s^2} \wt \b_o \ .
 \end{align}
 }
 }
 \ \\
 We may now use the explicit expressions \eqref{form} and \eqref{formino-bis} for  $\sigma$ and $\wt \b_o$ in   terms of  $\f(x,y)$ and $\cappa$.  With the help of  a symbolic manipulation program such as {\it Maple}, one can immediately check that  these substitutions  yield
 the identity
\beq  \Ric^{(m_o)}({\qs_o},{\qs_o})\  = m_o \frac{2 r \f^2(B+1)}{(r^2+\f^2)^3} \ .
\eeq
Since $B = -1$,  it follows that $ \Ric^{(m_o)}({\qs_o},{\qs_o}) \equiv 0$, as we needed to prove.  The computations  of the three steps  show also that there is no constraint on  $m_o$ and that  if   a metric $g^{(m_o)}$ is Ricci flat, then  automatically  any other metric $g^{(m)}$, $m \neq m_o$, is Ricci flat.
\par
\medskip
\section{Explicit expressions of   Ricci flat  Kerr families  and    real analytic extensions of some of such metrics} \label{expli-sol}
\subsection{Potentials  for Kerr type   backgrounds and  corresponding Kerr families} 
By Corollary \ref{333}, if $M$ is   a Kerr $4$-manifold admitting a Kerr type background,  then the underlying  Riemann surface $(N, J)$ is  non-compact, admitting a K\"ahler metric $\wt g_o$  of constant Gaussian curvature $\k = \pm 1 $ and   equipped with a K\"ahler metric $g_o$  with a  global potential $\f: N \to \bR$  satisfying \eqref{42}.   In  suitable  isothermal coordinates $(x, y)$ for $\wt g_o$  and  $g_o$ in an open subset $\cV \subset N$,  the conditions \eqref{42}  are equivalent to either 
\begin{align} \label{firsteq} &  \Delta \f  + \frac{8 \f}{\left(1 +(	x^2 + y^2)\right )^2} =0\quad \&\quad \f < 0\ 
\qquad \text{or to}\\
 \label{secondeq}  &  \Delta \f  - \frac{8\f }{\left(1 - ( x^2 + y^2)\right )^2} =0 \quad \&\quad \f > 0\ . \end{align}
The  solutions     to these equations  on  bounded open sets  have the following properties.  \par 
 \begin{lem} \label{Lemma61}
  Let $\cV$ be  a relatively compact region in  $\mathbb D= \{|z|^2 < 1\} \subset \bC = \bR^2$ with boundary of class $\cC^{k, \a}$, $k \geq 2$, $\a \in (0,1)$.  Any  solution to  \eqref{firsteq} (resp. \eqref{secondeq}) on $\overline\cV$ is  real analytic on $\cV$ and  the correspondence 
  $$\f \longmapsto f = \f|_{\p \cV}$$
   is a bijection  between the solutions $\f$ that are in  $\cC^\o(\cV) \cap \cC^{k, \a}(\overline \cV)$
  and  the   non identically vanishing  functions  $f \leq 0$    (resp.  $f \geq 0$)  that are in  $\cC^{k, \a}(\p \cV)$.
     \end{lem}
     \begin{pf} By \cite[Thms. 3.5, 6.13 \&  6.19]{GT} and by  \cite{Ho} (see also \cite{Jo}), there exists a unique solution to \eqref{secondeq} on $\overline \cV$ of class $\cC^\o(\cV) \cap \cC^{k, \a}(\overline \cV)$ for any prescribed boundary value $f = \f|_{\cV}$ in   $\cC^{k, \a}(\partial \cV)$ and such a solution satisfies the Maximum Principle. This implies the claim for the conditions  \eqref{secondeq}. 
     This argument does not work for \eqref{firsteq} because, in this case,  the zero order term  is positive and \cite[Thms. 3.5 \&  6.13]{GT} do not apply. Nonetheless, we may consider the special  negative solution to \eqref{firsteq} 
     \beq \f_o: \bD \to  (0, + \infty) \ ,\qquad \f_o(x,y) \= -\frac{1 - x^2 - y^2}{1 + x^2 + y^2}\eeq
  and observe that the correspondences $\f \mapsto \wt \f \= \frac{\f}{\f_o}$ and $\wt \f \mapsto \f = \wt \f \f_o$ establish a bijection between the solutions to  \eqref{firsteq}  and the solutions to 
      \beq \label{firsteq*} \D \wt \f + \frac{2}{\f_o}\left( \frac{\p \f_o}{\p x}    \frac{\p \wt \f}{\p x} + \frac{\p \f_o}{\p y}    \frac{\p \wt \f}{\p y} \right)= 0\eeq
      on $\overline \cV \subset \bD$. Since the previous argument now works for \eqref{firsteq*}, the lemma is proved for \eqref{firsteq} as well. 
       \end{pf}
      Lemma \ref{Lemma61},  Corollary \ref{333} and Theorem \ref{aaa} immediately yield  the following 
       \begin{cor} Any background metric $\h$ of Kerr type  and any  metric $g^{(m)} $, $m \neq 0$, of a Kerr family of gravitational fields on a  Kerr $4$-manifold $M$ are real analytic in  the coordinates $(x, y, v, r)$ considered in Corollary \ref{333} and  Theorem \ref{aaa}. 
       In particular, each  such metric   is uniquely determined by its restriction to any  open subset $\cV \subset N$ of  the open cover $\{\cV_A\}$ of $N$,  described in  Corollary \ref{333}.
       \end{cor}
      We note that, as  a direct consequence of  \cite{MN},  the potentials that satisfy  \eqref{firsteq}   (resp. \eqref{secondeq}) on  a  disk $\overline{\bD(0,\gr_o)}\= \{|z|\leq \gr_o\} $, $0 < \gr_o < 1$,   are real analytic up to the boundary, thus  in one-to-one correspondence with  the Fourier series  of real analytic  boundary data on $\p \bD(0,\gr_o)$. 
     In fact, using  polar coordinates $(\gr , w)$  with  $(x = \gr \cos w, y =  \gr \sin w)$ and  considering the Fourier expansion $f(w) = a_0 + \sum^\infty_{n = 1}( a_n  \cos(n w) +  b_n \sin(n w))$ of a non-negative real analytic  function  $f(w)$ on $\partial \bD(0,  \gr_o) $, $\gr_o < 1$, one can check that  the corresponding solution $\f$  to \eqref{firsteq}  (resp. \eqref{secondeq})  is  the sum of the  series 
     \beq\label{explicita}
     \begin{split}
     & \f(\gr, w) = -  \left( a_0  \frac{\f_0(\gr)}{\f_0(\gr_o)} + \sum^\infty_{n = 1}\frac{\f_n(\gr)}{\f_n(\gr_o)}\big(  a_n \cos(n w) +   b_n \sin(n w)\big)\right)\ ,\\
     & \text{where}\ \ \ \f_n(\gr) \=   \left(1 - \frac{2 \gr^2}{ (1 + \gr^2)}\frac{1}{n+1} \right)  \gr^n \\
     \bigg(& \text{resp. }\ \ \  \f_n(\gr) \=\left(\frac{1+\gr^2}{1 - \gr^2} + n\right) \gr^n\ \ \bigg)\ .
     \end{split}
     \eeq
   Each of such solutions  can be conveniently  expressed in terms of a single holomorphic function as follows.  Using the complex coordinate   $z=\gr e^{iw}$ and setting  $c_0 \=\frac{- a_0}{\f_0(\gr_o)}$ and $c_n\=\frac{-  (a_n-ib_n)}{\f_0(\gr_o)}$, $n \geq 1$, we have that \eqref{explicita} can be re-written as 
  \begin{align*}
  \f(z)&=\Re \left(\sum_{n=0}^\infty   \left(1 - \frac{2 |z|^2}{ (1 + |z|^2)}\frac{1}{n+1} \right) c_n z^n\right)\ \ &&  \text{ if } \k=1\ , \\
  \f(z)&=\Re \left(\sum_{n=0}^\infty  \left(\frac{1+|z|^2}{1 - |z|^2} + n\right) c_n z^n\right) \ \ && \text{ if } \k=-1\ . 
  \end{align*}
Considering the holomorphic function    $F(z) \=\sum_{n=0}^\infty c_n z^n$,  we get that    
\beq \label{68}
\begin{split}
 \f(z)&= \Re \left(F(z)- \frac{2|z|^2}{1+|z|^2}\frac1{z} \int_0^z F(\z) d\z\right)\qquad  \text{ if } \k=1\ ,\\
  \f(z)&=\Re \left( \frac{1+|z|^2}{1-|z|^2} F(z) + z\frac{dF(z)}{dz} \right) \qquad\qquad  \text{ if } \k=-1\ .\\
\end{split}
\eeq
Conversely,  for any  holomorphic function $F$,
   for which the corresponding function  in  \eqref{68} is non-positive or  non-negative (according to $\k = 1$ or $-1$),  
determines a solution  $\f(z)$ to \eqref{firsteq}  or \eqref{secondeq}, respectively. 
\par
\medskip
\subsection{Explicit  examples of Ricci flat Kerr families of gravitational fields}
Let $M = (N \times \bR) \times \bR$  (or $M =(N \times \bR) \times I$,   $I \subset \bR$) be a Kerr $4$-manifold admitting a Kerr family of Ricci flat gravitational fields $g^{(m)}$.  As usual, we denote by  $(N, J, g_o)$   the (non-compact)  Riemann surface with K\"ahler metric $g_o$ which underlies the Kerr manifold.  \par
Note that, denoting by  $\wt g_o$  the K\"ahler metric on $(N, J)$ with  constant Gaussian curvature $\k = \pm 1$ that corresponds  to the family $g^{(m)}$ by  Theorem \ref{aaa},  using the relation   $g_o=-2 \k \f \wt g_o$, the  metrics of the given family $g^{(m)}$  have  locally the form 
\beq \label{6666} g^{(m)} =  -  \frac{ \k}{2} \left( r^2 +   \f^2 \right)  \wt g_o  
    +  (d v + d^c \f )\vee \bigg( d(v + r) + \frac{1}{2}  \left( - 1 +   \frac{\cappa m r}{r^2 + \f^2} \right)( d v + d^c \f )\bigg)\ ,
\eeq
    where: 
    \begin{itemize}[leftmargin = 15 pt]
    \item $\cappa$  is  a  scaling constant  (for  the   Kerr metrics \eqref{standardKerr},  such a constant is  $\cappa = 2$); 
    \item $\f$ is the global K\"ahler potential for $g_o$, defined in  claim (2) of Theorem \ref{the51}; 
    \item $r$ is a  coordinate for the standard fiber of  the trivial bundle $\pi^\cS: M = \cS \times \bR \to N \times \bR$; 
    \item  $v: \cU \subset M \to \bR$ is a locally defined function,   which completes  a triple $(x, y, r)$, given by a pair of local coordinates  $(x, y)$ for $\cV \subset N$ and the globally defined  function $r$,   to a coordinate chart $(x, y, v, r): \cU \subset M  \to \bR^4$, in which  $\qs_o = \frac{\p}{\p v} + B \frac{\p}{\p r}$.  
    \end{itemize}
 Note  also that   the real  function  $v: \cU \to \bR$ is  a local  potential for the globally defined closed $1$-form  $\q - d^c \f$ (i.e., $d v = \q - d^c \f$). Thus,  in case $(N, J)$ (and consequently  also $M = N \times \bR \times \bR$) is simply connected,  such local potential  $v$ can be  extended to a globally defined  potential  $v: N \times \bR \times \bR \to \bR$.  On the other hand,  if  $(N, J)$ is {\it not} simply connected,  several non equivalent   connection $1$-forms $\theta$ satisfying    $d \theta = \pi^\cS \o_o$ might occur, each of them  corresponding to  non equivalent  local potentials $v$. For an exhaustive discussion of such non-equivalent  possibilities, see e.g.  \cite{Ko}.
   \par
 \smallskip
 Plugging   explicit coordinate expressions for  solutions $\f$ to \eqref{firsteq} and \eqref{secondeq} into  \eqref{6666},  one   gets a huge number of explicit coordinate expressions for Lorentzian Ricci flat metrics $g^{(m)}$. In the next two subsections,  we briefly discuss    a few  features of these metrics.  \par
 \smallskip  
\subsubsection{A quick look at   the metrics of the  Kerr families  with  $\k = 1$}
Consider the   expressions  \eqref{explicita} for  possible potentials $\f$ in case $\k = 1$. There, they are   given   in terms of the  polar coordinates  $(\gr, w)$ associated with a pair of  isothermal coordinates $(x = \gr \cos w,y = \gr \sin w)$ in   which $\wt g_o$ takes the form $\wt g_o = \frac{4}{(1 +(x^2 + y^2))^2} (d x^2 + d y^2)$.   Replacing  $(\gr, w)$ by     coordinates $(\xi,  \psi)$  defined     by the relations
 $\gr = \tan \left( \frac{\xi}{2}\right)$,  $w = \psi$, 
  the metric   $\wt g_o$  and the potentials $\f$  take the  form 
    \beq\label{explicita**}
     \begin{split}
     & \wt g_o  = d \xi^2 + \sin^2 \xi d \psi^2\ ,\\
     & \f(\xi, \psi) = \wt a_0 \cos \xi   +   \wt a_1   \sin \x \cos \psi  +   \wt b_1   \sin \x \sin \psi  +   \\
     & \hskip 1 cm + \sum^\infty_{n = 2}  \frac{1 + \frac{n - 1}{n+1}  \tan^{2} \left( \frac{\xi}{2}\right)}{1 +  \tan^2 \left( \frac{\xi}{2}\right)} \tan^n \left( \frac{\xi}{2}\right) \big(  \wt a_n \cos(n \psi) +   \wt b_n \sin(n \psi)\big)
     \end{split}
     \eeq
     for  appropriate constants $\wt a_n, \wt b_n \geq 0$.  In this way,  $\wt g_o$ takes the usual  form  of  the standard  metric of $S^2$  in  spherical coordinates and $\f$ can be (locally) identified with 
     a function on an open subset of the unit sphere.  \par
\smallskip
Let us  now introduce the following  notation.  Given a pair of  real numbers $(\a_1, \a_2)$,  let $R^{(\a_1, \a_2)}$  be the isometry of  $(S^2, \wt g_o)$ given by the double  rotation 
 \begin{multline*} R^{(\a_1, \a_2)}: S^2 \longrightarrow  S^2\ ,\\
  R^{(\a_1, \a_2)}\left( \begin{array}{c} y^1\\y^2\\y^3\end{array} \right) \= \left( \begin{array} {ccc} \cos \a_2 & -\sin \a_2 & 0\\
 \sin \a_2 & \cos \a_2 & 0\\
 0 & 0 & 1 \end{array} \right){\cdot} \left( \begin{array}{ccc}  \cos \a_1 & 0 & -\sin \a_1\\
  0 & 1 & 0 \\
 \sin \a_1 &0 &  \cos \a_1
 \end{array}  \right){\cdot} \left( \begin{array}{c} y^1\\y^2\\y^3 \end{array} \right)\ .\end{multline*}
We may thus consider the  corresponding  ``rotated'' spherical coordinates  $\big(\xi' \= \xi(\a_1, \a_2)$, $\psi' \= \psi(\a_1, \a_2)\big)$ of $S^2$, i.e. the coordinates  that  are related with the usual spherical  coordinates  $(\xi, \psi)$ by  
\beq  \left( \begin{array}{c} \sin \xi'  \cos \psi'  \\ \sin \xi'  \sin \psi' \\ \cos \xi' \end{array} \right)=  R^{(\a_1, \a_2)}\left( \begin{array}{c} \sin \xi \cos \psi \\ \sin \xi \sin \psi\\ \cos \xi\end{array} \right)\ . \eeq
In different words,  $(\xi', \psi')$ are just new spherical coordinates, which one gets  after   an appropriate  replacement of   the north-south axis  and  of the reference meridian $\psi = 0$. \par
\smallskip
A straightforward check shows that  for any $\f(\xi, \psi) = \wt a_0 \cos \xi   +   \wt a_1   \sin \x \cos \psi  +   \wt b_1   \sin \x \sin \psi$ (i.e.  a potential   \eqref{explicita**} with    $\wt a_n,  \wt b_n  = 0$, $n \geq 2$), there exists a choice of  ``rotated'' spherical coordinates $\big(\xi' \= \xi(\a^o_1, \a^o_2)$, $\psi' \= \psi(\a^o_1, \a^p_2)\big)$  in which   $\f$ takes  the form 
\beq \label{oldkerr} \f(\xi', \psi') = a \cos \xi'\ , \qquad \text{with}\ \ a = \sqrt{\wt a_0^2 + \wt a_1^2 + \wt b_1^2}\ .\eeq
Plugging \eqref{oldkerr} into \eqref{6666},  we   see that the $3$-parameter family of (local) potentials $\f$  with    $\wt a_n = \wt b_n = 0$, $n \geq 2$, correspond to  Lorentzian metrics that are locally isometric to the classical Kerr metrics  \eqref{standardKerr}.  The  parameters of such a  family of potentials (i.e. the  constants $\wt a_0$, $\wt a_1$,  $\wt b_1$) correspond to the three space-like components (with respect to a fixed  frame) of the    total angular momentum $\overset{\to}\cJ = m \overset{\to} a  $  of the gravitational field. \par
\smallskip
From this  we get  that, whenever   $\wt a_n $, $\wt b_n$, $n \geq 2$, are  equal to $0$, the corresponding  metrics -- which  a priori  are defined just on an open subset of $S^2_{(+)} \times \bR_+ \times \bR$ --    extend in a real analytic way to  the whole  $(S^2 \times \bR_+) \times \bR =( \bR^3 \setminus \{0\}) \times \bR$.  All such extensions are isometric to  the classical Kerr metrics for appropriate values for   $m$ and  $\cJ = m |\overset{\to} a|$.  On the contrary,  for    non-zero choices   of a finite collection of   constants $\wt a_{n_\ell}, \wt b_{n_\ell}$,   $2 \leq n_1 <  \ldots < n_s$,  by just  looking  at the  terms corresponding to the  highest integer $n_s$,   one   realises   that  the corresponding    metrics $g^{(m)}$  do not extend to the whole  $( \bR^3 \setminus \{0\}) \times \bR $.   Actually,   for a few  examples of  this kind,  the explicit computations  of the Carminati-McLenaghan invariants  of the Riemann curvature \cite{CM} indicate  that  such  metrics develop  intrinsic curvature singularities in approaching the south pole of $S^2$ and it is reasonable to expect that this occurs for all of them. \par
\smallskip
These remarks  give evidence that   the  potentials   \eqref{explicita**} provide  a very large number of Ricci flat metrics  that  are not   locally isometric to any of the   classical Kerr metrics.  It would be  interesting to  study  their  geometric properties  and explore the possibilities  to generalise  and/or glue them with  appropriate solutions of the Einstein equations in order to obtain  new solutions to the  Einstein equations, representing gravitational fields  with non-zero total angular momentum and determined by   stress-energy tensors    that vanish  on  large (but not spherically symmetric) regions of the   space-time. We hope to address   these issues in the near future. 
\par
\smallskip
 As we mentioned in the Introduction,   it is reasonable to     expect that also when there is an {\it infinite} collection of  non-vanishing coefficients $\wt a_{n_\ell}$, $\wt b_{n_\ell} \neq 0$, $n_\ell \geq 2$, the corresponding   metric  admits no real analytic extension to  $(\bR^3 \setminus\{0\}) \times \bR$. A proof of  this quite  sensible  conjecture  would give  a new   characterisation of the  classical Kerr metrics (they  would be the only  metrics   \eqref{6666}  that are  extendible to   $(\bR^3 \setminus\{0\}) \times \bR$),   with   interesting  relations with   the  so far known    uniqueness theorems for the Kerr metrics (see e.g. \cite{Ro1}).
\par
\smallskip
\subsubsection{A quick look at   the metrics of the  Kerr families    with  $\k = -1$} Let us now focus on the   potentials in    \eqref{explicita}   for $\k = -1$.   Replacing  the polar coordinates $(\gr, w)$ by     the coordinates $(\xi,  \psi)$   given by   the   relations 
 $\gr = \tanh \left( \frac{\xi}{2}\right)$,  $w = \psi$, 
  the metric   $\wt g_o$  and the potentials $\f$  take the  form 
    \beq\label{explicita*}
     \begin{split}
     & \wt g_o  = d \xi^2 + \sinh^2 \xi d \psi^2\ ,\\
     & \f(\xi, \psi) = \wt a_0 \cosh \xi   +   \wt a_1   \sinh \x \cos \psi  +   \wt b_1   \sinh \x \sin \psi  +   \\
     & \hskip 1 cm + \sum^\infty_{n = 2} \tanh^n \left( \frac{\xi}{2}\right)  \left( \cosh \xi  +  n -1 \right)  \big(  \wt a_n \cos(n \psi) +   \wt b_n \sin(n \psi)\big)\ .     \end{split}
     \eeq
 Note also that, if the coefficients $\wt a_\ell$ are taken sufficiently small so that $\f$ is sufficiently close to $0$ on some relatively compact subset of $\bD$, the corresponding Lorentzian metric is close to the metric \eqref{6666} with $\f = 0$, namely to 
     \beq \label{6666*} g^{(m)} =   \frac{r^2}{2} ( d \xi^2 + \sinh^2 \xi d \psi^2)  
    +  d v \vee \bigg( d r + \frac{ r + m  }{2 r} d v \bigg)\ ,
\eeq

This metric can be taken as an analog of the classical Schwarzschild metric, but with  the crucial difference  of having   the  group $\bR \times \SO^o_{1, 2}$ instead of  $\bR \times \SO_3$ as connected component of the group of  isometries.  This fact gives evidence  that for    potentials sufficiently close to  $\f = 0$,  if $\k = -1$, the corresponding Lorentzian metrics  cannot be   locally isometric to the classical Kerr metrics.  A rigorous proof would however demand a detailed study of the isometric invariants of those  metrics, a task which goes beyond the scope of this paper.  
       \par 
\medskip
For what concerns the geometric properties,    the Ricci flat metrics \eqref{6666}  with  $\k = -1$ mirror certain crucial    features of those with $\k = 1$. Indeed,  in both settings, each   metric  is invariant under the action of the vector field $\bT \=   \qs_o -  B \ps_o$. 
By  the discussion in Remark \ref{importantrem},   with respect to   the background metric $\h$,   the  vector field $\bT$ is time-like if $\k = 1$, while it  is space-like if   $\k = -1$. Moreover, passing from the case 
 $\k = 1$ to the case $\k = -1$,  according to   \eqref{back1}  the signature of  $\h$   changes from $(-, - , -, +)$ to $(+, + , -, +)$. This means that the roles of the last two coordinates appearing in \eqref{back1} switch from space-like to time-like  and vice versa. \par
This phenomenon greatly helps in   guessing      features of the metrics  with  $\k = -1$. For instance,  let us first  recall that  the  classical Kerr metrics  (i.e. the cases with  $\k = 1$ and $\wt a_n = \wt b_n = 0$, $n \geq 2$)  determine  distinguished  hypersurfaces  (characterised by an equation in the space-like coordinate $r$) which are  time-invariant  and have the fundamental role of being  {\it event horizons}. According to  the above described switch  between time-like and space-like roles,   it is reasonable to  expect that the analogs  in the class with  $\k = -1$ determine   distinguished hypersurfaces (characterised by an equation in the  {\it time-like} coordinate $r$), which are  invariant under a space-like vector field and  with  features   similar to  those of the  event horizons, provided that   the meanings of  the words  ``time'' and   ``radius''  are switched at all places.  \par
\smallskip
Guided by  this kind of intuition, while the Kerr metrics tend to  stationary  spherically symmetric weak  gravitational fields  at large distances,  the newly determined Ricci flat  metrics  on $M^{\text{Lob}} \= (\bD \times \bR) \times \bR$ are expected to  tend to radially-independent small deformations of the Minkowski  metrics for large values of  the  ``time''.  Being not asymptotically spherically symmetric, these  Ricci flat metrics  are not  appropriate to  model the  gravity generated by  a  localised mass. Nonetheless, they might be considered  as   non-static cosmological models
or, more precisely, as  geometric  limits  (for  the  stress-energy tensors tending to $0$)   of a large  class of   non-static cosmological models. 

\par
\medskip
\appendix
\section{Details of the  computations that yield  \eqref{5.32ter}} \label{appendix-1}
Let us  recall   the following  properties and relations between  the Christoffel symbols and  components $J^i_j$ of   the complex structure: 
\begin{itemize}
\item[(a)]   $J^m_m = 0$  (it is the trace of a complex structure); 
\item[(b)]   $ \Ga i j r - \Ga j i r - c_{ij}^r = 0 $ ($\Ga i j r$ are  Christoffel symbols of the torsion free connection $\n^{g_o}$);  
\item[(c)]  $E_m(J^\ell_i) =  - \Ga m r \ell J^r_i + \Ga m i r J^\ell_r,
$ (in fact, since $(N, J, g)$ is K\"ahler and hence   $\n^{g_o} J = 0$, we have that  
$$E_m(J^\ell_i) =  \big( \n^{g_o}_{E_m} E^\ell\big)\left(J(E_i)\right) + E^\ell \left(\big(\n^{g_o}_{E_m} J\big)(E_i)\right) 
+  E^\ell\left(J(\n^{g_o}_{E_m}E_i) \right) = - \Ga m r \ell J^r_i + \Ga m i r J^\ell_r\ \bigg).
$$
\end{itemize}
Using   \eqref{theEhat*}, \eqref{519} -- \eqref{6.17}  and the above  relations (a) -- (c),   the Ricci curvature components \eqref{6.21}  can be written as \\
\scalebox{0.82}{\vbox{
\begin{align*}
\nonumber  \Ric(\wh E_i, \ps_o) =   \wh E_m (\GGa i {\ps_o} m) &-  \wh E_i (\GGa m {\ps_o} m) -  \GGa m {\ps_o}  \ell \GGa i  \ell  m  +  \GGa i {\ps_o} \ell  \GGa m  \ell  m  - c^r_{mi} \GGa r {\ps_o} m  = \\
\nonumber  & = E_m \left(  \frac{J^m_{i} }{4 \s}  +  \ps_o(\s) \frac{\d^m_i}{2 \s}\right) + B  J_m^\ell \f_\ell \ps_o \left(  \frac{J^m_{i} }{4 \s}  +  \ps_o(\s) \frac{\d^m_i}{2 \s}\right)  - \\
 \nonumber &  \hskip 1 cm -   E_i \left(\xcancel{  \frac{J^m_{m} }{4 \s} } +  \ps_o(\s) \frac{\d^m_m}{2 \s}\right) - B  J_i^\ell \f_\ell \ps_o \left(  \xcancel{\frac{J^m_{m} }{4 \s} } +  \ps_o(\s) \frac{\d^m_m}{2 \s}\right)  -\\
\nonumber  & \hskip 1 cm -   \left(  \frac{J^\ell_{m} }{4 \s}  +  \ps_o(\s) \frac{\d^\ell_m}{2 \s}\right) \GGa i  \ell  m   
   +\left(   \frac{J^\ell_{i} }{4 \s}  +  \ps_o(\s) \frac{\d^\ell_i}{2 \s}  \right) \GGa m  \ell  m  - c^r_{mi} \GGa r {\ps_o} m = \\
   \nonumber  & =     \frac{E_m \left(J^m_{i} \right)}{4 \s}  -  \frac{J^m_{i} \s_m}{4 \s^2}   +\frac{1}{2 \s}  E_i \left( \ps_o(\s) \right)-  \frac{1}{2 \s^2} \ps_o(\s)\s_i - \\
\nonumber  & \hskip 1 cm  - B  J^m_{i} J_m^\ell \f_\ell   \frac{\ps_o \left(\s\right) }{4 \s^2}    + B  J_i^\ell \f_\ell \frac{1}{2 \s}  \ps_o \left(\ps_o(\s) \right) - B  J_i^\ell \f_\ell \frac{1}{2 \s^2} \left( \ps_o \left(\s\right) \right)^2  - \\
\nonumber  & \hskip 1 cm   -   E_i \left( \ps_o(\s)\right) \frac{(n-2)}{2 \s} +    \s_i \ps_o(\s) \frac{(n-2)}{2 \s^2} -  \\
\nonumber & \hskip 1 cm - B J_i^\ell \f_\ell \frac{(n-2)}{2 \s} \ps_o\left( \ps_o(\s)\right)  + B  J_i^\ell \f_\ell \frac{(n-2)}{2 \s^2}  \left( \ps_o(\s) \right)^2 +\\
\nonumber  & \hskip 1 cm -  \frac{J^\ell_{m} }{4 \s}\GGa i  \ell  m    -   \ps_o(\s) \frac{1}{2 \s} \GGa i  m  m     
  +   \frac{J^\ell_{i} }{4 \s} \GGa m  \ell  m +  \ps_o(\s) \frac{1}{2 \s} \GGa m  i  m   - c^r_{mi} \GGa r {\ps_o} m = \\
  \nonumber  & =  \underset{ E_m(J^m_i)}{\underbrace{ \xcancel{ - \frac{\Ga m r m J^r_i}{4 \s} } +\frac{\Ga m i r J^m_r}{4 \s} }}    -  \frac{J^\ell_{i} \s_\ell}{4 \s^2}   - (n-3)\frac{E_i \left( \ps_o(\s) \right)}{2 \s} + (n-3) \frac{ \ps_o(\s)\s_i}{2 \s^2} + \\
\nonumber  & \hskip 1 cm  + B   \frac{\f_i \ps_o \left(\s\right)}{4 \s^2}     
 - B   J_i^\ell  \f_\ell \ps_o\left( \ps_o(\s)\right) \frac{(n-3)}{2 \s}  + B   J_i^\ell \f_\ell  \left(\ps_o(\s)\right)^2 \frac{(n-3)}{2 \s^2}  -\\ \nonumber & \hskip 1 cm - \frac{J^\ell_{m} \Ga i \ell m}{4 \s} 
- \xcancel{\frac{J^m_{m} }{8 \s^2}    \s_i }  -   \frac{J^\ell_{i} \s_\ell }{8\s^2}   
 +   \frac{\overset{- J^k_i}{\overbrace{J^\ell_{m} g_{i\ell}     g^{mk}}  }\s_k}{8 \s^2} -\\
  \nonumber  & \hskip 1 cm - \xcancel{ \frac{J^m_{m} }{4 \s}\frac{  B  \ps_o(\s)\f_r }{2 \s } J^r_i} + \frac{  B  \ps_o(\s) \f_i }{8 \s^2 } 
 +  \frac{ B  \ps_o(\s) \f_r }{8 \s^2 } \overset{- J_i^k}{\overbrace{  J^\ell_m g_{i \ell} g^{m k}}} J_k^r - \\
 \nonumber  & \hskip 1 cm  -   \frac{  \ps_o(\s) \Ga i m m }{2 \s} 
 - \xcancel{ \frac{(n-2) \s_i  \ps_o(\s) }{4 \s^2}  }    -  \xcancel{ \frac{ \s_i  \ps_o(\s)  }{4 \s^2 }  }
 +   \xcancel{   \frac{ \ps_o(\s)  \s_i}{4 \s^2} } -   \xcancel{  \frac{  (n-2)     B   J^\ell_i  \f_\ell \left(\ps_o(\s)\right)^2 }{4 \s^2 }  }
   +
 \end{align*}
}
}
\par
\scalebox{0.82}{\vbox{
 \begin{align*}
 \nonumber   & \hskip 1 cm + \xcancel{ \frac{J^\ell_{i} \Ga m \ell m}{4 \s} }
+ \xcancel{  \frac{J^\ell_{i}  \s_\ell  }{8 \s^2}  }  +   \frac{(n-2)J^\ell_{i} \s_\ell  }{8 \s^2} 
 -\xcancel{  \frac{J^\ell_{i}  \s_\ell }{8 \s^2}  }   -   \frac{ (n-2)    B  \ps_o(\s) \f_i }{8 \s^2 }
 +  \\
\nonumber    & \hskip 1 cm +  \frac{ \ps_o(\s)\Ga m i m }{2 \s} 
+ \xcancel{ \frac{ \ps_o(\s)  \s_i }{4 \s^2 } }  +  \xcancel{\frac{  \ps_o(\s) \s_i (n-2) }{4 \s^2} }
 -    \xcancel{\frac{  \ps_o(\s) \s_i}{4 \s^2} }  +\\
\nonumber  & \hskip 1 cm + \xcancel{ \frac{(n-2) \left(\ps_o(\s)\right)^2    B   J_i^\ell\f_\ell }{4 \s^2 }  }
 -   c^r_{mi}    \frac{J^m_{r} }{4 \s}   -  \frac{c^m_{mi}   \ps_o(\s)}{2 \s} = \\
\nonumber  & =   \xcancel{  \frac{J^m_\ell}{4 \s} \left(\Ga m i \ell -  \Ga i m \ell  - c^\ell_{mi}  \right)  }   +   (n-6) \frac{J^\ell_{i} \s_\ell}{ 8\s^2}          -  \frac{(n-3)E_i \left( \ps_o(\s) \right)}{2 \s}    +  \frac{(n-3) \ps_o(\s) \s_i}{ 2 \s^2}       + \\
\nonumber  & \hskip 1 cm   - \frac{(n-6)   B  \f_i \ps_o \left(\s\right)}{8 \s^2}       -\frac{(n - 3)   B  J_i^\ell \f_\ell }{ 2 \s}\ps_o \left(\ps_o(\s) \right)+   \frac{(n - 3)   B  J_i^\ell \f_\ell }{ 2 \s^2} \left( \ps_o \left(\s\right) \right)^2      + \\
\nonumber  & \hskip 1 cm +\xcancel{  \frac{J^m_\ell}{4 \s} \left(\Ga m i \ell -  \Ga i m \ell  - c^\ell_{mi}  \right)  }   
  -\xcancel{ \frac{\ps_o(\s)}{2 \s} \left( \Ga m i m -  \Ga i m m 
  -  c^m_{mi}  \right)}.\
\end{align*}
}
}
\\
From the last equality and    factoring  out $\frac{n-3}{8\s^2}$, \eqref{5.32ter} follows .  \par
\medskip
\section{The proof of Lemma \ref{lemmappend}} \label{appendixA}
In this appendix we  provide the  details of the computations that  prove Lemma \ref{lemmappend}. In what follows we  denote $\lambda = \D \f$ and we use  the   notation $\d \Ga i j k$ for  the differences 
  $\dGa i j k = \GGa i j k - \Ga i j k$.  
  Using the symmetry  $\Ga i j k = \Ga j i k$ (because the  $E_i$ commute) and the identity
  \begin{multline} 0 =E_i\left(g(E_j, E_k)\right) - g(\n_{E_i}E_j, E_k)  -  g( E_j, \n_{E_i}E_k) = \\
  = E_i(\lambda) \d_{jk} - \lambda \Ga i j \ell \d_{\ell k}  - \lambda \Ga i k \ell \d_{j \ell} =  E_i(\lambda) \d_{jk} - \lambda \Ga i j k   - \lambda \Ga i k j \ ,\end{multline}
 we obtain the following expressions for the $\Ga i j k$   and for the differences $\dGa i j k$ \\
\beq \label{6.78} 
\begin{split} &  \Ga 1 1 1 = \frac{1}{2} \frac{\lambda_x}{\lambda} \ , \qquad 
 \Ga 2 2 2 = \frac{1}{2} \frac{\lambda_y}{\lambda}\ ,\qquad
  \Ga 2 1 1 =   \Ga 2 2 2\ , \qquad
   \Ga 1 2 2 =   \Ga 1 1 1 \ ,\\[10pt]
&\Ga 1 2 1 = \Ga 2 2 2\ , \qquad
\Ga 1 1 2 =  -\Ga 2 2 2\ , \qquad
 \Ga 2 1 2 =  \Ga 1 1 1\ ,\qquad
  \Ga 2 2 1 =  -  \Ga 1 1 1 \ .
\end{split}
  \eeq
\beq \begin{split}
&  \dGa 1 1 1 = \frac{\s_x}{2\s } + B \frac{ \s_t \f_y}{2 \s} \ , \qquad 
 \dGa 2 2 2 =\frac{1}{2 \s} \s _y - B \frac{\s_t \f_x}{2 \s}\ ,\qquad 
  \dGa 2 1 1 =   \dGa 2 2 2 \ ,\qquad
   \dGa 1 2 2 =   \dGa 1 1 1 \ ,\\[10pt]
& \dGa 1 2 1 =  \dGa 2 2 2 \ ,\qquad
\dGa 1 1 2 =  - \dGa 2 2 2\ ,\qquad 
 \dGa 2 1 2 =  \dGa 1 1 1 \ ,\qquad
  \dGa 2 2 1 =  -  \dGa 1 1 1 \ .
  \end{split}
\eeq
Note that each of the  Christoffel symbols $\Ga i j k$ (resp. $\GGa i j k$) is (up to a sign) equal either to $\Ga 1 1 1 = \frac{\lambda_x}{\lambda} $ (resp. $\GGa 1 1 1$) or to $\Ga 2 2 2   = \frac{\lambda_y}{\lambda} $ (resp. $\GGa 2 2 2$).   
  Let us  now focus on the equation \eqref{6.75}. Expanding the summations,  the left hand side becomes  (in the formula, pairs of terms that cancel out   are marked with equal underlying  boxed symbols)\par
  \scalebox{0.8}{\vbox{
  \begin{align*}
  \Ric_{12}=& \underset{\boxed{A}}{\underbrace{\xcancel{\left(E_1 + J_1^\ell E_\ell(\f) B \frac{\p}{\p r}\right)   (\GGa 1 2 1)}}} -   \underset{\boxed{A}}{\underbrace{\xcancel{\left(\frac{\p}{\p x} + B \f_y \frac{\p}{\p r}\right)(\GGa 1 2 1)}} }-  \GGa 1 2  \ell \GGa 1  \ell  1    -\underset{\boxed{B}}{\underbrace{ \xcancel{ \GGa 1 2  {\qs_o} \GGa 1 {\qs_o}  1}} }+\\
  &\hskip 1.5 cm  +   \GGa 1 2 \ell  \GGa 1 \ell  1+ \underset{\boxed{D}}{\underbrace{ \xcancel{\GGa 1 2 {\qs_o}  \GGa 1 {\qs_o}  1}}}  - \lambda  \GGa {\qs_o} 2 2   + B\frac{\p}{\p r} (\GGa 1 2 {\qs_o}) - \underset{\boxed{C}}{\underbrace{ \xcancel{\GGa {\qs_o} 2  1 \GGa 1  1 {\qs_o}} } }- \underset{\boxed{D}}{\underbrace{\xcancel{\GGa {\qs_o} 2  2 \GGa 1  2 {\qs_o}} }} +\\
  & -\left(E_2 + J_2^\ell E_\ell(\f)    B \frac{\p}{\p r}\right)   (\GGa 1 2 2) -  \left(\frac{\p}{\p x} + B \f_y \frac{\p}{\p r}\right)(\GGa 2 2 2) -  \GGa 2 2  \ell \GGa 1  \ell  2    -  \underset{\boxed{C}}{\underbrace{  \xcancel{\GGa 2 2  {\qs_o} \GGa 1 {\qs_o}  2 }}}  +   \GGa 1 2 \ell  \GGa 2 \ell  2+ \underset{\boxed{B}}{\underbrace{ \xcancel{ \GGa 1 2 {\qs_o}  \GGa 2 {\qs_o}  2}}}   = 
\\
  =& - \xcancel{ \GGa 1 2  1 \GGa 1  1  1}   - \xcancel{ \GGa 1 2  2 \GGa 1  2  1}  +  \xcancel{ \GGa 1 2 1  \GGa 1 1  1}+  \xcancel{ \GGa 1 2 2  \GGa 1 2  1}    -  \lambda  \GGa {\qs_o} 2 2   + B\frac{\p}{\p r} (\GGa 1 2 {\qs_o})  +\\
  &+ \left(\frac{\p}{\p y} - B \f_x     \frac{\p}{\p r}\right)   (\GGa 1 2 2) -  \left(\frac{\p}{\p x} + B \f_y \frac{\p}{\p r}\right)(\GGa 2 2 2) -  \GGa 2 2  1 \GGa 1  1  2  - \xcancel{  \GGa 2 2  2 \GGa 1  2  2 }  +   \GGa 1 2 1  \GGa 2 1  2+ \xcancel{  \GGa 1 2 2  \GGa 2 2  2} =
 \end{align*}
  }}\\[10 pt]
  \scalebox{0.8}{\vbox{
  \begin{align*}
  &= -  \lambda  \GGa {\qs_o} 2 2  + B\underset{= \frac{\p}{\p r}\left(\frac{\lambda}{2}\right) = 0}{\xcancel{\frac{\p}{\p r} (\GGa 1 2 {\qs_o})}} + \frac{\p}{\p y}(\Ga 1 2 2) +\frac{\p}{\p y}(\d \Ga 1 2 2) - \f_x B \frac{\p}{\p r}(\d \Ga 1 2 2)    -\frac{\p}{\p x}(\Ga 2 2 2)  -\frac{\p}{\p x}(\d\Ga 2 2 2)-\\
  & - B \f_y \frac{\p}{\p r}(\d\Ga 2 2 2) 
  -\Ga 2 2  1\Ga 1  1  2 -\Ga 2 2  1 \d \Ga 1  1  2 -\d \Ga 2 2  1 \Ga 1  1  2 -\d \Ga 2 2  1\d \Ga 1  1  2    + \Ga 1 2 1 \Ga 2 1  2+\Ga 1 2 1\d \Ga 2 1  2+\d  \Ga 1 2 1\Ga 2 1  2+ \d  \Ga 1 2 1 \d \Ga 2 1  2
 \end{align*}
 \begin{align*}
  =&\  \underset{= \Ric_{12}^N}{\underbrace{  \frac{\p}{\p y}(\Ga 1 2 2)  -\frac{\p}{\p x}(\Ga 2 2 2) -\Ga 2 2  1\Ga 1  1  2+ \Ga 1 2 1 \Ga 2 1  2}}  -  \lambda  \GGa {\qs_o} 2 2 +\frac{\p}{\p y}(\d \Ga 1 2 2)- B \f_x    \frac{\p}{\p r}(\d \Ga 1 2 2)      -\frac{\p}{\p x}(\d\Ga 2 2 2)  - B \f_y \frac{\p}{\p r}(\d\Ga 2 2 2) -\\
  &
  \xcancel{ -\Ga 2 2  1 \d \Ga 1  1  2} -\xcancel{\d \Ga 2 2  1 \Ga 1  1  2} -\xcancel{\d \Ga 2 2  1\d \Ga 1  1  2 }    +\xcancel{\Ga 1 2 1\d \Ga 2 1  2}+\xcancel{\d  \Ga 1 2 1\Ga 2 1  2}+ \xcancel{\d  \Ga 1 2 1 \d \Ga 2 1  2} = \\
   =&\ \Ric_{12}^N  - B \lambda  \frac{\s_t}{2\s}  +\frac{\p}{\p y}(\d \Ga 1 1 1)- B\f_x     \frac{\p}{\p r}(\d \Ga 1 1 1)      -\frac{\p}{\p x}(\d\Ga 2 2 2) - B \f_y \frac{\p}{\p r}(\d\Ga 2 2 2) \\
     =&\ \Ric_{12}^N   - B \lambda  \frac{\s_t}{2\s}  + \\
  & +\frac{\p}{\p y}\left (\frac{\s_x}{2\s } + B \frac{ \s_t \f_y}{2 \s}\right ) - B \f_x  \frac{\p}{\p r}\left (\frac{\s_x}{2\s } + B \frac{ \s_t \f_y}{2 \s}\right)      -\frac{\p}{\p x}\left (\frac{\s _y}{2 \s}  - B \frac{\s_t \f_x}{2 \s}\right) - B \f_y \frac{\p}{\p r}\left (\frac{\s _y}{2 \s}  - B \frac{\s_t \f_x}{2 \s}\right )  = \\
   = &\ \Ric_{12}^N  - B\lambda  \frac{\s_t}{2\s}  +\xcancel{\frac{\s_{xy}}{2\s}}-\xcancel{\frac{\s_x \s_y}{2\s^2}}+\xcancel{B \frac{\s_{ty}}{2\s}\f_y}-\xcancel{B\frac{\s_t \s_y}{2\s^2}\f_y}+ B \frac{\s_t}{2\s}\f_{yy}	- \\
   & -\xcancel{B\f_x \frac{\s_{xt}}{2\s}} +\xcancel{B\f_x \frac{\s_x \s_t}{2\s^2}	} -\xcancel{B^2 \f_x \f_y \frac{\s_{tt}}{4\s}}  +\xcancel{B^2 \f_x \f_y \frac{\s_t^2}{4\s^2} } -  \\
    &-\xcancel{\frac{\s_{xy}}{2\s}}+\xcancel{\frac{\s_y \s_x}{2\s^2}}+\xcancel{B \frac{\s_{tx}}{2\s}\f_x}-\xcancel{B\frac{\s_t \s_x}{2\s^2}\f_x}+ B \frac{\s_t}{2\s}\f_{xx} - \xcancel{B \f_y \frac{\s_{yt}}{2\s}}+\xcancel{B \f_y \frac{\s_y \s_t}{2\s^2}}+ \xcancel{B^2 \f_y \f_x \frac{\s_{tt}}{4\s}}-\xcancel{B^2\f_y \f_x \frac{\s_t^2}{4\s^2}} = \\
  = &\ \Ric_{12}^N  - B\lambda  \frac{\s_t}{2\s} 
  + B \underset{= \lambda}{\underbrace{(\f_{yy} + \f_{xx}) }}  \frac{\s_t}{2\s}  =  \Ric_{12}^N \ .
 \end{align*}
 }}\\
This means that \eqref{6.75} is  equivalent to the  equation  $  \Ric^N(E_1, E_2) = 0$.  The latter is identically satisfied   because, by assumption, $g_o$ is K\"ahler and $E_2 = J E_1$. 
We now need to show that   the equations \eqref{6.74} and \eqref{6.76} are both  equivalent to \eqref{commonequation}. Expanding  the summations in the indices $m = 1,2$, the left hand sides of those equations become\par
\scalebox{0.8}{\vbox{
\begin{align*} 
 \Ric_{11}  = &\  \xcancel{ \frac{\p}{\p x}(\GGa 1 1 1)} + J_1^\ell E_\ell(\f)   B\frac{\p}{\p r}(\GGa 1 1 1)    - \xcancel{ \frac{\p}{\p x}(\GGa 1 1 1) }- B \f_y \frac{\p}{\p r}(\GGa 1 1 1) - \xcancel{ \GGa 1 1  \ell \GGa 1  \ell  1 }   -  \xcancel{\GGa 1 1  {\qs_o} \GGa 1 {\qs_o}  1} +\\
  &\hskip 1.5 cm  +  \xcancel{ \GGa 1 1 \ell  \GGa 1 \ell  1}+   \xcancel{\GGa 1 1 {\qs_o}  \GGa 1 {\qs_o}  1}   - \lambda   \GGa {\qs_o} 1 2   + B \frac{\p}{\p r} (\GGa 1 1 {\qs_o}) -  \GGa {\qs_o} 1  \ell \GGa 1  \ell {\qs_o} +\\
  +&  \frac{\p}{\p y}(\GGa 1 1 2) + J_2^\ell E_\ell(\f)    B\frac{\p}{\p r}(\GGa 1 1 2)    -  \frac{\p}{\p x}(\GGa 2 1 2) - B \f_y \frac{\p}{\p r}(\GGa 2 1 2) -  \GGa 2 1  \ell \GGa 1  \ell  2    -  \GGa 2 1  {\qs_o} \GGa 1 {\qs_o}  2  +   \GGa 1 1 \ell  \GGa 2 \ell  2+   \GGa 1 1 {\qs_o}  \GGa 2 {\qs_o}  2   =  \end{align*}
  \begin{align*}
  =&\ \xcancel{ \f_y   B \frac{\p}{\p r}(\GGa 1 1 1) }   -\xcancel{   B \f_y \frac{\p}{\p r}(\GGa 1 1 1) }   -   \lambda   \GGa {\qs_o} 1 2  + B \frac{\p}{\p r} (\GGa 1 1 {\qs_o}) -\underset{\boxed{B}}   {\xcancel{\GGa {\qs_o} 1  1 \GGa 1  1 {\qs_o}} }-\underset{\boxed{A}} {\xcancel{ \GGa {\qs_o} 1  2 \GGa 1  2 {\qs_o} }}+  \frac{\p}{\p y}(\GGa 1 1 2) -\f_x  B \frac{\p}{\p r}(\GGa 1 1 2)    - \\
  -& \ \frac{\p}{\p x}(\GGa 2 1 2) - B \f_y \frac{\p}{\p r}(\GGa 2 1 2) -  \GGa 2 1  1 \GGa 1  1  2  -  \GGa 2 1  2 \GGa 1  2  2     -\underset{\boxed{-A}} {\xcancel{  \GGa 2 1  {\qs_o} \GGa 1 {\qs_o}  2}} +  \GGa 1 1 1  \GGa 2 1  2 +   \GGa 1 1 2  \GGa 2 2  2+  \underset{\boxed{B}}  { \xcancel{ \GGa 1 1 {\qs_o}  \GGa 2 {\qs_o}  2}}  
  \end{align*}
  }}\par
  \scalebox{0.8}{\vbox{
  \begin{align*} 
 \Ric_{22}  = &\   \frac{\p}{\p x}(\GGa 2 2 1) + J_1^\ell E_\ell(\f)   B\frac{\p}{\p r}(\GGa 2 2 1)    -  \frac{\p}{\p y}(\GGa 1 2 1) + B \f_x \frac{\p}{\p r}(\GGa 1 2 1) - \GGa 1 2  \ell \GGa 2  \ell  1    -  \GGa 1 2  {\qs_o} \GGa 2 {\qs_o}  1 +\\
  + & \    \GGa 2 2 \ell  \GGa 1 \ell  1+   \GGa 2 2 {\qs_o}  \GGa 1 {\qs_o}  1   +  \lambda   \GGa {\qs_o} 1 2  + B \frac{\p}{\p r} (\GGa 2 2 {\qs_o}) -  \GGa {\qs_o} 2  \ell \GGa 2  \ell {\qs_o} +  \xcancel{ \frac{\p}{\p y}(\GGa 2 2 2)} + J_2^\ell E_\ell(\f)    B\frac{\p}{\p r}(\GGa 2 2 2)    -  \\
  -&\  \xcancel{\frac{\p}{\p y}(\GGa 2 2 2)} + B \f_x \frac{\p}{\p r}(\GGa 2 2 2) - \xcancel{ \GGa 2 2  \ell \GGa 2  \ell  2}    - \xcancel{ \GGa 2 2  {\qs_o} \GGa 2 {\qs_o}  2 }+  \xcancel{ \GGa 2 2  \ell  \GGa 2 \ell  2} +   \xcancel{\GGa 2 2 {\qs_o}  \GGa 2 {\qs_o}  2}   =  \\
= &   \frac{\p}{\p x}(\GGa 2 2 1) + \f_y   B\frac{\p}{\p r}(\GGa 2 2 1)    -  \frac{\p}{\p y}(\GGa 1 2 1) + B \f_x \frac{\p}{\p r}(\GGa 1 2 1) - \GGa 1 2  \ell \GGa 2  \ell  1    - \xcancel{ \GGa 1 2  {\qs_o} \GGa 2 {\qs_o}  1 }+\\
  + &\   \GGa 2 2 \ell  \GGa 1 \ell  1+  \xcancel{ \GGa 2 2 {\qs_o}  \GGa 1 {\qs_o}  1}   + \lambda   \GGa {\qs_o} 2 1   + B \frac{\p}{\p r} (\GGa 2 2 {\qs_o}) - \xcancel{ \GGa {\qs_o} 2  1 \GGa 2  1 {\qs_o}} - \xcancel{ \GGa {\qs_o} 2  2 \GGa 2  2 {\qs_o} }  -  \xcancel{\f_x   B\frac{\p}{\p r}(\GGa 2 2 2)}    + \xcancel{  B \f_x \frac{\p}{\p r}(\GGa 2 2 2) }. \ 
  \end{align*}
  }}\\
After all cancellations,   the expressions reduce to \par
 \scalebox{0.8}{\vbox{
\begin{align} 
 \nonumber \Ric_{11}  = & -  \lambda   \GGa {\qs_o} 1 2   + B \frac{\p}{\p r} (\GGa 1 1 {\qs_o})  +   \frac{\p}{\p y}(\GGa 1 1 2) - B \f_x   \frac{\p}{\p r}(\GGa 1 1 2)    -  \frac{\p}{\p x}(\GGa 2 1 2) - B \f_y \frac{\p}{\p r}(\GGa 2 1 2) -  \\
\label{6.79} & - \GGa 2 1  1 \GGa 1  1  2  -  \GGa 2 1  2 \GGa 1  2  2      +   \GGa 1 1 1  \GGa 2 1  2 +   \GGa 1 1 2  \GGa 2 2  2\\
 \nonumber \Ric_{22}  =  &   \lambda   \GGa {\qs_o} 2 1   + B \frac{\p}{\p r} (\GGa 2 2 {\qs_o}) +   \frac{\p}{\p x}(\GGa 2 2 1) + B \f_y   \frac{\p}{\p r}(\GGa 2 2 1)    -  \frac{\p}{\p y}(\GGa 1 2 1) + B \f_x \frac{\p}{\p r}(\GGa 1 2 1) -\\ 
\label{6.80} & -  \GGa 1 2  1 \GGa 2  1  1   -  \GGa 1 2  2 \GGa 2  2  1   +    \GGa 2 2 1  \GGa 1 1  1+     \GGa 2 2 2  \GGa 1 2  1\ .
  \end{align}
  }}
  \\
  Let us now expand these expressions using the decompositions   $\GGa i j k = \Ga i j k + \dGa i j k$.  We get:  \\
   \scalebox{0.8}{\vbox{
  \begin{align*}
   \Ric_{11} & =   \frac{\p}{\p y}(\Ga 1 1 2)   -  \frac{\p}{\p x}(\Ga 2 1 2)  -  \Ga 2 1  1 \Ga 1  1  2  -  \Ga 2 1  2 \Ga 1  2  2      +   \Ga 1 1 1  \Ga 2 1  2 +   \Ga 1 1 2  \Ga 2 2  2 +  \frac{\p}{\p y}\dGa 1 1 2   -  \frac{\p}{\p x}\dGa 2 1 2 - \\
      &  -\Ga 2 1  1\dGa 1  1  2    -  \Ga 2 1  2 \dGa 1  2  2  +   \Ga 1 1 1  \dGa 2 1  2  +   \Ga 1 1 2   \dGa 2 2  2 -  \dGa 2 1  1  \Ga 1  1  2  - \dGa 2 1  2 \Ga 1  2  2    + \dGa 1 1 1   \Ga 2 1  2 +   \dGa 1 1 2  \Ga 2 2  2 + \\
       & - \dGa 2 1  1\dGa 1  1  2    -  \dGa 2 1  2 \dGa 1  2  2   +   \dGa 1 1 1  \dGa 2 1  2  +   \dGa 1 1 2   \dGa 2 2  2 -   \lambda   \GGa {\qs_o} 1 2  + B \frac{\p}{\p r} \left(\GGa 1 1 {\qs_o} -  \f_y \GGa 2 1 2    -  \f_x \GGa 1 1 2   \right) \ , \\
          \Ric_{22} & =   \frac{\p}{\p x}(\Ga 2 2 1)   -  \frac{\p}{\p y}(\Ga 1 2 1)   -  \Ga 1 2  1 \Ga 2  1  1   -  \Ga 1 2  2 \Ga 2  2  1   +    \Ga 2 2 1  \Ga 1 1  1+     \Ga 2 2 2  \Ga 1 2  1 + \frac{\p}{\p x}(\dGa 2 2 1)   -  \frac{\p}{\p y}(\dGa 1 2 1) -\\
            &  -  \dGa 1 2  1 \Ga 2  1  1   -  \dGa 1 2  2 \Ga 2  2  1   +    \dGa 2 2 1  \Ga 1 1  1+     \dGa 2 2 2  \Ga 1 2  1 -  \Ga 1 2  1 \dGa 2  1  1   -  \Ga 1 2  2 \dGa 2  2  1   +    \Ga 2 2 1  \dGa 1 1  1+     \Ga 2 2 2  \dGa 1 2  1 - \\
              &-  \dGa 1 2  1 \dGa 2  1  1   -  \dGa 1 2  2 \dGa 2  2  1   +    \dGa 2 2 1  \dGa 1 1  1+     \dGa 2 2 2  \dGa 1 2  1  + \lambda   \GGa {\qs_o} 2 1  + B \frac{\p}{\p r} \left(\GGa 2 2 {\qs_o} + \f_y   \GGa 2 2 1 +  \f_x \GGa 1 2 1\right) \ .
        \end{align*}
       }
       }
       \ \\
 In both expansions, the first six summands give the components $\Ric^N_{11}$  and  $\Ric^N_{22}$  of the Ricci tensor of the  K\"ahler metric $g_o$ of $(N, J)$.  Hence, 
  using the  expressions for the $\Ga i j k$ and $\dGa i j k$  given in \eqref{6.78} and the fact that the functions $\f$  and $\Ga i j k$ are independent of the coordinate  $r$, the  equations $\Ric_{ii} = 0$, $ i = 1,2$, become \par
    \scalebox{0.8}{\vbox{
  \begin{align}
 \nonumber  \Ric_{11}  & = \Ric^N_{11}    -  \frac{\p}{\p y}(\dGa 2 2 2)  -  \frac{\p}{\p x}(\dGa 1 1 1)    +\xcancel{\Ga 2 2  2\dGa 2  2  2}    -  \xcancel{\Ga 1 1  1 \dGa 1  1  1}  +   \xcancel{\Ga 1 1 1  \dGa 1 1  1 } - \xcancel{  \Ga 2 2 2   \dGa 2 2  2} + \\
  \nonumber      & + \xcancel{ \dGa 2 2  2  \Ga 2  2  2}  -\xcancel{ \dGa 1 1  1 \Ga 1  1  1 }   + \xcancel{\dGa 1 1 1   \Ga 1 1  1} -\xcancel{  \dGa 2 2  2  \Ga 2 2  2} +  \xcancel{\dGa 2 2  2\dGa 2 2  2}    -\xcancel{  \dGa 1 1  1 \dGa 1  1  1}   +  \xcancel{ \dGa 1 1 1  \dGa 1 1  1}  -   \xcancel{\dGa 2 2 2   \dGa 2 2  2} - \\
 \label{6.68} &-  \lambda   \GGa {\qs_o} 1 2  + B  \frac{\p}{\p r} \left(\GGa 1 1 {\qs_o} -  \f_y \GGa 1 1 1    +  \f_x \GGa 2 2 2   \right)  
  \end{align}
   \begin{align}
 \nonumber   \Ric_{22}  & = \Ric^N_{22}     - \frac{\p}{\p x}(\dGa 1 1 1)   -  \frac{\p}{\p y}(\dGa 2 2 2)      -  \xcancel{\dGa 2 2  2 \Ga 2  2  2}   + \xcancel{ \dGa 1 1  1 \Ga 1  1  1 }  -  \xcancel{ \dGa 1 1 1  \Ga 1 1  1}+    \xcancel{ \dGa 2 2 2  \Ga 2 2  2} -\\
  \nonumber      & -  \xcancel{\Ga 2 2  2 \dGa 2  2  2}   +\xcancel{  \Ga 1 1  1 \dGa 1  1  1}   -    \xcancel{\Ga 1 1 1  \dGa 1 1  1} +   \xcancel{  \Ga 2 2 2  \dGa 2 2  2} -  \xcancel{\dGa 2 2  2 \dGa 2  2  2}  + \xcancel{ \dGa 1 1  1 \dGa 1  1  1}   -   \xcancel{ \dGa 1 1 1  \dGa 1 1  1}+    \xcancel{ \dGa 2 2 2  \dGa 2 2  2} + \\
   \label{6.69} & + \lambda   \GGa {\qs_o} 2 1  + B \frac{\p}{\p r} \left(\GGa 2 2 {\qs_o} - \f_y   \GGa 1 1 1 +  \f_x \GGa 2 2 2\right)  
  \end{align}
  }}\\
 We now  recall that 
  $ \GGa {\qs_o} 2 1 = - \GGa {\qs_o} 1 2$,  $\GGa 2 2 {\qs_o} = \GGa 1 1 {\qs_o}$ and that  $\Ric^N_{22} = \Ric^N_{11}$ (because $g_o$ is K\"ahler). Inserting these identities into  \eqref{6.69} we get \eqref{6.68}, showing that  
the equations  \eqref{6.74} and \eqref{6.76} are equivalent.\par
 Let us now focus on the  equation \eqref{6.74} using   \eqref{6.68}.  Using the property   that $\f$ and $\lambda$ are independent of the coordinate $r$, the equation can be written as\par
     \scalebox{0.8}{\vbox{
    \begin{align*}
   & \Ric^N_{11}     - \frac{\p}{\p x}\left( \frac{\s_x}{2\s } + B \frac{ \s_t \f_y}{2 \s} \right)   -  \frac{\p}{\p y}\left(\frac{1}{2 \s} \s _y - B \frac{\s_t \f_x}{2 \s}  \right)     -  \lambda   B  \frac{1}{4 \s}  + \\
  & + B \frac{\p}{\p r} \left( - \lambda   \s_t  - \f_y  \left(\Ga 1 1 1
+ \frac{1}{2 \s }   \s_x    + B \frac{\s_t \f_y }{2  \s } \right)    +  \f_x\left(\Ga 2 2  2
+ \frac{1}{2 \s }  \s_y  - B \frac{\s_t \f_x }{2 \s }\right) \right)  = \\
   = &\ \Ric^N_{11}     -\frac{1}{2} (\log \s)_{xx}    - \frac{B}{2}  ( \log \s)_{tx} \f_y  - \xcancel{\frac{B}{2}  (\log \s)_t \f_{xy}}   -  \frac{1}{2} (\log \s)_{yy} + \frac{B}{2} (\log\s)_{ty} \f_x + \xcancel{ \frac{B}{2} (\log\s)_{t} \f_{xy}  }    - \\
   &\ - \frac{  \lambda B }{4 \s}  - B\lambda   \s_{tt}  -    \f_y   \frac{B}{2}   (\log \s)_{xt}    -      B^2 \frac{(\log \s)_{tt} (\f_y)^2 }{4 }   
+ \frac{B}{2 } (\log \s)_{ty} \f_x   - B^2 \frac{(\log \s)_{tt}( \f_x)^2}{2 }   =\\
= & \   \Ric^N_{11}     -\frac{1}{2} (\log \s)_{xx}  - \frac{1}{2} (\log \s)_{yy}  - B  ( \log \s)_{tx} \f_y    + B (\log\s)_{ty} \f_x  - \frac{  \lambda  B }{4 \s}  - B\lambda   \s_{tt}   -   \\
&\ -   B^2 \frac{(\log \s)_{tt} \left(( \f_x)^2 +(\f_y)^2 \right)}{4 }     = 0\ .
  \end{align*}
  }}\\
Since\par
    \scalebox{0.8}{\vbox{
  $$(\log \s)_x = \frac{\s_x}{\s} \ ,\quad (\log \s)_{xx} = \frac{\s_{xx}}{\s} - \frac{\s_x^2}{\s^2}\ ,\quad
 (\log \s)_y = \frac{\s_y}{\s} \ ,\quad (\log \s)_{yy} = \frac{\s_{yy}}{\s} - \frac{\s_y^2}{\s^2}\ ,$$
   $$(\log \s)_{tx} = \frac{\s_{tx}}{\s} - \frac{\s_t \s_x}{\s^2}\ ,\quad  (\log \s)_{ty} = \frac{\s_{ty}}{\s} - \frac{\s_t \s_y}{\s^2}\ , \quad  (\log \s)_{tt} = \frac{\s_{tt}}{\s} - \frac{\s^2_t }{\s^2}\ ,$$
   }}\\ 
   the same equation can be also written as \par
       \scalebox{0.8}{\vbox{
     \begin{align*}
&   \Ric^N_{11}     - \frac{\s_{xx}}{2\s} + \frac{\s_x^2}{2\s^2}  -  \frac{\s_{yy}}{2\s} + \frac{\s_y^2}{2\s^2}  - B \frac{\s_{tx}}{\s}  \f_y + B \frac{\s_t \s_x}{\s^2} \f_y     + B \frac{\s_{ty}}{\s} \f_x - B \frac{\s_t \s_y}{\s^2} \f_x   -\\
  &\  - B \frac{ \lambda     }{4 \s}  - B \lambda   \s_{tt}   -  B^2  \left( \frac{\s_{tt}}{2\s} - \frac{\s^2_t }{2\s^2}\right) \left((  \f_x)^2 +( \f_y)^2 \right)       = 0\ .
  \end{align*}
  }}\\
 Note that   $\frac{\p}{\p r} =   \ps_o$. Hence,  re-arranging the summands in an appropriate way,  we  may write this equation as 
      \begin{align}
 \nonumber     &\Ric^N_{11}        -\frac{\s_{xx}}{2\s } +\frac{\s_x^2}{2\s^2}       -\frac{\s _{yy}}{2 \s}+\frac{\s_y^2}{2\s^2}   
 - B  \frac{ \lambda}{4 \s}      - B \lambda \ps_o(\ps_o(\s))-  \\
\nonumber   &- B \left(  \frac{E_1(\ps_o(\s))}{\s}-\frac{\ps_o(\s) E_1(\s)}{\s^2}  + B \frac{\ps_o(\ps_o(\s))}{2\s} \f_y- B \frac{(\ps_o(\s))^2}{2\s^2} \f_y\right) \f_y+\\
\label{77}  & + B \left( \frac{E_2(\ps_o(\s))}{\s}-  \frac{\ps_o(\s)E_2( \s)}{\s^2}   - B \frac{\ps_o(\ps_o(\s))}{2\s} \f_x + B \frac{(\ps_o(\s))^2}{2\s^2} \f_x\right)  \f_x   = 0 \ .
  \end{align}
  We now recall  $n = 4$ and $\f$ and $\s$ satisfy  \eqref{5.32ter}, which is  the same as the system\par
    \scalebox{0.8}{\vbox{
  \begin{align*}
  &  -\frac{1}{2 \s^2} E_2 \left(\s\right)       -  \frac{1}{\s} E_1 \left( \ps_o(\s) \right)    +  \frac{1}{\s^2} \ps_o(\s) E_1 \left(\s\right)     +   \frac{B}{2 \s^2} \f_x \ps_o \left(\s\right)      -   \frac{B}{\s}     \f_y \ps_o \left(\ps_o(\s) \right)+   \frac{B}{\s^2}      \f_y \left( \ps_o \left(\s\right) \right)^2 
 = 0\ ,
\end{align*}
 \begin{align*}
  &    \frac{1}{2 \s^2}   E_1 \left(\s\right)       -  \frac{1}{\s}  E_2 \left( \ps_o(\s) \right)    +  \frac{1}{\s^2} \ps_o(\s) E_2 \left(\s\right)    + \frac{B}{ 2 \s^2}    \f_y \ps_o \left(\s\right)      +  \frac{B}{\s}     \f_x \ps_o \left(\ps_o(\s) \right) - \frac{B}{\s^2}     \f_x  \left( \ps_o \left(\s\right) \right)^2 
 = 0\ ,
\end{align*}
}}\\
or, equivalently, to \par
    \scalebox{0.8}{\vbox{
   \begin{align*}
   &   \frac{ E_1 \left( \ps_o(\s) \right) }{ \s}    -  \frac{ \ps_o(\s) E_1 \left(\s\right)  }{ \s^2}   =  -\frac{1}{2 \s^2} E_2 \left(\s\right)    + \frac{ B}{2 \s^2}   \f_x \ps_o \left(\s\right)      -      \frac{ B}{ \s}   \f_y \ps_o \left(\ps_o(\s) \right)+      \frac{B}{ \s^2}   \f_y\left( \ps_o \left(\s\right) \right)^2 \ ,
 \end{align*}
 \begin{align*}
 &    \frac{ E_2 \left( \ps_o(\s) \right)}{ \s}     -  \frac{ \ps_o(\s) E_2 \left(\s\right) }{\s^2}   =   \frac{1}{2 \s^2} E_1 \left(\s\right)    +  \frac{ B}{2 \s^2}   \f_y \ps_o \left(\s\right)      +     \frac{B }{\s}    \f_x \ps_o \left(\ps_o(\s) \right) -      \frac{ B}{ \s^2}    \f_x \left( \ps_o \left(\s\right) \right)^2 \ .
 \end{align*}
 }}\\
 Replacing this in \eqref{77},  we get\par
     \scalebox{0.8}{\vbox{
    \begin{align}
 \nonumber     &\Ric^N_{11}        -\frac{\s_{xx}}{2\s } +\frac{\s_x^2}{2\s^2}       -\frac{\s _{yy}}{2 \s}+\frac{\s_y^2}{2\s^2}   
  - B  \frac{ \lambda}{4 \s}     - B \lambda \ps_o(\ps_o(\s))-  \\
\nonumber   &- B \left(   -\frac{1}{2 \s^2} \s_y    +\xcancel{ \frac{ B}{2 \s^2}   \f_x \ps_o \left(\s\right)}      -     \frac{ B}{ 2\s}   \f_y \ps_o \left(\ps_o(\s) \right) +       \frac{B}{ 2\s^2}   \f_y\left( \ps_o \left(\s\right) \right)^2 
\right) \f_y+\\
\label{77*}  & + B \left(   \frac{1}{2 \s^2}\s_x   +  \xcancel{\frac{ B}{2 \s^2}   \f_y \ps_o \left(\s\right)}      +       \frac{B }{2\s}    \f_x \ps_o \left(\ps_o(\s) \right) -       \frac{ B}{ 2\s^2}    \f_x \left( \ps_o \left(\s\right) \right)^2   \right)  \f_x   = 0 \ .
  \end{align}
  }}\\
  Since
  $-2\s\ps_o(\ps_o(\s)) +\ps_o(\s)^2+ \frac{1 }{4} = 0$ (which  allows us to replace  $  -\frac{B}{2\s^2}\ps_o(\s)^2 = -  \frac{B}{\s} \ps_o(\ps_o(\s)) +    \frac{B}{8\s^2 } $ in all terms), multiplying left and right hand side by $\s^2$, 
the equation transforms into \eqref{commonequation}, as we needed to prove.

\section*{Declarations}
The authors declare that they have no conflict of interest.

\vskip 1.5truecm
\hbox{\parindent=0pt\parskip=0pt
\vbox{\baselineskip 11 pt \hsize=3.1truein
\obeylines
{\smallsmc
Masoud Ganji and Gerd Schmalz
School of Science and Technology
University of New England,
Armidale NSW 2351
Australia
}\medskip
{{\smallit E-mail}\/: {\smalltt mganjia2@une.edu.au}}
{\smallit E-mail}\/: {\smalltt schmalz@une.edu.au}
}
\vbox{\baselineskip 11 pt \hsize=3.1truein
\obeylines
{\smallsmc
Cristina Giannotti and Andrea Spiro
Scuola di Scienze e TecnolSogie
Universit\`a di Camerino
I-62032 Camerino (Macerata)
Italy
}\medskip
{{\smallit E-mail}\/: {\smalltt cristina.giannotti@unicam.it}
}
{\smallit E-mail}\/: {\smalltt andrea.spiro@unicam.it
}
}
}
\end{document}